\documentclass[12pt,a4paper]{article}

\usepackage{a4wide,amssymb,amsmath,amsthm,xspace,epsfig, amsfonts,amscd}
\usepackage{tikz-cd}

\usepackage[left=1.5cm, right=1.5cm, top=1.5cm]{geometry} 

\usepackage{multirow}

\usepackage{hyperref}

\usepackage{mathrsfs}

\usepackage{float}

\begin{document}

\def\checkmark{\tikz\fill[scale=0.4](0,.35) -- (.25,0) -- (.8,.7) -- (.25,.15) -- cycle;}

\newcommand{\N}{\mathbb{N}}
\newcommand{\R}{\mathbb{R}}
\newcommand{\Z}{\mathbb{Z}}
\newcommand{\Q}{\mathbb{Q}}
\newcommand{\C}{\mathbb{C}}
\newcommand{\PP}{\mathbb{P}}

\newcommand{\LL}{\Bbb L}
\newcommand{\OO}{\mathcal{O}}

\newcommand{\esp}{\vskip .3cm \noindent}
\mathchardef\flat="115B

\newcommand{\lev}{\text{\rm Lev}}

\def\ut#1{$\underline{\text{#1}}$}
\def\CC#1{${\cal C}^{#1}$}
\def\h#1{\hat #1}
\def\t#1{\tilde #1}
\def\wt#1{\widetilde{#1}}
\def\wh#1{\widehat{#1}}
\def\wb#1{\overline{#1}}

\def\restrict#1{\bigr|_{#1}}

\def\hu#1#2{\mathsf{U}_{fin}\bigl({#1},{#2}\bigr)}
\def\ch#1#2{\left(\begin{array}{c}#1 \\ #2 \end{array}\right)}

\def\EC#1#2{\mathsf{EC}\bigl({#1},{#2}\bigr)}
\def\S#1#2{\mathsf{S}\bigl({#1},{#2}\bigr)}
\def\BR#1#2{\mathsf{BR}\bigl({#1},{#2}\bigr)}
\def\L#1#2{\mathsf{L}\bigl({#1},{#2}\bigr)}

\newtheorem{lemma}{Lemma}[section]

\newtheorem{thm}[lemma]{Theorem}
\newtheorem*{thm*}{Theorem}
\newtheorem*{lemma?}{Lemma ??}

\newtheorem{defi}[lemma]{Definition}
\newtheorem{conj}[lemma]{Conjecture}
\newtheorem{cor}[lemma]{Corollary}
\newtheorem{prop}[lemma]{Proposition}
\newtheorem{prob}[lemma]{Problem}
\newtheorem{qu}[lemma]{Question}
\newtheorem{q}[lemma]{Question}
\newtheorem*{rem}{Remark}
\newtheorem{rem_numbered}[lemma]{Remark}
\newtheorem{examples}[lemma]{Examples}
\newtheorem{example}[lemma]{Example}

\title{Eventually Constant and stagnating functions in non-Lindel\"of spaces}
\author{Mathieu Baillif}
\maketitle

\abstract{\footnotesize Inspired by recent work of A. Mardani
          which elaborates on the elementary fact that for any continuous function $f:\omega_1\times\R\to\R$,
          there is an $\alpha\in\omega_1$ such that 
          $f(\langle\beta,x\rangle) = f(\langle\alpha,x\rangle)$ for all $\beta\ge\alpha$ and $x\in\R$,
          we introduce four properties
          $\mathsf{P}(X,Y)$, $\mathsf{P}\in\{\mathsf{EC},\mathsf{S},\mathsf{L},\mathsf{BR}\}$, 
          which are different formalizations of the idea vaguely
          stated as
          ``given a continuous $f:X\to Y$, there is a small subspace of $X$ outside of 
          which $f$ does not do anything much new''.
          More precisely, we say that the spaces $X,Y$ satisfy the property 
          $\mathsf{EC}(X,Y)$ [resp. $\mathsf{S}(X,Y)$]
          (resp. $\mathsf{L}(X,Y)$)
          iff given $f:X\to Y$, then there is a Lindel\"of $Z\subset X$ such that $f(X-Z)$ is a singleton [resp.
          there is a retraction $r:X\to Z$ such that $f\circ r = f$] (resp. $f(Z) = f(X)$). 
          $\mathsf{BR}(X,Y)$ is defined similarly. 
          We investigate the relations between these four generalizations of Lindel\"ofness
          and other classical topological properties.
          Actually, two more variants $\mathsf{P_{cl}},\mathsf{P_{cpt}}$ of each 
          property are given depending on whether $Z$ can be chosen to be closed or compact.
          To get an idea of what our results look like, here is a sample.\\
          An uncountable subspace $T$ of a tree of height $\omega_1$ is 
          $\omega_1$-compact iff $\mathsf{S}(T,Y)$ holds for any metrizable space $Y$ of uncountable cardinality. 
          (The case $Y=\R$ and $T$ a Suslin tree was proved by Stepr\=ans long ago.) 
          If $M$ is a $\aleph_1$-strongly 
          collectionwise Hausforff non-metrizable manifold satisfying either
          (a weakening of) $\mathsf{S}(M,\R)$ or $\mathsf{EC}(M,\R)$, then $M$
          is $\omega_1$-compact. 
          The property $\mathsf{L}(M,\R)$ holds for any manifold while $\mathsf{L}(M,\R^2)$ does not.
          Under {\bf PFA}, a locally compact countably tight space $Y$ for which 
          $\mathsf{EC}(\omega_1,Y)$ holds is isocompact, while
          there are counterexamples under $\clubsuit_C$.
          There is a pseudocompact space $X$ such that $\mathsf{L}(X,\R)$ does not hold.
          Some of our results are (more or less elaborate) restatements of other researchers work 
          put in our context. 
          }

\tableofcontents

\section{Introduction}
In this note, by `space' we mean 
topological Hausdorff space (in particular `regular' and 'normal' imply Hausdorff) and every function is assumed continuous, unless specified. 
We denote ordered pairs with brackets $\langle \cdot,\cdot\rangle$, 
reserving parenthesis for intervals in (totally) ordered spaces.
We refer to \cite{Engelking} for standard topological notions\footnote{Notice however
that we sometimes depart from the conventions of \cite{Engelking} about including or not separation axioms into properties
like pseudocompactness. We will specify when this is the case.}
and to the end of this introduction for less standard ones. 
Our starting point is the following well known elementary result.
(For a proof, combine Lemma B.30 in \cite{GauldBook} with the fact that $\R$ is separable.)  
As usual $\omega_1$ is the first uncountable ordinal endowed with the order topology.

\begin{thm}\label{thm0}
   Let $f:\omega_1\times\R\to\R$ be continuous.
   Then there is $\alpha\in\omega_1$ such that 
   $f(\langle\beta,x\rangle) = f(\langle\alpha,x\rangle)$ for all $\beta\ge\alpha$ and $x\in\R$.
\end{thm}

To summarize in a
very imprecise catchphrase: {\em Outside of a small subspace, $f$ does not do anything much 
new.
}
Here of course the small subspace is $[0,\alpha]\times\R$.
The purpose of this note is to investigate
the notions obtained when `small' and `nothing much new' are interpreted in various ways.
Well, not so various as to the former:
`Small' will almost always mean `Lindel\"of' (the only real exception being Section \ref{sec:cpct}).
This might be surprising at first since compactness is arguably a more natural mesure of smallness of a subspace;
but
we believe that our choice yields more interesting results, as replacing it by compactness 
almost always
confines us to pseudocompact spaces, at least for real valued maps. 
\\
The phrase
`$f$ does not do anything much new' gives way to more varied interpretations (at least to the extent permitted by our
imagination).
A rather general interpretation is that there is a small subset whose image is equal to that of the whole domain space.
In Theorem \ref{thm0}, $f([0,\alpha]\times\R) = f(\omega_1\times\R)$.
Figuratively speaking, $f$ (or $\omega_1\times\R$) is a lazy explorer of $\R$:
after a small exploration it remains at the same places 
forever.\footnote{We could have invoked the figure of a regular customer in a western european local bar, since these individuals
tend to keep to the same place as well after choosing their favorite table and chair. We 
refrained from it in view of the confusion induced by the terminology collision with regular spaces, 
and by broken-record-like conversations which tend to emerge in these settings, which would mistakenly allude 
to the property defined just next.
}
Looking at it from the other side, 
we can say that outside of a small subset, any value taken by $f$ will be taken again and again, like a broken record, 
as we go further away (horizontally in this case) in the domain 
space.
In Theorem \ref{thm0}, if $\beta>\alpha$, there is $\gamma>\beta$ such that 
$f(\langle\gamma,t\rangle)=f(\langle\beta,t\rangle)$
(of course, this actually holds for each $\gamma>\beta$ in that case).
Another more restrictive (a priori)
interpretation is that
the map $f$ {\em stagnates} outside of a small subspace;
that is, there is a retraction $r$ of $X$ onto a small subset such that $f\circ r =f$. 
A {\em retraction} is a map $r:X\to Z\subset X$ whose restriction to $Z$ is the identity. 
In Theorem \ref{thm0} we may define it as 
$r(\langle\beta,t\rangle) = \langle\min\{\alpha,\beta\},t\rangle$.
Lastly,
looking only at horizontal lines, we see that any $g:\omega_1\to\R$ is 
{\em eventually constant}, in the sense that $g$ is constant outside of
the small subset
$[0,\alpha]$. 
These four interpretations are formalized in the following 
definition.

\begin{defi} \label{defgen}
   Let $X, Y$ be spaces. 
   Then $X,Y$ satisfy the left side of the table below iff 
   for each $f:X\to Y$, there is a Lindel\"of subset $Z\subset X$ such that the right side holds.
   The middle boxes contain the shorthands for each property.
   \begin{center}
   \begin{tabular}{|l|c|l|}
   \hline
   $X$ is eventually constant in $Y$ & $\mathsf{EC}(X,Y)$ &  $f(X-Z)$ is a singleton \\
   \hline
   \multirow{2}{*}{$X$ stagnates in $Y$} & \multirow{2}{*}{$\mathsf{S}(X,Y)$} & there is a retraction \\
   & & $r:X\to Z$ satisfying $f= f\circ r $ \\
   \hline
   $X$ is a lazy explorer of $Y$ & $\mathsf{L}(X,Y)$  & $f(Z) = f(X)$ \\
   \hline
   \multirow{2}{*}{$X$ is a broken record in $Y$} & \multirow{2}{*}{$\mathsf{BR}(X,Y)$}  & for each Lindel\"of $W\supset Z$, \\
   & & $f(X-W) = f(X-Z)$\\
   \hline
   \end{tabular}
   \end{center}
   Given $\mathsf{P}\in\{\mathsf{EC},\mathsf{S},\mathsf{L},\mathsf{BR}\}$,
   if $Z$ can be chosen to be closed, we write $\mathsf{P_{cl}}(X,Y)$ for the stronger property,
   and we write $\mathsf{P_{cpt}}(X,Y)$
   if $Z$ can be chosen to be compact. 
   We abbreviate
   $\mathsf{P}(X,\R)$, $\mathsf{P_{cl}}(X,\R)$, $\mathsf{P_{cpt}}(X,\R)$
   by $\mathsf{P}(X)$, $\mathsf{P_{cl}}(X)$, $\mathsf{P_{cpt}}(X)$.
\end{defi}
(We already stress that there are spaces $X$ for which, really, nothing much new happens outside of some closed 
Lindel\"of subset for any real valued map, but
none of these properties with $Y=\R$ hold, see Example \ref{ex:pasAS} (c).)
To help concision, the letter $\mathsf{P}$ denotes any of the properties $\mathsf{EC},\mathsf{S},\mathsf{L},\mathsf{BR}$
in this paragraph.
Since our spaces are Hausdorff, we immediately have 
$\mathsf{P_{cpt}}(X,Y)\Longrightarrow\mathsf{P_{cl}}(X,Y)\Longrightarrow\mathsf{P}(X,Y)$
for all spaces $X,Y$.
Notice in passing that any space $X$ with the property that any real valued fonction on $X$ is constant satisfies 
$\mathsf{P_{cpt}}(X)$ for each $\mathsf{P}$.
Such spaces do exist, some being even regular
(for instance \cite[Example 92]{CEIT}).
We also notice that 
every functionally countable space $X$ (i.e. a space whose image under a real valued continuous map is countable)
satisfies $\mathsf{L}(X)$. More generally, compact spaces trivially satisfy each
$\mathsf{P_{cpt}}(X,Y)$ and Lindel\"of spaces each $\mathsf{P}(X,Y)$ \& $\mathsf{P_{cl}}(X,Y)$ for any $Y$.
Hence, seen as properties of $X$, they are a generalizations of compactness and
Lindel\"ofness.

To our knowledge, all the properties in Definition \ref{defgen} are new, at least at this level of generality.
Some particular instances have been already investigated before.
For example,
spaces $Y$ satisfying $\mathsf{EC}(\omega_1,Y)$ are called $\omega_1$-squat by A. Mardani, whose PhD thesis 
\cite{MardaniThesis}
contains results about this and related classes of spaces which motivated
this work. The term `squat' was first used by D. Gauld (in particular in \cite{GauldBook}) for almost the same property.
Many examples of spaces $X$ satisfying $\mathsf{EC}(X)$ appear in the litterature, we shall enconter some
of them later on. The author became interested in these properties while investigating which non-metrizable manifolds
could be contractible. It happens that $X$ cannot be contractible
when 
real-valued maps have features similar to 
$\mathsf{EC}(X)$ and $\mathsf{S}(X)$  (see \cite{mesziguesContractibility}).  

\vskip .3cm
\noindent
This note contains our musings about the interplays these twelve 
notions have with each other and other topological properties
in various classes of spaces.
It is organized in sections of variable lengths whose titles hint to their contents.
We were driven by pure curiosity and the pleasure of wandering in this 
landscape\footnote{We believe that it is important for the reader to keep that in mind.
Despite the length of this paper,
we have a leisurely approach and do not claim any grand result.}. 
In some cases we have tried to obtain quite general results 
(for instance in sections \ref{sec:interplays}--\ref{sec:isocompact}), 
while in others 
we have concentrated our efforts on particular classes of spaces, such as set-theoretic trees 
(the entire section \ref{sec:Suslin}) and non-metrizable manifolds 
(all of sections \ref{sec:lackretr}--\ref{sec:SscwH} and some of section \ref{sec:BRL}).
As such, this note is more akin to a little stroll in the garden -- rusty old tools in hands,
peeking below scattered rocks and looking for strange insects -- 
than to securing the foundations for a (future) twenty storeys tower. 
Some of our results are (more or less elaborate) restatements of other researchers work 
put in our context, and we shall also re-encounter many classical examples of general topology.
It is quite possible that 
effective shortcuts escaped us while we were meandering.\footnote{Since our methods 
are mainly elementary, it is also possible that some of our results
already appeared in the litterature, disguised in different terminology.}
We do not know if our results will seem appealing to other researchers who tend to enjoy the scenery differently from us,
but we 
hope that the reader's mood is bucolic enough to at least enjoy some of them.
The sections are somewhat independent, although the later ones tend to refer to the previous ones 
(big surprise). 
There are three exceptions:
Section \ref{sec:interplays}, which contains 
basic results used pervasively, Section \ref{sec:tables}, which contains tables that briefly
summarize the properties of the strange insects we shall encounter, 
and Section \ref{sec:questions} where we compile the questions that appear along the text
and ask other ones about the relashionships between our properties and
recent generalizations of Lindel\"ofness and compactness involving cellular families.

We note that despite their similarities, 
$\mathsf{EC}$, $\mathsf{L}$, $\mathsf{S}$ and $\mathsf{BR}$
seem to lead to different techniques
of investigation.
For instance, $\mathsf{EC}$ seems to be linked to countable compactness, $\omega_1$-compactness and isocompactness
(see Section \ref{sec:isocompact} for definitions), that is: general topology considerations.
In the case of trees, $\mathsf{S}$ is also linked to $\omega_1$-compactness, but 
we borrow from more geometric techniques (even very basic covering spaces theory) when looking at 
$\mathsf{S}(M)$ for manifolds $M$. Finally, 
$\mathsf{L}$ and $\mathsf{BR}$ seem to be well adapted to spaces that are increasing
chains of ``nices subspaces''. These differences end up scattering the contents
of this article, in a way quite similar, 
we recall again, to the impressions imbued by a walk in a not very well kept garden.

We end this introduction by recalling definitions that are either 
not-so-standard or slightly different from the usual ones.

\vskip .3cm
By {\em cover} of a space $X$ is understood a family of open sets whose union contains $X$.
A cover is a {\em chain cover} if it is linearly ordered by the inclusion relation.
Any non-trivial chain cover (that is, one without the whole space as a maximal element)
has a subcover indexed by a regular cardinal whose members are pairwise distinct. 
We always implicitely take such a subcover and use such indexing.
A space $X$ is of {\em Type I} (Nyikos \cite{Nyikos:1984}) 
iff $X=\cup_{\alpha<\omega_1}X_\alpha$, where $X_\alpha$ is open, $\wb{X_\alpha}\subset X_\beta$ whenever $\alpha<\beta$,
and $\wb{X_\alpha}$ is Lindel\"of for each $\alpha$. 
Although it it not included in the usual definition, for simplicity we assume in this note that Type I spaces
are {\em not} Lindel\"of, that is, $X\not= X_\alpha$ for each $\alpha$.
If $\cup_{\alpha <\beta} X_\alpha = X_\beta$ whenever $\beta$ is limit, the cover 
$\{X_\alpha\,:\,\alpha<\omega_1 \}$ is called {\em canonical}. Any chain cover 
of a Type I space has a subcover indexed by $\omega_1$ which can be made canonical by adding the missing $X_\alpha$'s.
Any two canonical covers agree on a {\em club} (i.e.
closed and unbounded) subset of $\omega_1$, as easily seen. Hence,
if $X$ is a Type I space, $X_\alpha$ will always denote the $\alpha$th member of some (often implicit) canonical cover.
We borrow the usual vocabulary used in $\omega_1$ for Type I spaces: a subset of $X$ is {\em bounded} 
iff it is contained in some $X_\alpha$
and {\em unbounded} otherwise, and
{\em club} means closed and unbounded. 
The closure of a Lindel\"of subset is Lindel\"of in a Type I space.

We denote the cardinality of the continuum by $\mathfrak{c}$.
By {\em manifold} or {\em $n$-manifold} (when we want to emphasize the dimension) 
we mean a connected space locally homeomorphic to $\R^n$. A {\em surface} is a $2$-manifold.
Connectedness is often an indispensable property in our results about manifolds.
{\em Manifolds with boundary} have some points 
(those in the manifold boundary) with open neighborhoods homeomorphic to $\R_{\ge 0}\times\R^{n-1}$ 
(and no open neighborhood homeomorphic to $\R^n$).
(There is an unfortunate terminology collision with the topological boundary, that is, the closure minus
the interior of a subset. When confusion is possible, we shall specify which meaning is intended.)
The longray $\LL_{\ge 0}$ is the $1$-manifold with boundary 
$\omega_1\times[0,1)$ with lexicographic order topology ({\em not} the product topology).
We often view $\omega_1$ as a subset of $\LL_{\ge 0}$ by identifying $\alpha$ with $\langle \alpha,0\rangle$,
and hence write for instance $\LL_{\ge_0} - \omega_1$ instead of $\LL_{\ge_0} - \omega_1\times\{0\}$.
The open longray $\LL_+$ is obtained by deleting the $0$ point in $\LL_{\ge 0}$.
A space is {\em $\omega$-bounded} iff any countable subset has a compact closure. For Type I spaces this 
is equivalent to being countably compact.
A manifold is $\omega$-bounded iff it is countably compact and of Type I \cite[Thm 4.10]{GauldBook}.
A {\em longpipe} \cite[Def. 4.11]{GauldBook} 
is an $\omega$-bounded (hence Type I) surface with boundary $S=\cup_{\alpha\in \omega_1}S_\alpha$ such that
$S_{\alpha+1}$ is homeomorphic to the cylinder $\mathbb{S}^1\times[0,1)$ and the topological boundary 
of $S_{\alpha+1}$ in $S_{\beta}$ is homeomorphic to the circle for each $\alpha$ and $\beta\ge\alpha+2$. 
(This may not be true if $\alpha$ is a limit ordinal, though.)
For more on non-metrizable manifolds and longpipes, see \cite{GauldBook} and \cite{Nyikos:1984}.

Finally, 
recall that a space is {\em $\omega_1$-compact} or has {\em countable extent}
iff its closed discrete subspaces are at most countable.
The following classical facts (see e.g. \cite[Thm 4.1.15, Thm 4.1.17]{Engelking}) will be
useful in some proofs.

\begin{lemma}\label{lemma:metriceq}
   In a metrizable space, the properties $\omega_1$-compact, Lindel\"of, 
   hereditarily Lindel\"of, separable and hereditarily separable 
   are all equivalent, and compactness is equivalent to countable compactness.
\end{lemma}

We close this introduction with the following elementary remark, which will be used implicitely many times.
\begin{rem_numbered}
   \label{rem:elementary}
   In the definitions of properties $\mathsf{EC},\mathsf{L},\mathsf{BR}$ in \ref{defgen},
   if the right side holds for the Lindel\"of subset $Z$ then it holds for any Lindel\"of subset containing $Z$.
\end{rem_numbered}


\section{Does a lazy constant broken record stagnate~? 
         }
\label{sec:interplays}

An astute reader probably suspects that 
the answer to the question in this section's title is {\em absolutely yes, but sometimes no, and reciprocally}.
This section is devoted to giving more details about this answer (and the question itself)
for properties $\mathsf{P}$ and $\mathsf{P_{cl}}$ (those of the form $\mathsf{P_{cpt}}$ are dealt with
in Section \ref{sec:cpct}).
When there is no risk of confusion, we abbreviate the statement 
`for all Hausdorff spaces $X,Y$, $\mathsf{P}_1(X,Y)\Rightarrow\mathsf{P}_2(X,Y)$'
by `$\mathsf{P}_1\Rightarrow\mathsf{P}_2$'. `$\mathsf{P}_1\not\Rightarrow\mathsf{P}_2$' and
`$\mathsf{P}_1\Leftrightarrow\mathsf{P}_2$' are to be understood similarly.
Notice that for $\mathsf{S}(X,Y)$,
$Z$ has to be closed in Definition \ref{defgen} since it is the set of fixed points of the retraction and the spaces are Hausdorff, 
so $\mathsf{S}\Leftrightarrow\mathsf{S_{cl}}$.
The others implications in the lemma below are immediate from the definitions.
(For aesthetical reasons, we tend to denote implication 
by single arrows $\longrightarrow$ in somewhat complicated diagrams, and
by double arrows $\Longrightarrow$ in one-line formulas.)
\begin{lemma}\label{lemma:arrows}
  The following implications hold.
  \begin{center}
  \begin{tikzcd}
  \mathsf{EC_{cl}}
    \arrow[rrr, bend left, rightarrow]
    \arrow[r, rightarrow]
    \arrow[d, rightarrow]
   & 
   \mathsf{BR_{cl}}
     \arrow[r, rightarrow]
   &
   \mathsf{BR}
     \arrow[r, leftarrow]
   &
   \mathsf{EC}
     \arrow[d, rightarrow]
   \\
   \mathsf{L_{cl}}
    \arrow[rrr, bend right, rightarrow]
    \arrow[r, leftarrow]
   &
   \mathsf{S_{cl}}
     \arrow[r, leftrightarrow]
   &
   \mathsf{S}
     \arrow[r, rightarrow]
   & 
   \mathsf{L}
   \end{tikzcd}
   \end{center}
  \end{lemma} 
Implications not shown do not hold for all spaces $X$, even when $Y=\R$. 
In particular, none of $\mathsf{L}$, $\mathsf{EC}$, $\mathsf{BR}$
imply their $\mathsf{cl}$-counterpart in general.
It is however the case for Type I spaces.

\begin{lemma}\label{lemma:typeIcl}
    For Type I spaces $X$, $\mathsf{P}(X,Y)\Leftrightarrow \mathsf{P_{cl}}(X,Y)$ for any space $Y$ and 
   $\mathsf{P}\in\{\mathsf{EC},\mathsf{L},\mathsf{BR},\mathsf{S}\}$.
\end{lemma}
\proof
   For $\mathsf{S}$, this always holds regardless of whether $X$ is Type I.
   If $X$ is Type I, 
   a Lindel\"of subspace of $X$ is contained in $\wb{X_\alpha}$ (Lindel\"of itself) for some $\alpha$, which yields
   the result by Remark \ref{rem:elementary}.
\endproof

Let us now have a look at arrows not in Lemma \ref{lemma:arrows}.
Our first example is trivial and shows that none of 
$\mathsf{S_{cl}}(X),\mathsf{L_{cl}}(X)$ and $\mathsf{BR_{cl}}(X)$
imply $\mathsf{EC}(X)$.
\begin{example}\label{thm:omega_1double}
  Let $X$ be the disjoint union of two copies of $\omega_1$. 
  Then $\mathsf{S_{cl}}(X),\mathsf{L_{cl}}(X)$ and $\mathsf{BR_{cl}}(X)$ hold but 
  $\mathsf{EC}(X)$ does not.
\end{example}
\proof[Details] 
$\mathsf{BR_{cl}}(X),\mathsf{S_{cl}}(X)$ and thus $\mathsf{L_{cl}}(X)$ hold by by 
Theorem \ref{thm0}. 
Taking a map sending one copy of $\omega_1$ on $0$ and the other on $1$ shows that $\mathsf{EC}(X)$ does not hold.
\endproof
Other examples as \ref{thm:omega_1double} are the {\em long line} $\LL$ made of 
two copies of $\LL_{\ge 0}$ glued at their $0$-point,
and the space $\omega_1\times\R$ of Theorem \ref{thm0}.
We now show that $\mathsf{P}(X)\not\Rightarrow\mathsf{P_{cl}}(X)$ for 
$\mathsf{P}\in\{\mathsf{EC},\mathsf{L},\mathsf{BR}\}$.
The idea is actually quite simple and can be infered from the two next lemmas. The first one is immediate.


\begin{lemma}
   \label{lemma:obvious2}
   Let $\mathsf{P}\in\{\mathsf{EC},\mathsf{L}\}$ 
   and
   $X$, $Y$ be spaces such that $X=A\cup B$ where $A$ is Lindel\"of and $\mathsf{P}(B,Y)$ holds. 
   Then $\mathsf{P}(X,Y)$ holds. 
\end{lemma}

\begin{lemma}\label{lemma:exnotcl}
   Let $X$ be a space with a dense countable subspace of isolated points which we identify with the integers $\omega$.
   If $X-\omega$ is not Lindel\"of, then $\mathsf{L_{cl}}(X)$ and $\mathsf{BR_{cl}}(X)$ do not hold.
\end{lemma}
\proof
   Let $f:X\to\R$ be defined by $f(n)=1/n$ for $n\in\omega$ and $f(x) = 0$ for $x\not\in\omega$.
       Then $f$ is continuous. Any $Z\subset X$ such that $f(Z) = f(X)$ must contain all of $\omega$,
       hence its closure is all of $X$, which is not Lindel\"of. It follows that $\mathsf{L_{cl}}(X)$ does not hold.
       The same function shows that $\mathsf{BR_{cl}}(X)$ does not hold: given $Z\subset X$, if 
       $\omega\not\subset Z$, for any $n\in\omega-Z$ we have 
       $$f(n)\in f(X-Z) \not= f(X-(Z\cup\{n\}))\not\ni f(n).$$
       Hence, $Z\supset\omega$ and its closure must again be the entire space.
\endproof

Hence, spaces as in this lemma are good candidates to show that $\mathsf{P}\not\rightarrow\mathsf{P_{cl}}$.
To complete the task, the following result is another tool.
Recall that $\beta X$ is the \v Cech-Stone compactification of the space $X$.

\begin{thm}\label{rem:betaX}
    Let $X$ be Tychonoff and let $\beta X$ be its \v Cech-Stone compactification.
    \\
    (a)
    If $|\beta X - X| = 1$, then $\mathsf{EC}(X)$ holds. 
    \\
    (b) If $\beta X - X$ is at most countable and $X$ is locally compact, then $\mathsf{L}(X)$ holds.
\end{thm}
Recall that a {\em $0$-set} (resp. a {\em co-$0$-set}) in a space $X$ is a preimage of $\{0\}$ (resp. of $(0,1]$)
under a map $X\to[0,1]$.
\proof
    If $A,B$ are disjoint $0$-sets in $X$, then their closure is disjoint in $\beta X$.
    (All claimed properties of $\beta X$ in this proof can be found in \cite[Section 3.6]{Engelking}.)
    This implies 
    that $|\beta X - X| = 1$ is equivalent to the property that given two
    disjoint $0$-sets in $X$, then at least one is compact. We show in Theorem \ref{thmDir} below that 
    this implies $\mathsf{EC}(X)$, 
    which proves (a).\\
    Assume as in (b) that $\beta X - X$ is at most countable, with $X$ locally compact, and let $f:X\to\R$ be given.
    We may assume that the range of $f$ is contained in $[0,1]$. 
    Let $\beta f:\beta X\to [0,1]$ be the extension of $f$ to all of $\beta X$ (which always exists,
    see e.g. \cite[Theorem 3.6.1]{Engelking}). 
    By local compactness, $X$ is open in $\beta X$ (see e.g. \cite[Theorem 3.5.8]{Engelking}).
    Then $C = \beta f( \beta X - X )$ is compact and countable in $[0,1]$.
    If $E\subset [0,1]$ is closed and disjoint from $C$, then $f^{-1}(E)$ is compact in $X$,
    otherwise its closure in $\beta X$ intersects $\beta X-X$,
    and hence $E\cap C\not=\varnothing$. Since any open subset of $[0,1]$ is an $F_\sigma$,
    $Z_0 = f^{-1}([0,1]-C)$ is Lindel\"of. 
    Add to $Z_0$ one preimage of each $c\in C$ such that $f^{-1}(C)\cap X\not = \varnothing$,
    to obtain a Lindel\"of subset $Z$ with full image. 
\endproof

\begin{example}[{S. Mr\'owka, in effect}]\label{ex:notclbis}
  There is a locally compact first countable pseudocompact separable space $X$ such that:\\
  (a)  $\mathsf{EC}(X)$ holds and hence so do $\mathsf{L}(X)$ and $\mathsf{BR}(X)$; \\
  (b)  $\mathsf{L_{cl}}(X)$ does not hold and thus neither do $\mathsf{EC_{cl}}(X)$ and $\mathsf{S}(X)$;\\
  (c)  $\mathsf{BR_{cl}}(X)$ does not hold.
\end{example}

A space $X$ is {\em pseudocompact} iff every real valued map defined on $X$ has a bounded image.
We do {\em not} assume that pseudocompact spaces are Tychonoff, departing from the definition given 
in \cite[p. 208]{Engelking}.
Recall also that a family of sets is {\em almost disjoint} iff the intersection of any
two members is finite.
\proof[Details]
Recall that a $\psi$-space is the union of an open countable discrete space (which we may take to be $\omega$) and
an uncountable discrete subspace whose points are given by 
a maximal family $\mathcal{R}$ of almost disjoint subsets of $\omega$.  
A neighborhood basis of $A\in \mathcal{R}$ is given by $\{A\}\cup(A-F)$ where $F$ is finite.
All $\psi$-spaces are pseudocompact and locally compact.
S. Mr\'owka \cite[Theorem 3.11]{Mrowka:1977} shows how to construct a $\psi$-space $X$ whose \v Cech-Stone compactification
is the one point compactification. 
This implies that $\mathsf{EC}(X)$ holds by Theorem \ref{rem:betaX}.
Of course, $\mathsf{L_{cl}}(X)$ and $\mathsf{BR_{cl}}(X)$ do not hold by Lemma \ref{lemma:exnotcl}.
\endproof

Notice in passing that a disjoint union of two spaces as in Example \ref{ex:notclbis}
yields a space satisfying $\mathsf{L}(X)$ and $\mathsf{BR}(X)$ but neither their cl-versions nor $\mathsf{EC}(X)$.
Another way of obtaining $\mathsf{P}\not\Rightarrow\mathsf{P_{cl}}$ is the following.
For the definition of the uncountable cardinal
$\mathfrak{p}$ and more on the subject, see see e.g. \cite[Chapter 3]{vanDouwen:1984}. 

\begin{example}[{P. Nyikos}]\label{ex:notcl}
  There is a locally compact first countable separable normal 
  space $X$ such that 
  $\mathsf{EC}(X)$ holds but $\mathsf{L_{cl}}(X)$ and $\mathsf{BR_{cl}}(X)$ do not.
  If $\mathfrak{p}=\omega_1$, the space can be made to be moreover countably compact.
\end{example}
\proof[Details]
  If $X$ is as in Lemma \ref{lemma:exnotcl} but such that $X-\omega$ is homeomorphic to $\omega_1$ (in the order topology),
  then $\mathsf{EC}(X)$ holds (by Lemma \ref{lemma:obvious2}) while $\mathsf{L_{cl}}(X)$ and $\mathsf{BR_{cl}}(X)$ do not.
  There are many ways to obtain normal first countable and locally compact such spaces, 
  see for instance \cite{Nyikos:HereditaryNormality}.
  Nyikos (\cite{Nyikos:countablycompact1986}, Theorem 2.1 and Example 3.4) shows how to 
  obtain a countably compact example whenever $\mathfrak{p}=\omega_1$.
\endproof

A variation of this space 
can be adapted to obtain a surface, see Section \ref{sec:cpct}.
The next examples show in particular that 
$\mathsf{L}\not\Rightarrow\mathsf{BR}$ and $\mathsf{BR}\not\Rightarrow\mathsf{L}$.
It is convenient to first introduce the following notation.
If $f:X\to Y$ is a map, set:
\begin{align*} 
   \mathsf{Bd}(f) &= \{y\in Y\,:\,f^{-1}(\{y\})\text{ is Lindel\"of in }X\}.
\end{align*}

\begin{lemma}\label{lemma:unraveling}
  Let $X,Y$ be spaces. Then 
  $$ \mathsf{BR}(X,Y)\text{ holds } \Longleftrightarrow f^{-1}(\mathsf{Bd}(f))
     \text{ is contained in a Lindel\"of subset for any }f:X\to Y .$$
\end{lemma}
\proof
Essentially by unraveling the definitions.
Here are some details for the skeptics.
Suppose that $\mathsf{BR}(X,Y)$ holds. 
Let $f:X\to Y$ be given and let $Z$ be Lindel\"of such that for any Lindel\"of $W\supset Z$, $f(X-W) = f(X-Z)$.
Let $y\in\mathsf{Bd}(f)$, then $W=f^{-1}(\{y\})\cup Z$ is Lindel\"of.
If $y$ has a preimage outside of $Z$, then $f(X-W)$ does not contain $y$ while $f(X-Z)$ does,
a contradiction. Hence $f^{-1}(\mathsf{Bd}(f))\subset Z$.\\
Conversely, let $f:X\to Y$ be given such that $f^{-1}(\mathsf{Bd}(f))$ is contained in a Lindel\"of subset $Z$ of $X$. 
Let $W\supset Z$ be also Lindel\"of.
If $y\in f(X-Z)$, then $f^{-1}(\{y\})$ is non-Lindel\"of, and hence 
$f^{-1}(\{y\})\not\subset W$. It follows that $y\in f(X-W)$.
This shows that $ \mathsf{BR}(X,Y)$ holds.
\endproof

\begin{example}\label{ex:pasAS}
  There are  Type I locally metrizable spaces $H^-$ and $H^+$ such that:\\
  (a) $\mathsf{S_{cl}}(H^-)$ and thus $\mathsf{L_{cl}}(H^-)$ hold
      but $\mathsf{BR}(H^-)$ does not;\\
  (b) $\mathsf{BR_{cl}}(H^+)$ holds  but $\mathsf{L}(H^+)$ does not;\\
  (c) $H=H^-\sqcup H^+$, the topological disjoint sum of
      $H^-$ and $H^+$, has a partition into closed non-Lindel\"of subsets such that for
      each map $f: H \to\R$ there is a Lindel\"of subset $Z$ such that 
      $f$ is constant on each member of the partition outside of $Z$
      but none of $\mathsf{P}( H )$ holds for $\mathsf{P}\in\{\mathsf{EC},\mathsf{L},\mathsf{BR},\mathsf{S}\}$.
\end{example}
(Recall that for Type I spaces the properties are equivalent to their $\mathsf{cl}$-counterparts, we worded the theorem in the strongest form.)
\proof[Details]
The examples are simple subspaces of $\omega_1\times [0,1]$.
Take an uncountable subset $S=\{s_\alpha\,:\,\alpha\in\omega_1\}$ of $[0,1]$ with dense complement (all the $s_\alpha$ are distinct). 
Then set:
\begin{align*} 
   H^+ &= \bigcup_{\alpha\in\omega_1}[\alpha,\omega_1)\times\{s_\alpha\}\\
   H^- & = \omega_1\times [0,1] - H^+.
\end{align*}
(The intervals $[\alpha,\omega_1)$ are taken in $\omega_1$.)
Then both are a Type I spaces with canonical covers $H^-_\alpha,H^+_\alpha$ 
given by the intersection of the space with $[0,\alpha)\times[0,1]$.\\
(a)
Let $f:H^-\to[0,1]$ be the projection on the second coordinate. 
Then $\mathsf{Bd}(f)=S$ and 
$f^{-1}(\mathsf{Bd}(f))$ is unbounded (it contains $[0, \alpha)\times\{s_\alpha\}$ for each $\alpha$);
hence, $\mathsf{BR}(H^-)$ does not hold.
To see that $\mathsf{S}(H^-)$ does hold, take a countable dense subset 
$Q=\{q_n\,:\,n\in\omega\}$ of $[0,1]-S$. Since $S$ has dense complement, $\wb{Q} = [0,1]$. 
Given a map $g:H^-\to\R$, there is $\beta$ such that $g$ is constant above $\beta$ on each horizontal $\omega_1\times\{q_n\}$,
and thus by density on every horizontal line, even those who do not go the entire length.
Hence the retraction $r(\langle\alpha,t\rangle) = \langle \min\{\alpha,\beta\} , t\rangle$ satisfies $g\circ r = g$.\\
(b) 
The fact that $S$ has dense complement is actually irrelevant for this part.
The projection on the second coordinate shows that
$\mathsf{L}(H^+)$ does not hold: new values are introduced as far as one wants.
To see that $\mathsf{BR}(H^+)$ does hold, use the same argument as in (a) with
a countable dense subset $\{s_{\alpha_n}\,:n\in\omega\}$ of $S$ to show that there is $\beta$ such that
$f$ is eventually constant above $\beta$ on every horizontal line.
Take $Z = \wb{H^+_\beta}$, then any Lindel\"of $W\supset Z$ satisfies $f(H^+-Z) = f(H^+-W)$.
\\
(c) Notice that while $H^+\cup H^- = \omega_1\times [0,1]$,
    the topological disjoint sum $ H =H^+\sqcup H^-$
    has a finer topology, as both $H^+$ and $H^-$ are clopen in $ H $.
    Partition $ H $ into its horizontal lines. Then the claimed properties are immediate from the arguments above.
\endproof

\vskip .3cm
It is not possible to find a manifold having the same properties as $H^+$: 
Theorem \ref{thm:manifoldsL} implies that $\mathsf{L}(X)$ holds for any manifold $M$.
We however do not know the answer to the following question.

\begin{q}
    \label{q:manifoldSBR}
    Is there a manifold $M$ such that $\mathsf{S}(M)$ holds
    but not $\mathsf{BR}(M)$~? 
\end{q}

See section \ref{sec:BRL} for more on $\mathsf{BR}$ and sections \ref{sec:lackretr}--\ref{sec:SscwH} 
for more on $\mathsf{S}$. 
While $\mathsf{EC}(X,Y)$ (and a fortiori $\mathsf{EC_{cl}}(X,Y)$)
seems to be a stronger property than $\mathsf{S}(X,Y)$, 
it does not imply it since $X$ might lack retractions on sufficiently large subspaces
(see Example \ref{ex:notcl}). Actually, $\mathsf{EC}(X)$
does imply the negation of a property weaker than $\mathsf{S}(X)$ if $X$ is a longpipe.
This is proved in Section \ref{sec:lackretr},
see Theorem \ref{thm:ECnotS} and Example \ref{ex:noret}. This is
the only implication not in Lemma \ref{lemma:arrows} and not ruled out by our examples so far.

Let us end this section by settling some general (almost) trivialities.
Firstly, the following obviously holds (as already noted in \cite[Lemma 4.3.34]{MardaniThesis} for a space 
eventually constant in another):
\begin{lemma}\label{obviouslemma}
   Let $X,Y,Z$ be spaces such that there is a continuous $1$-to-$1$ map $Y\to Z$.
   Then $\mathsf{P}(X,Z) \Longrightarrow  \mathsf{P}(X,Y)$ for 
   $\mathsf{P}\in\{\mathsf{EC},\mathsf{L},\mathsf{BR},\mathsf{S}\}$ or the $\mathsf{cl}$ and $\mathsf{cpt}$ versions.
\end{lemma}

\begin{cor}\label{obviouscor} 
  Let $\tau\supset\rho$ be Hausdorff topologies on $Y$.
  Then $\mathsf{P}(X,\langle Y,\rho\rangle) \Longrightarrow  \mathsf{P}(X,\langle Y,\tau\rangle)$ for 
  $\mathsf{P}\in\{\mathsf{EC},\mathsf{L},\mathsf{BR},\mathsf{S}\}$ or the $\mathsf{cl}$ 
  and $\mathsf{cpt}$ versions.
\end{cor}
\proof $id:\langle Y,\tau\rangle\to\langle Y,\rho\rangle$ is a $1$-to-$1$ continuous map.
\endproof

Finally, the next lemmas are also almost immediate.

\begin{lemma}
   \label{lemma:Limage}
   Let $X,Y,Z$ be spaces, $f:X\to Z$ be continuous and onto, and 
   $\mathsf{P}\in\{\mathsf{EC},\mathsf{L}\}$.
   Then $\mathsf{P}(X,Y)\Rightarrow\mathsf{P}(Z,Y)$,
   and $\mathsf{P_{cpt}}(X,Y)\Rightarrow\mathsf{P_{cpt}}(Z,Y)$.
\end{lemma}
\begin{proof}
   Let $g:f(X)\to Y$ be continuous,
   $W\subset X$ and $V=f(W)$.
   If $g\circ f$ is constant outside of $W$, then $g$ is constant outside of $V$.
   If $g\circ f(W)=g\circ f(X)$ then $g(V)=g(Z)$.
   The result follows since $V$ is Lindel\"of [resp. compact] whenever $W$ is Lindel\"of [resp. compact].
\end{proof}

\begin{lemma}
   \label{lemma:QBRcountable}
   Let $X=\cup_{\alpha\in\omega_1}X_\alpha$ 
   be a Type I space and $Y$ be a countable space with no separation axiom assumed. 
   Then $\mathsf{L}(X,Y)$ and $\mathsf{BR}(X,Y)$ hold.
\end{lemma}
\proof
   Let $f:X\to Y$ be given. Then $f(X_\alpha)$ and
   $f(X-X_\alpha)$ are respectively increasing and decreasing $\omega_1$-sequences of countable sets,
   they must then stagnate above some $\alpha$.   
\endproof

\begin{lemma}\label{lemma:clopendiscrete}
   Let $X$ be a space containing a clopen uncountable discrete subspace $A$ and let $Y$ be a space. Then the following hold.\\
   (a) If $|Y|\ge\aleph_1$, then none of $\mathsf{L}(X,Y)$, $\mathsf{S}(X,Y)$ and $\mathsf{BR}(X,Y)$ do hold.\\
   (b) If $|Y|\ge 2$, then $\mathsf{EC}(X,Y)$ does not hold.\\
\end{lemma}
\proof
   \ \\
   (a) By assumption, we can define a map $f:X\to Y$ which is $1$-to-$1$ on $A$ and constant on $X-A$.
       Then $f$ contradicts $\mathsf{L}(X,Y)$ (and thus $\mathsf{S}(X,Y)$) and $\mathsf{BR}(X,Y)$.
       \\
   (b) Partition $A$ into two clopen discrete uncountable subsets $A_0,A_1$.
    Define $f:X\to Y$ that 
    sends $A_0$ to one of the points of $Y$ and $X-A_0$ to another one.
    This defines a continuous map which contradicts $\mathsf{EC}(X,Y)$
\endproof


\section{$\mathsf{EC}(X)$ and $\mathsf{EC}(\omega_1,Y)$}
\label{sec:isocompact}

In this section, we mainly investigate the $\mathsf{EC}$ property.
Theorem \ref{thm0} shows in particular that $\mathsf{EC}(\omega_1)$ holds.
We first investigate which spaces satisfy $\mathsf{EC}(X)$,
and when does $\mathsf{EC}(X)$ imply $\mathsf{EC}(X,Y)$.
(Our results have similarities with those in \cite[Section 7]{Nyikos:1992}.)
Then, we look for properties of $Y$ that
imply or are implied by $\mathsf{EC}(Z,Y)$ when $Z$ is similar (in some way) to $\omega_1$.
The notions of C-closedness and isocompactness will be central (see definitions below).


\subsection{Consequences of $\mathsf{EC}(X)$}

It is well known and easy to prove that
finite unions and at most countable intersections of $0$-sets are $0$-sets,
and finite intersections and countable unions of co-$0$-sets are co-$0$-sets.
If $Y$ is perfectly normal (in particular, metric), then the preimage of a closed subset of $Y$ is a $0$-set.

\begin{lemma}\label{lemma:closeddisjoint}
   Let $X$ be a space, and $E\subset D\subset X$ be subspaces such that $E$ is Lindel\"of and $D$ is non-Lindel\"of.
   Then the following hold.\\
   (a) There is an open $U\supset E$ such that $D-U$ is non-Lindel\"of.\\
   (b) If $X$ is Tychonoff and $D$ is a $0$-set, 
       then there is an open $U\supset E$ such that $D-U$ is a non-Lindel\"of $0$-set.
\end{lemma}
\proof
   Since $D$ is non-Lindel\"of and $E$ Lindel\"of, let 
   $\mathcal{U}$ be a cover of $D$ without countable subcover and let $\mathcal{U}_0\subset\mathcal{U}$
   be a countable subcover of $E$. Then $D-\cup\mathcal{U}_0$ is non-Lindel\"of, which proves (a).
   For (b), since $X$ is Tychonoff for each $x\in E$ we may fix $g_x:X\to [0,1]$
       such that $g_x(x) = 1$ and $g_x$ is $0$ outside of $\cup\mathcal{U}_0$.
       Let $\{x_n\,:\,n\in\omega\}$ be such that $\mathcal{W}=\{g_{x_n}^{-1}((0,1])\,:\,n\in\omega\}$ is a cover of $E$.
       Then $\cup\mathcal{W}$ is a countable union of co-$0$-sets and hence a co-$0$-set and is 
       included in $\cup\mathcal{U}_0$.
       It follows that $D-\cup\mathcal{W}$ is a $0$-set which is non-Lindel\"of since
       it contains the closed non-Lindel\"of subset $D-\cup\mathcal{U}_0$.
\endproof

\begin{defi}\label{def:ICI0}
  We say that a space $X$ has property $\mathcal{DC}$ [resp. $\mathcal{DO}$]
  iff given two disjoint closed subsets [resp. $0$-sets] of $X$, at least one of them is Lindel\"of.
\end{defi}


\begin{thm}
  \label{thmDir}
  Let $X$ be a space.
  Then the properties below are related as follows.
  \begin{center}
  (1)$\Longleftrightarrow$(2)$\Longleftrightarrow$(3a)$\Longleftarrow$(3b)$\Longleftarrow$(4a)$\Longleftrightarrow$(4b)
  \end{center}
  Moreover, if $X$ is Tychonoff then (3a)$\Longleftrightarrow$(3b) and if $X$ is normal, 
  all properties are equivalent.

   (1) $\mathsf{EC}(X)$ holds,

   (2) $\mathsf{EC}(X,Y)$ holds when $Y$ is a metric space,

   (3a) $X$ satisfies $\mathcal{DO}$,

   (3b) given two non-Lindel\"of $0$-sets of $X$, their intersection is non-Lindel\"of,

   (4a) $X$ satisfies $\mathcal{DC}$,

   (4b) given two closed non-Lindel\"of subsets of $X$, their intersection is non-Lindel\"of.
\end{thm}
\proof 
   If $X$ is Lindel\"of, we have nothing to do, hence we assume that $X$ is not Lindel\"of in what follows. \\
   (2)$\Rightarrow$(1),  
    (3b)$\Rightarrow$(3a),
    (4b)$\Rightarrow$(3b),
    (4b)$\Rightarrow$(4a) are all immediate.\\
   (1)$\Rightarrow$(3a). 
       Let $f:X\to\R$ be given.
       If $A=f^{-1}(\{0\})$ is non-Lindel\"of, 
       then since $\mathsf{EC}(X)$ holds $f$ must be eventually constant on $0$.
       Hence, $f^{-1}(\R-\{0\})$ is contained in a Lindel\"of subset $Z$.
       Any closed $B$ disjoint from $A$ is contained in $Z$ and hence Lindel\"of.\\
   (3a)$\Rightarrow$(2). 
      Recall that the preimage of a closed set of a metric space is a $0$-set.
      Let $f:X\to Y$ be given, with $Y$ a metric space.
      Then $f(X)$ has countable extent since any uncountable closed discrete subspace
      $D\subset f(X)$ can be partitioned in two disjoint such subspaces whose
      preimages yield disjoint non-Lindel\"of $0$-sets of $X$.
      By Lemma \ref{lemma:metriceq}
      $f(X)$ is hereditarily Lindel\"of.
      Let $B(y,\epsilon)$ denote the open ball of radius $\epsilon$ around $y$ in $f(X)$.
      For each $\epsilon>0$, by Lindel\"ofness of $f(X)$ and non-Lindel\"ofness of $X$ there is at least one $y\in f(X)$
      such that $f^{-1}(\wb{B(y,\epsilon)})$ is non-Lindel\"of.
      Moreover, if $f^{-1}(\wb{B(y,\epsilon)})$ is non-Lindel\"of, then $f^{-1}(f(X)-\wb{B(y,\epsilon)})$
      is Lindel\"of, because it is a countable union of $0$-sets disjoint from 
      $f^{-1}(\wb{B(y,\epsilon)})$.
      Let $y_0$ be such that $f^{-1}(\wb{B(y,1)})$ is non-Lindel\"of.
      For each $n\in\omega$, choose $y_{n+1}\in f(X)\cap\wb{B\left(y_n,\frac{1}{n}\right)}$
      such that $f^{-1}\left(\wb{B\left(y_{n+1},\frac{1}{n+1}\right)}\right)$ is non-Lindel\"of.
      Then $f(X)\bigcap\cap_{n\in\omega}\wb{B\left(y_n,1/n\right)}$ 
      is non-empty (otherwise $f(X) = \cup_{n\in\omega} (f(X)-\wb{B(y_n,1/n)})$ would have Lindel\"of
      preimage) and thus contains exactly one point $y$.
      By construction, the complement of $y$ has a Lindel\"of preimage,
      hence $f$ is eventually constant on $y$.
    \\
    (4a)$\Rightarrow$(4b) Let $C_1,C_2\subset X$ be closed and non-Lindel\"of. 
       By way of contradiction assume that $C_1\cap C_2$ is (at most) Lindel\"of.
       By Lemma \ref{lemma:closeddisjoint} (a), there are open $U,V$ both containing $C_1\cap C_2$
       such that $C_1-U$ and $C_2-V$ are disjoint non-Lindel\"of closed subsets of $X$, contradicting $\mathcal{DC}$.\\
    (3a)$\Rightarrow$(3b) when $X$ is Tychonoff.
       Same proof as (4a)$\Rightarrow$(4b), using Lemma \ref{lemma:closeddisjoint} (b).
       \\
    (3a)$\Rightarrow$(4a) if $X$ is normal. Given two disjoint closed sets $E,F$ in $X$ we obtain two disjoint $0$-sets 
       $A\supset E, B\supset F$
       with a Urysohn fonction. If both $E$ and $F$ are non-Lindel\"of, then so are $A$ and $B$, contradicting (3a).
\endproof

Half of the next result is a direct consequence of Theorem \ref{thmDir}.

\begin{thm}\label{thm:normal}
   If $\mathsf{EC}(X)$ holds and $X$ is either normal or $\aleph_1$-scwH and Tychonoff, then $X$ is $\omega_1$-compact.
\end{thm}

Recall that a collection of subsets of a space $X$ is {\em discrete} 
iff for each $x\in X$ there is an open set containing it which intersects at most one member of the collection.
In a discrete collection $\{U_\alpha\,:\,\alpha\in\lambda\}$, $\cup_{\alpha\in E}\wb{W_\alpha}$ is 
closed for any $E\subset\lambda$ and $W_\alpha\subset U_\alpha$. A space is
{\em $\kappa$-[strongly] collectionwise Hausdorff 
(abbreviated $\kappa$-[s]cwH)}  iff any closed discrete subset $\{x_\alpha\,:\,\alpha\in\lambda\}\subset X$
of size $\lambda\le\kappa$
can be expanded to a disjoint [resp. discrete] 
collection of open sets $\mathcal{U}=\{U_\alpha\,:\,\alpha\in\lambda\}$ with $x_\alpha\in U_\alpha$.
Such a $\mathcal{U}$ is called a {\em disjoint [resp. discrete] expansion} of $D$.
If $X$ is $\kappa$-[s]cwH for each $\kappa$, we say that $X$ is [s]cwH.
A normal $\kappa$-cwH space is $\kappa$-scwH, as well known.

\begin{proof}
   A closed discrete subset $D$ (of cardinality $\aleph_1$) can be partitioned in two such subsets $D_1,D_2$, 
   which contradicts $\mathcal{DC}$. If $X$ is normal, the result follows by Theorem \ref{thmDir}.
   If $X$ is Tychonoff and $\aleph_1$-scwH, then $D=D_1\cup D_2$ is contained in a $0$-set $E$
   which is the union of a discrete family $\{E_d\,:\,d\in D\}$, where $d\in E_d$ for each $d\in D$.
   Indeed, take a discrete expansion $\mathcal{U}=\{U_d\,:\,d\in D\}$ of $D$ and
   $g_d:X\to[0,1]$ which is $1$ on $d$ and $0$ outside of $U_d$.
   By discreteness, $g_i=\sum_{d\in D_i}g_d$ ($i=1,2$) is a continuous function,
   and $D_i=g_i^{-1}(\{1\})$ ($i=1,2$) contradict $\mathcal{DO}$,
   which yields the result by Theorem \ref{thmDir} again.
\end{proof}

Notice that there are non-$\omega_1$-compact spaces $X$ such that 
$\mathsf{EC}(X,\R)$ holds;
Example \ref{ex:notclbis} is one, Example \ref{ex:ECnotomega_1cpct} below is another, which is moreover Type I.

\begin{figure}[h]
  \begin{center}
    \epsfig{ figure = 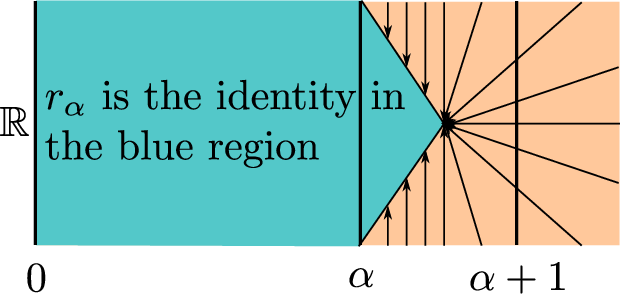, width = .4\textwidth} \quad 
    \caption{The retraction $r_\alpha$ of Example \ref{ex:ECnotomega_1cpct}.}
    \label{fig:r}
  \end{center}
\end{figure}

\begin{example}[{Nyikos, in effect}]
   \label{ex:ECnotomega_1cpct}
   There is a non-$\omega_1$-compact and non-normal Type I surface $M$ such that $\mathsf{EC_{cl}}(M)$ and 
   $\mathsf{S}(M)$ hold.
\end{example}
\proof[Details]
   This example is due to Nyikos and is described in \cite[section 6 \& 7]{Nyikos:1992}. 
   Every property we will claim in this proof is proved in Nyikos' paper to which we refer for details. 
   Let us give a quick description.
   We first consider a tangent bundle $T\LL_+$ of $\LL_+$ given by a smoothing.
   The removal the $0$-section 
   $L_0$ of $T\LL_+$, which is a copy of $\LL_+$, 
   disconnects $T\LL_+$ into two homeomorphic connected submanifolds $T^+$ and $T^-$. 
   Since $T\LL_+$ is not trivial, $T^+$ does not contain a copy of $\LL_+$, 
   actually there is no unbounded map $\LL_+\to T^+$.
   Write $\pi:T\LL_+\to\LL_+$ for the bundle projection.
   Both $T\LL_+$ and $T^+$ are Type I manifolds, a canonical cover being given 
   by the fibers $\{U_\alpha = \pi^{-1}(0,\alpha)\,:\,\alpha\in\omega_1\}$
   (and their intersections with $T^+$ for this latter space).
   The choice of the smoothing is important since $T\LL_+$ and $T^+$ can
   exhibit quite different topological properties depending on it.
   Here, as we need 
   $\mathsf{EC_{cl}}(T^+)$ (and thus $\mathsf{EC_{cl}}(T^-)$) to hold, 
   we take what Nyikos calls a {\em smoothing of class 7}, see \cite[Section 6 \& 7]{Nyikos:1992}
   for the construction and the proof. In particular, Nyikos shows that $T^+$ is not normal.
   Set $M$ to be $T\LL_+ - \omega_1$, where $\omega_1$ is seen as a subset of $L_0$.
   Then the subspace $\{\langle\alpha,0.5\rangle\,:\,\alpha\in\omega_1\}$ of $L_0$ is closed
   discrete in $M$.
   It follows that $M$ is neither $\omega_1$-compact nor normal, but $\mathsf{EC_{cl}}(M)$ holds. Indeed,
   any $f:S\to\R$ is eventually constant on both $T^+$ and $T^-$, by density and connectedness
   it is eventually constant
   on all of $M$. 
   \\
   By construction, $U_\alpha$ is homeomorphic to $(0,\alpha)\times(-1,1)$ for each $\alpha$.
   Given $\alpha$, one may easily define $r_\alpha:T\LL_+\to U_{\alpha+1}$
   which is the identity on $\wb{U_\alpha}$ and sends all of $T\LL_+ - U_{\alpha+1}$ to a point
   in $U_{\alpha+1} - \wb{U_\alpha}$, as seen in Figure \ref{fig:r}.
   Then $r_\alpha$ is well defined on $M$ for each $\alpha$.
   Together with $\mathsf{EC_{cl}}(M)$, it implies that $\mathsf{S}(M)$ holds as well.
\endproof

\vskip .3cm
We now investigate for which spaces $X$ and $Y$ we have that $\mathsf{EC}(X)$ does imply $\mathsf{EC}(X,Y)$.
Theorem \ref{thmDir} deals with the case where $Y$ is metric and shows that we
may assume that $X$ satisfies $\mathcal{DO}$.
It might be the right time for a trivial remark: if $Y\supset X$ and $X$ is not Lindel\"of, then 
the inclusion $i:X\to Y$ shows that $\mathsf{EC}(X,Y)$ does not hold. Since any Tychonoff space 
is a subspace of a compact space (its \v Cech-Stone compactification),
it would be silly to try to show that $\mathsf{EC}(X,Y)$ holds 
for some class of non-Lindel\"of Tychonoff spaces $X$ by relying only on 
covering properties of $Y$, we also need some local properties.
We shall see that having $G_\delta$ points plays a central role.

We shall concentrate on cases where $X$ is countably compact (or ``almost'').
Recall first that it is not true (at least consistently) that if $X$ is countably compact
then $\mathsf{EC}(X)$ implies $\mathsf{EC_{cpt}}(X)$, as shown by Example \ref{ex:notcl} under 
$\mathfrak{p}=\mathfrak{c}$. We show in Theorem \ref{thm:cl-cpt2} below that 
$\mathsf{EC_{cl}}(X)$ implies $\mathsf{EC_{cpt}}(X)$ when $X$ is pseudocompact.

\begin{q}
   Is there an example in {\bf ZFC} of a countably compact space $X$ such that
   $\mathsf{EC}(X)$ holds but $\mathsf{EC_{cl}}(X)$ does not~?
\end{q}

A space is {\em C-closed} iff any countably compact subspace is closed.
The class of C-closed spaces contains in particular the hereditarily meta-Lindel\"of spaces,
the sequential spaces and the regular spaces with $G_\delta$ points, see e.g. \cite{IsmailNyikos}.
The following lemma is immediate, although it implies for instance \cite[Lemma A.33]{GauldBook}. 
\begin{lemma}
   \label{lemmaCclosed}
   Let $X$ be countably compact and $Y$ be C-closed.
   Then any 
   $f:X\to Y$ is a closed map.
\end{lemma}

Theorem \ref{thm1EC} below is a slight generalization of \cite[Lemmas B.29--30]{GauldBook}, with a simpler proof.
We first state a useful technical lemma from \cite[Lemma 2.2]{mesziguesContractibility}.

Given a chain cover $\mathcal{U}=\{U_\alpha\,:\,\alpha\in\kappa\}$ of $X$,
a subset of $X$ is $\mathcal{U}$-unbounded iff it is not contained in any $U_\alpha$.
It is well known that 
a space is compact iff any chain cover is trivial, that is, has a maximal element which is the whole space.
(To our knowledge, this is due to Alexandroff and Urysohn \cite{AlexandroffUrysohn} in 1929.)

\begin{lemma}[{\cite[Lemma 2.2]{mesziguesContractibility}}]
   \label{lemma:cunbd}
   Let $f:X\to Y$ be such that $X$ is countably compact and non-compact, $f(X)\subset Y$ is compact, and $Y$ is C-closed.
   Let $\mathcal{U}=\{U_\alpha\,:\,\alpha\in\kappa\}$, with $\kappa\ge\omega_1$ a regular cardinal, be a chain cover of $X$.
   Then, there is $c\in Y$ such that $f^{-1}(\{c\})$ is $\mathcal{U}$-unbounded (and hence non-Lindel\"of).
\end{lemma}
We include a proof for completeness.
\proof
   For each $\alpha\in\kappa$ choose $x_\alpha\in X-U_\alpha$ and set $E_\alpha = \{x_\beta\,:\,\beta\ge\alpha\}$.
   Then $\wb{E_\alpha} \cap U_\alpha = \varnothing$ and each $E_\alpha$ is $\mathcal{U}$-unbounded.
   Since $\wb{E_\alpha}$ is countably compact and $Y$ C-closed, 
   $f(\wb{E_\alpha})\subset f(X)$ is closed and thus compact in $Y$.
   Any finite intersection of $E_\alpha$ is nonempty, hence
   by compactness $\cap_{\alpha\in\kappa} f(\wb{E_\alpha})\not=\varnothing$.
   Take $c$ is this intersection, then $f^{-1}(\{c\})\cap \wb{E_\alpha}\not=\varnothing$ for each $\alpha$, which shows that
   $f^{-1}(\{c\})$ is $\mathcal{U}$-unbounded.
\endproof

We state (and prove) the next theorem in an almost absurd level of generality.
The {\em linear Lindel\"of number} $\ell L(X)$ of a space $X$ is the smallest $\kappa$
such that any chain cover of $X$ has a subcover of cardinality $\le\kappa$.
It is well known that a countably compact space with $\ell L(X)=\aleph_0$ is compact.
A space is {\em $[\aleph_0,\kappa]$-compact} iff any cover by at most $\kappa$ open sets has a countable subcover.

\begin{thm}\label{thm1EC}
  Let $X$ be a countably compact non-compact space with linear Lindel\"of number $\kappa\ge\omega_1$.
  Let $Y$ be a $[\aleph_0,\kappa]$-compact C-closed space. Then the following holds.
  \\
  (a) For each $f:X\to Y$, $f(X)$ is compact and there is $c\in Y$ with a non-Lindel\"of preimage in $X$. \\
  (b) If $X$ satisfies $\mathcal{DC}$ and
      $Y$ has $G_\delta$ points, then $\mathsf{EC}(X,Y)$ holds.\\
  (c) If $X$ satisfies $\mathcal{DO}$ and 
      $Y$ is Tychonoff with $G_\delta$ points, then $\mathsf{EC}(X,Y)$ holds.
\end{thm}

This theorem is very similar to Lemma 6.11 in \cite{mesziguesContractibility}. We will 
however repeat the proof for convenience of the reader.
Recall that as said earlier a regular space with $G_\delta$ points is C-closed, hence in (c) there is some redundancy
in the assumptions. 
Actually, the image of $X$ in $Y$ is first countable in (b) and (c) since a countably compact space with 
$G_\delta$ points is first countable.
Typical examples of $X$ satisfying the hypothesis of (b) are $\omega_1$ and $\LL_{\ge 0}$.
Notice that if $X$ is Type I, then $\kappa=\omega_1$, hence $Y$ may be assumed to be $[\aleph_0,\aleph_1]$-compact,
which is equivalent to the property that for any uncountable subset $B\subset Y$ there is a point of $Y$ 
any of whose neighborhoods
contains uncountably many members of $B$; such a point is called a {\em condensation point of $B$}.
An immediate corollary of Theorems \ref{thmDir} \& \ref{thm1EC} is the following.

\begin{cor}
   \label{cor:ccGdelta}
   Let $X$ be countably compact and $Y$ be a Tychonoff Lindel\"of space with $G_\delta$ points.
   Then $\mathsf{EC}(X)$ implies $\mathsf{EC}(X,Y)$.
\end{cor}

\proof If $X$ is compact, then both $\mathsf{EC}(X)$ and $\mathsf{EC}(X,Y)$ hold.
       Otherwise, $X$ is countably compact and non-compact
       hence $\ell L(X)\ge\omega_1$. If $\mathsf{EC}(X)$ holds then 
       by Theorem \ref{thmDir} $X$ satisfies $\mathcal{DO}$.
       $Y$ is $[\aleph_0,\kappa]$-compact for each $\kappa$ by Lindel\"ofness, and thus
       $\mathsf{EC}(X,Y)$ follows by Theorem \ref{thm1EC} (c).
\endproof

\proof [Proof of Theorem \ref{thm1EC}]
\ \\
   (a) Clearly, $f(X)$ is countably compact and closed in $Y$ and $\ell L(f(X))\le\ell L(X)=\kappa$.
       Given a chain cover $\mathcal{U}$ of $f(X)$, there is a subcover of cardinality $\le\kappa$. 
       By $[\aleph_0,\kappa]$-compactness of $Y$,
       there is a countable subcover and hence a finite one, so $f(X)$ is compact. 
       Lemma \ref{lemma:cunbd} yields the result.\\
   (b) 
       Given $f:X\to Y$, by (a) there is $c\in Y$ with a non-Lindel\"of preimage.
       Let $V_n\subset Y$ be open sets such that $\{c\}=\cap_{n\in\omega}V_n$. 
       By $\mathcal{DC}$, 
       $f^{-1}(Y-V_n)$ must be Lindel\"of, and hence 
       $f^{-1}(Y-\{c\}) = \cup_{n\in\omega} f^{-1}(Y-V_n)$ is also Lindel\"of.
       \\
   (c) Define $V_n$ as in (b) and take $g:Y\to\R$ which is $0$ on $c$ and $1$ on $Y-V_n$.
       Since $(g\circ f)^{-1}(\{0\})\supset f^{-1}(\{c\})$ is non-Lindel\"of,
       $(g\circ f)^{-1}(\{1\})\supset f^{-1}(Y-V_n)$ must be Lindel\"of, and we conclude as in (b).
\endproof

Another quick corollary of Theorem \ref{thm1EC} is the following.
\begin{cor} \label{corgenMardani}
   Let $X$ be a countably compact space satisfying $\mathcal{DC}$.
   Then $\mathsf{EC}(X,Y)$ holds whenever $Y$ is a C-closed isocompact space with $G_\delta$ points.
\end{cor}

A space is {\em isocompact} iff every closed countably compact subset is compact. 
Hence a space is C-closed and isocompact iff each countably compact subset is compact.
Such a space is sometimes called {\em hereditarily isocompact}.
For more on isocompact spaces, see for instance \cite{Bacon:1970,ChoPark}. Note for instance that
a meta-Lindel\"of space is isocompact.
\proof[Proof of Corollary \ref{corgenMardani}]
   Again, $\ell L(X)\ge\omega_1$.
   Let $f:X\to Y$.
   Since $Y$ is C-closed, $Z=f(X)$ is a closed countably compact subspace of $Y$.
   Thus $Z$ is compact, and in particular $[\aleph_0,\kappa]$-compact.
   We may thus apply (b) of Theorem \ref{thm1EC}.
\endproof

\begin{rem}
   Any space of the form $Z\cup X$, where $Z$ is Lindel\"of and 
   $X$ satisfies the hypothesis of Theorem \ref{thm1EC} (or Corollaries \ref{cor:ccGdelta}--\ref{corgenMardani}) 
   also satisfies
   its conclusions.
\end{rem}

A particular case of Corollary \ref{corgenMardani} is the following:

\begin{cor}[A. Mardani, Prop. 4.3.12 and 4.3.23 in \cite{MardaniThesis}]
   \label{corMardani}
   Let $Y$ be either a realcompact space with $G_\delta$ points or a discrete space.
   Then $\mathsf{EC}(\omega_1,Y)$ holds.
\end{cor}
\proof
If $Y$ is discrete, each countably compact subspace is finite and hence closed.
Recall that a realcompact space is Tychonoff and hence C-closed if it has $G_\delta$ points.
A countably compact closed subset of a realcompact space is compact,
see e.g. \cite[3.11.1 \& 3.11.4]{Engelking}.
It follows that a realcompact space is isocompact.
\endproof

Let us now show that if $Y$ is confined to regular spaces, then there is
a way to obtain (b) (and (a)) of Theorem \ref{thm1EC}
without the assumption that $X$ is countably compact 
(union a Lindel\"of space)
whenever $X$ is `narrow' enough.
The proof of the next theorem is essentilally done in \cite[Lemmas 6.10--6.11]{mesziguesContractibility}, 
but again, we repeat it for completeness.

\begin{thm}\label{thm:ecstat}
   Let $X$ be a space 
   with $\kappa = \ell L(X)>\aleph_0$ and such that
   any countable family of closed non-Lindel\"of subsets of $X$ has a non-Lindel\"of intersection. 
   Let $Y$ be a regular $[\aleph_0,\kappa]$-compact 
   space with $G_\delta$ points.
   Then $\mathsf{EC}(X,Y)$ holds.
\end{thm}
A typical example for $X$ is a stationary subset of $\omega_1$ with subspace topology.

\proof
   By assumption $X$ satisfies $\mathcal{DC}$.
   The proof of (b) in Theorem \ref{thm1EC} above does not 
   use countable compactness and only relies on $\mathcal{DC}$, the fact that points are $G_\delta$ and
   that there is $c\in Y$ with non-Lindel\"of preimage.
   We now show that the latter holds in our case as well.
   Since $\kappa > \aleph_0$
   there is a chain cover $\mathcal{U}=\{U_\alpha\,:\,\alpha\in\lambda\}$ ($\omega_1\le\lambda\le\kappa$)
   of $X$ such that $X\not\subset U_\alpha$ for each $\alpha$, with $\lambda$ a regular cardinal.
   Take $x_\alpha$ in $X-U_\alpha$ and set $A_\alpha=\{x_\beta\,:\,\beta>\alpha\}$.
   Then $\wb{A_\alpha}\cap U_\alpha = \varnothing$, each $\wb{A_\alpha}$ is $\mathcal{U}$-unbounded,
   and a countable intersection of $\wb{A_\alpha}$'s is nonempty.
   If $|f(A_0)|<\lambda$, there is some $c$ such that $f^{-1}(\{c\})$ is $\mathcal{U}$-unbounded.
   Else, let $c$ be a point of $f(A_0)$ such that any neighborhood of $c$
   contains $\lambda$-many points of $f(A_0)$ (which exists by $[\aleph_0,\kappa]$-compactness).
   Since $Y$ is regular and points are $G_\delta$, there are open $U_n\subset Y$, $n\in\omega$, such that 
   $\{c\} = \cap_{n\in\omega}\wb{U_n}$. Then $f^{-1}(\wb{U_n})$ is $\mathcal{U}$-unbounded
   (and thus non-Lindel\"of) in $X$ for each $n$.
   Hence, their intersection is non-Lindel\"of as well and equal to $f^{-1}(\{c\})$.   
\endproof

It seems difficult to significantly weaken the hypothesis about $G_\delta$ points in Theorem \ref{thm1EC}, as the
simple next example shows.
\begin{example}\label{ex:quotientomega_1}
  The quotient space $Y=\omega_1/\Lambda$ (where $\Lambda$ is the subspace of limit ordinals) is Fr\'echet-Urysohn (hence $C$-closed),
  compact (hence isocompact), but $\mathsf{EC}(\omega_1,Y)$ does not hold, as shown by the quotient map.
\end{example}
This example is also treated in \cite[Ex. 4.2.42]{MardaniThesis}.
The only non-isolated point of $Y$ is not a $G_\delta$ due to the pressing down lemma 
(also known as Fodor's theorem, see any book on set theory, e.g. \cite[Lemma 6.15]{Kunen}).


\subsection{$\mathsf{EC}(\omega_1,Y)$ and isocompactness of $Y$}

A canonical case in Corollary \ref{corgenMardani} is $X=\omega_1$.
We now go the other way and ask: Does $\mathsf{EC}(\omega_1,Y)$, along with some `mild' properties of $Y$, 
imply isocompactness of $Y$~?
By a mild property we mean something that is shared by many spaces 
but does not render the question trivial\footnote{We are aware that this phrase is very similar to 
``By something yellow we mean something containing a lot of yellowness''.}.
An interesting case is that of
locally compact countably tight spaces, for which
the answer depends on the axioms of set theory. Firstly,
Theorem \ref{thm:iso_PFA} just below shows that a `yes' is consistent.
It is yet another example of the ``ubiquity'' of $\omega_1$ in countably compact non-compact spaces
under the proper forcing axiom {\bf PFA}, and
is actually only a restatement of other authors old results from our point of view. Secondly, 
Example \ref{exNy} below (also due to another author) shows that `no' is also consistent.
We first show a small lemma.

\begin{lemma}\label{lemma:noncountablycpct}
  If $X=\cup_{\alpha\in\omega_1}X_\alpha$ is of Type I and countably compact, then $X_\alpha \not= \wb{X_\alpha}$ for limit $\alpha$,
  so $\cup_{\alpha\in\Lambda} (\wb{X_\alpha} - X_\alpha)$ is a perfect preimage of
  $\omega_1$.
\end{lemma}
Recall that a map is {\em perfect} iff it is closed and points have compact preimages and that $\Lambda\subset\omega_1$ is the subset of limit ordinals.

\proof
  Take a strictly increasing sequence $\alpha_n$ whose 
  limit is $\alpha$ and for each $n$ some $x_n\in X_{\alpha_{n+1}} - X_{\alpha_{n}}$.
  This sequence cannot have an accumulation point 
  in $X_\alpha$. This shows that $X_\alpha$ is not countably compact when $\alpha$ is limit.
  If $X$ itself is countably compact, then $X_\alpha$ is not closed for limit $\alpha$, hence $\wb{X_\alpha} - X_\alpha \not= \varnothing$.
  The map 
  $\cup_{\alpha\in\Lambda}(\wb{X_\alpha} - X_\alpha) \to \omega_1$
  sending each point in $\wb{X_\alpha} - X_\alpha$ to $\alpha$ is perfect and its image is homeomorphic to $\omega_1$.
\endproof

\begin{thm}[{\bf PFA}]
   \label{thm:iso_PFA} \ \\
   Let $Y$ be a countably tight space such that $\mathsf{EC}(\omega_1,Y)$ holds.
   If $Y$ is moreover either of Type I or locally compact, then  $Y$ is isocompact.
\end{thm}
\proof
Notice that any copy of $\omega_1$ is closed in a countably tight space.
Recall also that, as shown by T. Eisworth in \cite{Eisworth:2002},
a countably tight space $Y$ containing a perfect preimage of $\omega_1$
contains a copy of $\omega_1$ under {\bf PFA}, hence $\mathsf{EC}(\omega_1,Y)$ does not hold.
\\
       We assume first that $Y=\cup_{\alpha\in\omega_1}Y_\alpha$ is of Type I.
       Suppose that $Y$ is not isocompact, then $Y$ has a countably compact closed subset $Z$ which is not compact.
       Since each $\wb{Y_\alpha}$ is Lindel\"of, $Z\not\subset\wb{Y_\alpha}$ so
       $Z$ is actually a Type I countably compact 
       space with $Z_\alpha = Z\cap Y_\alpha$. 
       By Lemma \ref{lemma:noncountablycpct} $\cup_{\alpha\in\Lambda}(\wb{Z_\alpha}-Z_\alpha)$ is a perfect preimage of
       $\omega_1$ which is moreover countably tight, and hence $\mathsf{EC}(\omega_1,Y)$ does not hold under {\bf PFA}.
       Assume now that $Y$ is locally compact.
       A closed subset of $Y$ is also locally compact. By \cite[Thm 2.6]{Balogh:1989}, under {\bf PFA}
       a locally compact countably compact non-compact space contains a 
       perfect preimage of $\omega_1$, 
       and we conclude as above.
\endproof
The full force of {\bf PFA} is probably not needed in the proof, but 
more than {\bf MA+$\neg$CH} is, as shown by the following examples.
Both are due to P. Nyikos and are described in a preliminary draft he made available on his webpage \cite{Nyikos:Antidiamond}.
Since P. Nyikos unfortunately passed away, this draft will never be completed by him, but we wrote a note 
which contains details of the proofs \cite{MesziguesNyikosExample}.
The axiom $\clubsuit_C$ follows from $\diamondsuit$ and is compatible with {\bf MA+$\neg$CH}.
Another example
under $\diamondsuit$ appears in another Nyikos paper \cite[Ex. 6.17]{Nyikos:1984}.

\begin{example}[$\clubsuit_C$]
   \label{exNy}
   There is a a longpipe (thus $\omega$-bounded $2$-manifold) $Y$ which is not isocompact but $\mathsf{EC}(\omega_1,Y)$ holds.
   Also, for each $n\in\omega$, $n\ge 2$, 
   there is a closed $n$-to-$1$ closed preimage of $\omega_1$ $P$ with the same properties.
\end{example}
\proof[Some details]
The surface is described in \cite[Section 5]{Nyikos:Antidiamond} 
and the version with $n=2$ (called `sprat' by Nyikos) $P$ is given by Theorem 2.1 in the same paper and the remarks after.
Both are also described in more detail in \cite{MesziguesNyikosExample}.
They are regular and first countable and hence C-closed, Type I, 
countably compact, non-compact, and neither contain a copy of $\omega_1$.
Actually, both satisfy a property stronger than (c) of Theorem \ref{thmDir}: 
any club subset of $Y$ (resp. of $P$) contains $\wb{Y_\alpha}-Y_\alpha$ (resp. $\wb{P_\alpha}-P_\alpha$) for a club set of $\alpha$.
Moreover any bounded subset of $Y$ or $P$ embeds in $\R^2$.
\\
The next lemma shows that a map $f$ from $\omega_1$ to $Y$ or $P$ must have a bounded image, and
since any bounded subset of $Y$ or $P$ embeds in $\R^2$ and $\mathsf{EC}(\omega_1\R^2)$ holds, 
$f$ must be eventually constant.
\endproof

\begin{lemma}
   If $Y$ is a Type I C-closed space and $f:\omega_1\to Y$ is unbounded, then there is a copy of $\omega_1$ in $Y$. 
\end{lemma}
  \proof
      The image of $f$ is club by Lemma \ref{lemmaCclosed}, hence by a routine argument similar to that of \cite[Lemma 1.19]{GauldBook}, 
      the set $C =\{\alpha\in\omega_1\,:\,f(\alpha)\in \wb{Y_\alpha}-Y_\alpha\}$ is club as well. 
      Thus $C$ embeds in $Y$ and is a copy of $\omega_1$.
  \endproof

We note that there are longpipes $Y$ (which are of course non-isocompact)
satisfying $\mathsf{EC}(\LL_{\ge 0},Y)$ in {\bf ZFC},
see Example \ref{ex:noret} below.


\section{$\mathsf{S}(T)$ and Suslin trees}
\label{sec:Suslin}

A tree $T$ is a partially ordered set such that each point has a well ordered set of predecessors.
We usually denote the order by $<$, and $>,\le,\ge$ are defined as usual.
Here all trees are endowed with the order topology (also called interval topology): a basis is given
by the intervals $\{z\in T\,:\, x<z\le y\}$ for each $x,y\in T$. 
We assume that our trees are Hausdorff, 
that is, if $x,y\in T$ are at a limit level and have the same predecessors, then $x=y$.
It follows that any tree is a $0$-dimensional space.\\
A {\em chain} is a totally ordered subset and an {\em antichain} a subset with pairwise incomparable elements.
The $\alpha$-th level of $T$ consists of the members whose set of predecessors has order type $\alpha$.
Points in the $\alpha$-th level are often said to be at height $\alpha$.
The height of $T$ is the smallest ordinal $\beta$ such that the $\beta$-th level of $T$ is empty.
\\
An {\em $\omega_1$-tree} has countable levels and height $\omega_1$.
A tree is {\em Suslin} if it has height $\omega_1$ and its chains and antichains are at most countable.
Recall that Suslin trees do not exist in {\bf ZFC} alone, but do exist under 
{\bf V=L} or $\diamondsuit$, for instance.
A subset $D\subset T$ is {\em order dense} iff given any $x\in T$ there is $y\in D$ with $y>x$.
For a tree $T$ and $t\in T$, write $T_{\ge t} = \{s\in T\,:\, s\ge t\}$,
$t\upharpoonright\alpha$ for the unique predecessor 
of $t$ at level $\alpha$ (if $t$ is below the $\alpha$-th level, $t\upharpoonright\alpha=t$)
and $T_{<\alpha}$ (resp. $T_{\le\alpha}$) for the subset of elements at level $<\alpha$ (resp. $\le\alpha$).  
If $T$ is an $\omega_1$-tree then it is a Type I space and the $T_{<\alpha}$ form a canonical cover, moreover
each $T_{<\alpha}$ embeds in $\R$. (For the latter claim, notice that $T_{<\alpha}$ is
a countable, second countable $0$-dimensional space, hence embeds in the Cantor set or the rationals.)\\
Given an ordinal $\alpha$, define $r_\alpha:T\to  T_{\le\alpha}$ as $r_\alpha(t)=t\upharpoonright\alpha$.
Notice that if $\beta\ge\alpha$ and $f\circ r_\alpha = f$ for some $f:T\to Y$, then $f\circ r_\beta = f$ as well.
If $A\subset T$, denote by $A^\downarrow$ its downward closure $\{x\in T\,:\,\exists y\in A\;x\le y\}$.
If $x\in A\subset T$ we let $A_{\ge x}$ be $A\cap T_{\ge x}$.
Since we are going to look at $\omega_1$-compact subspaces of trees, let us recall the following (more or less classical) facts.

\begin{lemma}\label{lemma:treefacts}
   Let $T$ be a tree of height $\omega_1$ and $S\subset T$ be uncountable and $\omega_1$-compact in the subspace topology. Then the following hold.
   \\
   (a) A closed discrete subset of $T$ is a countable union of antichains, and an antichain of $T$ is closed discrete.\\
   (b) If $A\subset T$, $A$ has an uncountable antichain if and only if $A^\downarrow$ has one.\\
   (c) $S^\downarrow$ is the union of a countable set, at most countably many copies of $\omega_1$ and a Suslin tree.\\
   (d) There is $\alpha\in\omega_1$ such that $|(S^\downarrow)_{\ge x}|=|S_{\ge x}|=\aleph_1$ when $x$ is above 
       level $\alpha$ in $S^\downarrow$.\\
   (e) $S$ intersects a stationary subset of levels of $S^\downarrow$. Moreover, any closed subset of $S^\downarrow$ intersects
       a stationary subset of levels.\\
   (f) Let $E,F\subset S$ be closed (in $S$). If $|E\cap F|\le\aleph_0$, then $|E^\downarrow\cap F^\downarrow|\le\aleph_0$.\\
   (g) If $A\subset S^\downarrow$ is uncountable, there is $x\in S^\downarrow$ such that 
       $(S^\downarrow)_{\ge x}\subset A^\downarrow$.\\
   (h) If $D\subset S^\downarrow$ is order-dense and upward-closed, then 
       $D\supset S^\downarrow-(S^\downarrow)_{\le\beta}$ for some $\beta$.\\
   (i) If $U \subset S$ is open in $S$ such that $U$ intersects stationary many levels of 
       $(S^\downarrow)_{\ge x}$ for each $x\in S$,
       then $U\supset (S^\downarrow-(S^\downarrow)_{\le\beta})\cap S$
       for some $\beta$.
\end{lemma}
\proof
   All are probably part of the folklore, actually they are well known facts when 
   $T=S$ is a Suslin tree. Proofs of items (b) to (g) can be found in 
   \cite[{Lemmas 6.6--6.7}]{mesziguesContractibility}, and (a) in \cite[Theorem 4.11]{Nyikos:trees}, for instance.
   By considering minimal elements, it is easy to see that
   (h) holds if $S^\downarrow$ is Suslin, hence it follows by (c).
   Let us prove (i). Suppose that $U$ is as in the statement of (i).
   Let $\alpha$ be given by (d). By removing the $\alpha$-th first levels, we may assume that 
   $S_{\ge x}$ is uncountable for all $x\in S$.
   By (h) it is enough to show that $F=\{x\in S\,:\,S_{\ge x}\subset U\}$ is order-dense in $S$ 
   (or equivalently in $S^\downarrow$) above level $\alpha$.
   Fix $z\in S$ and let $W$ be open in $S^\downarrow$ such that $W\cap S = U$.
   Then $W$ intersects stationary many limit levels of $(S^\downarrow)_{\ge z}$.
   For each $x\in W$ at a limit level with $x > z$, let $\sigma(x)<x$ be such that $\{u\,:\,\sigma(x)\le u\le x\}\subset W$.
   By the pressing-down lemma for $\omega_1$-trees (see, e.g., \cite[p. 154]{Hart-Souslin})
   there is some $y\in S^\downarrow$ such that $E = \sigma^{-1}(\{y\})$ meets stationary many levels of 
   $(S^\downarrow)_{\ge z}$.
   For each $x\in E$, $\{t\in T\,:\,y\le t\le x\}$ is contained in $W$. We may assume that $y\ge z$.
   By (g), there is some $x\ge y$ such that $(S^\downarrow)_{\ge x}\subset E^\downarrow$. 
   But this means that $(S^\downarrow)_{\ge x}\subset W$, hence $S_{\ge x}\subset U$. 
   We thus proved that there is a point $x$ of $F$ above $z$, and hence $F$ is order-dense.
\endproof

The following was essentially proved by Stepr\=ans \cite{StepransTrees} when $S=T$. (His statement is weaker, 
as he only claims that $T$ is functionally countable, but he gives two proofs
which actually show more than stated.) A proof for $S\not= T$, adapted from Stepr\=ans', is given in 
\cite[{Lemma 6.7 (i)}]{mesziguesContractibility}.
We present another proof, by forcing, also based on Stepr\=ans' ideas.
\begin{thm}\label{thm:Steprans}
    \ \\    
    Let $T$ be a tree of height $\omega_1$, $S\subset T$ be uncountable and $Y$ be a space. Then the following hold.
    \\
    (a) If $S$ is $\omega_1$-compact for the subspace topology and $Y$ is submetrizable, 
        then $\mathsf{S}(S,Y)$ holds. If $S=T$ then the retraction is given by $r_\alpha$ for some $\alpha$.\\
    (b) If $|Y|\ge \aleph_1$ and $\mathsf{L}(S,Y)$ holds, then $S$ is $\omega_1$-compact.\\
    (c) If $Y$ is an uncountable submetrizable space, then $\mathsf{L}(T,Y)\Leftrightarrow\mathsf{S}(T,Y)$.
\end{thm}
Recall that a space is {\em submetrizable} iff it has a coarser metrizable topology.
\proof We prove (b) first.\\
(b) Suppose that $S$ is not $\omega_1$-compact.
    By Lemma \ref{lemma:treefacts} (a)--(b) and (d), $S$ contains an uncountable antichain $A$ consisting
    of isolated points of $S$ and thus a clopen discrete uncountable subspace.
    Then, apply Lemma \ref{lemma:clopendiscrete} (a). 
\\
(a) We may assume that $Y$ is metrizable by Corollary \ref{obviouscor}.
    Let $f:S\to Y$ be given.
    We shall see that there is some $\alpha$ such that $f$ is constant on $S_{\ge x}$ 
    when $x$ is at level (in $T$) above $\alpha$. 
    It follows that, given $y$ at level $>\alpha$, we may set 
    $r(y) = \min\{ x \le y \,:\, x\in S,\, \text{height}(x)\ge\alpha\}$.
    If $S=T$, it is enough to set $r(y) = y\upharpoonright\alpha$ (i.e. $r=r_\alpha$) since there are points at every level.
    \\
    Let us prove that there is such an $\alpha$.
    By Lemma \ref{lemma:treefacts} (c), the set $E$ containing the minimal points of 
    $$\{x\in S\,:\, S_{\ge x}\text{ is uncountable and totally ordered}\}$$
    is at most countable. If $x\in E$, then $S_{\ge x}$ 
    is (homeomorphic to) an $\omega_1$-compact and hence stationary subset of $\omega_1$.
    Recall that hereditary Lindel\"ofness and extent are equal in a metrizable space (see Lemma \ref{lemma:metriceq});
    hence, $f(S)$ is Lindel\"of and by Theorem \ref{thm:ecstat} $f$ is eventually constant on 
    $S_{\ge x}$ for each $x\in E$.
    We may thus assume by Lemma \ref{lemma:treefacts} (c) that $S^\downarrow$ is a Suslin tree. 
    We take out the set of $x$ such that $|f(S_{\ge x})| = 1$. 
    If what remains is countable, we are over, otherwise by
    Lemma \ref{lemma:treefacts} (d) we may assume that 
    for each $x\in S$, $S_{\ge x}$ is uncountable and $|f(S_{\ge x})|>1$.
    We can also assume that $S^\downarrow$ is rooted by adding a common root below its minimal elements.
    We now use a forcing argument, and force with $S^\downarrow$ with the reverse order.
    Let $G$ be a generic filter.
    Recall that if $\langle X,\tau\rangle$ is a topological space in the ground model, then $\tau$ serves as a base for 
    the topology $\tau(G)$ of $X$ on a forcing extension by $G$. Thus, any function that is continuous in
    the ground model remains so in the forcing extension, and $Y$ remains metrizable.
    Since $D_\alpha = \{ x\in S^\downarrow\,:\,\text{height}(x)\ge\alpha\}$ is dense
    (and is in the ground model), 
    $G$ is an new uncountable branch in $S^\downarrow$.
    Forcing with a Suslin tree preserves stationarity and cardinals, 
    actually it does not add countable sets to the universe,
    see for instance \cite[{Exercices VII H1--H2 and Theorems VII.5.10 \& VII.8.4}]{Kunen}.
    Hence, $\omega_1^V = \omega_1^{V[G]}$. 
    Since $S$ intersects stationary many levels of $S^\downarrow$ in $V$ by Lemma \ref{lemma:treefacts} (e)
    it does so in $V[G]$ as well. Hence, $S$ intersects stationary many levels of $G$ and $S\cap G$ is 
    homeomorphic to a stationary subset of $\omega_1$.
    By Theorem \ref{thm:ecstat}, $f$ is constant on $G$ above some height; hence,
    there must be $\alpha\in\omega_1$, $s\in S$, $u\in Y$ with 
    $s\Vdash \check{f}\left(\check{G} - \check{S}_{\le\check{\alpha}}\right) = \{\check{u}\}$.
    But since $|f(S_{\ge x})|>1$ for any $x\in S$, $\{z\in S\,:\, f(z)\not= u\}$ is order-dense in $S^\downarrow$,
    and we have that
    $1\Vdash \check{f}\left( \check{G} - \check{S}_{\le\check{\alpha}}\right)\not=\{\check{u}\}$, a contradiction.
    It follows that
    $|f(S_{\ge x})| = 1$ if $x$ is above some level $\alpha$.
    \\
    (c) Follows immediately from (a) and (b).
\endproof

One consequence of Theorem \ref{thm:Steprans} is the following.
\begin{lemma}\label{lemma:retractr_alpha}
   Let $T$ be an $\omega_1$-compact tree of height $\omega_1$ (in particular a Suslin tree), 
   $Y$ be a space and assume that $\mathsf{S}(T,Y)$ holds.
   Then, given $f:T\to Y$, the retraction witnessing $\mathsf{S}(T,Y)$ 
   can be chosen to be $r_\alpha$ for some $\alpha$.
\end{lemma}
\proof
   Let $f:T\to Y$ be given and $r:T\to Z\subset  T_{\alpha+1}$ be a retraction 
   such that $f\circ r = f$. Since $T$ is $\omega_1$-compact, 
   by Theorem \ref{thm:Steprans} $\mathsf{S}(T)$ holds.
   Since
   $T_{\alpha+1}$ can be embedded in $\R$, there is some $\alpha$
   such that $r \circ r_\alpha = r$; hence, 
   $$ f\circ r_\alpha = (f\circ r)\circ r_\alpha = f\circ (r\circ r_\alpha)= f\circ r = f.$$
\endproof

This has the following corollary.
\begin{cor}\label{lemma:treesproduct}
   Let $Y_n$, $n\in \omega$, be spaces and $T$ be a tree of height $\omega_1$.
   Assume that $Y_n$ is uncountable for at least one $n$.
   Then $\mathsf{S}(T,Y_n)$ holds for each $n\in \omega$ iff $\mathsf{S}(T,\prod_{n\in \omega}Y_n)$ holds.
\end{cor}
\proof 
   The reverse implication is immediate. Assume that $\mathsf{S}(T,Y_n)$ holds for each $n\in \omega$.
   Let $\pi_n$ be the projection on the $n$-th factor.
   Since $\mathsf{S}(T,Y_n)$ holds and $|Y_n|\ge \aleph_1$ 
   for some $n$, by Theorem \ref{thm:Steprans} (b) $T$ is $\omega_1$-compact.
   Let $f:T\to\prod_{n\in \omega}Y_n$.
   By Lemma \ref{lemma:retractr_alpha}, for each $n\in\omega$
   there is $\alpha_n$ such that $\pi_n\circ f\circ r_{\alpha_n} = \pi_n\circ f$.
   Hence $f\circ r_\alpha = f$ for $\alpha = \sup_n \alpha_n$.
\endproof

Another simple consequence of Theorem \ref{thm:Steprans} and Lemma \ref{lemma:retractr_alpha} is the
following.

\begin{prop}\label{prop:stepransBR}
   Let $T$ be a tree of height $\omega_1$ and let $Y$ be an uncountable space. Then
   $\mathsf{S}(T,Y)\Rightarrow \mathsf{BR}(T,Y)$.
\end{prop}
\proof 
       If $\mathsf{S}(T,Y)$ holds, 
       then by Theorem \ref{thm:Steprans} (b)
       $T$ cannot have an uncountable antichain; hence, $T$ is $\omega_1$-compact.
       By Lemma \ref{lemma:treefacts} (d),
       there is $\gamma\in\omega_1$ such that $T_{\ge x}$ is uncountable if $x$ is above level $\gamma$.
       Given $f:T\to Y$, by Lemma \ref{lemma:retractr_alpha} there is some $\alpha$ such that $f\circ r_\alpha = f$.
       We can take $\alpha\ge\gamma$.
       Then, for each $\beta\ge\alpha$, $f(T-T_{<\beta}) = f(T-T_{<\alpha}) = f(\{x\,:\,\text{height}(x)= \alpha\})$.
       This shows that $\mathsf{BR}(T,Y)$ holds.
\endproof

Let us now briefly see what happens when the codomain is countable.

\begin{thm}\label{thm:stepranscountable}
   Let $S$ be an uncountable $\omega_1$-compact subset of a tree of height $\omega_1$ and $Y$ 
   be a countable $T_1$ (non necessarily Hausdorff) 
   space, then $\mathsf{S}(T,Y)$ holds. 
\end{thm}
\proof
   Let $f:S\to Y$ be given. As before we work in $S^\downarrow$.
   By Lemma \ref{lemma:treefacts} (f), $|(f^{-1}(\{y\}))^\downarrow \cap (f^{-1}(\{z\}))^\downarrow |\le\omega$ 
   for each distinct $y,z\in Y$.
   Since $Y$ is countable, there is some height $\alpha$ such that $S_{\ge x}$ 
   intersects the preimage of at most (and thus exactly) one point of $Y$
   if $\text{height}(x)\ge\alpha$.
\endproof

This theorem does not hold for every $\omega_1$-tree, and in fact,
Proposition \ref{prop:stepransBR} is false when $Y=\Q$.
Recall that a tree $T$ is {\em special} iff it is a countable union of antichains, or 
equivalently iff there is an order preserving (not necessarily continuous) map $f:T\to\Q$.
If there is a continuous such $f$, we say that 
$T$ is {\em c-special}. 
C-special Aronzsajn trees exist in {\bf ZFC} (for instance included in the compact
subsets of $\Q$ ordered by end-extension, see e.g. \cite[Exercice 14.32]{JustWeeseII}).
Notice that under {\bf MA + $\neg$CH}, all Aronszajn trees are c-special
and special non c-special Aronszajn trees exist under $\diamondsuit$
(see e.g. \cite{KunenLarsonSteprans:2012}).

\begin{lemma}
   \label{ex:BRnotS} 
   Let $T$ be a c-special Aronzsajn tree. 
   Then $\mathsf{BR}(T,\Q)$ and $\mathsf{L}(T,\Q)$ hold but $\mathsf{S}(T,\Q)$ does not. 
\end{lemma}
In the proof below, we denote $\{x\}^\downarrow$ by $x^\downarrow$ for $x\in T$.
We choosed an argument which shows that the weaker property $w\mathsf{S}(T,\Q)$ 
defined later (Definition \ref{def:wSaS}) does not hold either.
\begin{proof}
   $\mathsf{BR}(T,\Q)$ and $\mathsf{L}(T,\Q)$ hold by Lemma \ref{lemma:QBRcountable}.
   Fix some countinous strictly increasing $f:T\to\Q$, and let $\alpha\in\omega_1$.
   Suppose that there is some $\beta\in\omega_1$
   and a continuous $r:T\to T_{\le\alpha}$ such that $f\circ r(x) = f(x)$ whenever $x$ is above level $\beta$.
   (In particular, $r:T\to T_{\le\alpha}$ might be a retraction, in that case $\beta=\alpha$.)
   Since $f$ is strictly increasing, $f\upharpoonright x^\downarrow$ is $1$-to-$1$ for any $x\in T$,
   hence $r\upharpoonright (x^\downarrow - T_{\le\beta})$ must also be $1$-to-$1$.
   Moreover, by continuity of $r$, for each $x$ at a limit level there is $y(x)<x$ such that 
   $r(z)<r(x)$ whenever $y(x)<z<x$. 
   Thus, if $x$ is above level $\max(\alpha,\beta)$, then 
   $r\upharpoonright (x^\downarrow - y(x)^\downarrow)$
   is an order preserving embedding.
   The pressing-down lemma for $\omega_1$-trees 
   yields $y\in T$ such that $E = \{x\in T\,:\,y(x) = y\}$ meets stationary many levels, let $\gamma = \text{height}(y)$.
   We may assume that $\gamma\ge\max(\alpha,\beta)$.
   Take $x\in E$ whose height is strictly more than $\gamma + \alpha$ (ordinal addition).
   Then $r\upharpoonright (x^\downarrow - y^\downarrow)$ is an embedding from a well order type $>\alpha$
   into one $\le\alpha$, a contradiction.
   Hence, there is no such $r$, and $\mathsf{S}(T,\Q)$ does not hold.
\end{proof}

Going back to Suslin trees and Theorem \ref{thm:Steprans} (a), it is not 
clear whether it is possible to weaken the assumption about submetrizabily of $Y$.
For instance, notice the following.

\begin{example}
   \label{ex:suslinline}
   Let $T$ be a Suslin tree such that each element has infinitely many immediate successors.
   Then, there is a a $1$-to-$1$ continuous $f:T\to Y$ where $Y$ is a 
   hereditarily Lindel\"of first countable monotonically normal space.
   Hence, $\mathsf{L}(T,Y)$ and $\mathsf{BR}(T,Y)$ do not hold.
\end{example}
\proof[Details]
   $Y$ is actually the tree $T$ itself with the topology given by a lexicographical ordering of 
   $T$ chosen so that the resulting space is a Suslin line (see just below).
   The lexicographical order $\le_\ell$ in $T$ is obtained by first considering a total order $\prec$ on $T$
   and then letting $y<_\ell x$ if and only if $y< x$ or $y\perp x$ and $y\upharpoonright\alpha \prec  x\upharpoonright\alpha$,
   where $\alpha$ is minimal such that $y\upharpoonright\alpha \not= x\upharpoonright\alpha$.
   Denote the topologies induced by $<,<_\ell$ on $T$ by $\tau_<,\tau_{<_\ell}$.
   The $<$-minimal elements of an
   $<_\ell$ interval $\{y\in T\,:\,x<_\ell y<_\ell z\}$ cannot be 
   at limit height because the tree is Hausdorff. It follows that any such $<_\ell$-interval is 
   a union of branches starting at a successor level and is thus
   $\tau_<$-open in $T$.
   Hence $id:\langle T, \tau_<\rangle\to \langle T, \tau_{<_\ell}\rangle$ is continuous.
   It is well known that if $\prec$ orders the immediate successors of 
   any given point as $\Q$, then $\langle T,<_\ell\rangle$ is a Suslin line,
   see for instance \cite[Lemma 14.21]{JustWeeseII} for a proof.
   We recall that a Suslin line is a linearly ordered topological space which is ccc but not separable.
   Every linearly ordered space is monotonically normal and every Suslin line is first countable and hereditarily Lindel\"of.
\endproof

Perhaps having a small or a $G_\delta$-diagonal is enough for $Y$ is Theorem \ref{thm:Steprans} (a).
It is well known that Suslin lines do not have these properties, see \cite[Proposition 4.4]{BennettLutzer:1997}
for definitions and statements.
Example \ref{ex:suslinline} shows in particular that the variants of 
Theorems \ref{thm1EC} and \ref{thm:ecstat} with $\mathsf{S}$ instead of $\mathsf{EC}$
do not hold for Suslin trees. But we can still say something
if the conclusion of Lemma \ref{lemma:cunbd} holds.

\begin{lemma}\label{lemma:suslinpreim}
   Let $Y$ be a regular space with $G_\delta$ points such that for each Suslin tree $T$ and each $f:T\to Y$, there
   is $y\in Y$ with an uncountable preimage. Then $\mathsf{S}(T,Y)$ holds for any Suslin tree $T$.
\end{lemma}
\proof
  Let $T$ be a Suslin tree.
  By Lemma \ref{lemma:treefacts} (d) we may assume that $T_{\ge x}$ is uncountable for each $x\in T$.
  By Lemma \ref{lemma:treefacts} (h), it is enough to show that 
  $E = \{x\in T\,:\,\exists y\in Y\text{ with }f(T_{\ge x})=\{y\}\}$ is order-dense in $T$.
  Let thus $z\in T$. Since $T_{\ge z}$ is Suslin there is $y$ such that $A = f^{-1}(\{y\})$ is uncountable.
  Hence, by Lemma \ref{lemma:treefacts} (g) there is $w\in T_{\ge z}$ such that $T_{\ge w}\subset A^\downarrow$.
  Since $A$ is closed, by Lemma \ref{lemma:treefacts} (e) it intersects a stationary subset of levels of 
  $T_{\ge u}$ for each $u\ge w$.
  Let $U_n\subset Y$ be open such that $\cap_{n\in\omega}U_n=\{y\}$.
  Then each $f^{-1}(U_n)$ is open, order-dense and intersects stationary many levels of $T_{\ge u}$ for each $u\ge w$.
  By Lemma \ref{lemma:treefacts} (i), $A = f^{-1}(\{y\}) = \cap_{n\in\omega}f^{-1}(U_n)$ contains all of 
  $T_{\ge w}$ above some level. This shows that $E$ is order-dense.
\endproof

We do not know whether this lemma is of any use, that is, if there are spaces satisfying its hypotheses but 
not those of Theorems \ref{thm:Steprans}.


\section{$\mathsf{BR}(X,Y)$, $\mathsf{L}(X,Y)$ when $Y$ is metric and $X$ an increasing chain of `nice' subspaces}
\label{sec:BRL}

The goal of this section is to prove in particular Theorem \ref{thm:manifoldsL} below. 
We recall that manifolds are assumed to be connected.

\begin{defi}
   A space $X$ is said to be d-connected iff it is path connected and
   for each countable closed discrete subset
   $\{x_n\,:\,n\in\omega\}\subset X$, there is an infinite subsequence $\{x_{n_k}\,:\,k\in\omega\}$
   and $f:[0,+\infty)\to X$ with closed image such that $f(k) = x_{n_k}$ for each $k$.
\end{defi}

While it is true that any countable subset of a manifold can be put in one chart 
(see e.g. \cite[Corollary 3.4]{GauldBook}), some manifolds are not d-connected
as shown by Nyikos in \cite[Example 6.7]{Nyikos:1990}.
Of course, any countably compact manifold is vacuously d-connected, as it contains no closed discrete infinite subset.

\begin{thm}\label{thm:manifoldsL}
   Let $M$ be a manifold. Then $\mathsf{L}(M)$ holds.
   If $M$ is moreover d-connected, then $\mathsf{L_{cl}}(M)$ holds as well.
\end{thm}

The ``moreover'' part has actually a simple proof.
(Recall that manifolds are path-connected (see e.g. \cite[Corollary 3.4]{GauldBook}).)

\begin{lemma}
   \label{lemma:d-connected}
   Let $X$ be d-connected. Then $\mathsf{L_{cl}}(X)$ holds.
\end{lemma}
\begin{proof} 
   Let $f:X\to\R$ be given, we may assume that
   $f(X)\subset[0,1]$. By (path) connectedness of $X$, $f(X)\subset\R$ is an interval 
   (perhaps reduced to a point) $I$, denote its endpoints
   by $a,b$ with $a\le b$.
   If $a=b$, then $f$ is constant hence $f(X) = f(\{x\})$ for any $x\in X$.
   If not, choose $c\in (a,b)$.
   Suppose first that $a\in I$. Let $g_a:[0,1]\to X$ be such that $f(g(0))=a, f(g(1)) = c$.
   By connectedness again, $f(g([0,1])) = [a,c]$, and $G_a = g([0,1])$ is compact and hence closed Lindel\"of.
   Suppose now that $a\not\in I$ and let $a_n\searrow a$ be a strictly decreasing sequence in $I$ converging to $a$.
   Choose $x_n\in f^{-1}(\{a_n\})$, then the sequence $x_n$ must be discrete since $a\not\in I$.
   Let thus $g:[0,+\infty)\to X$ be such that $G_a = g([0,+\infty))$ is closed 
   and contains infinitely many $x_n$'s. Notice that $G_a$ is automatically Lindel\"of.
   We may assume that $c\in G_a$ (if not, take some $x_n\in G_a$, a path from $x_n$ to $c$ and argue as
   in the case $a\in I$).
   Then $f(G_a)\supset (a,c]$. \\
   In both cases, we found a closed Lindel\"of $G_a$ with $f(G_a)\supset [a,c]\cap I$.
   Proceed the same with $b$ to find a closed Lindel\"of $G_b\subset X$ whose image contains $[c,b]\cap I$.
   Then $G_a\cup G_b$ does the job.
\end{proof}

If $M$ is not d-connected,
the proof relies on properties of increasing subsets of $\R$, and more generally of metric spaces. 
We start by specifying some notation.
Let $\mathscr{P}$ be a topological property and $\gamma$ be an ordinal.
We say that a space $X$ is a {\em $\mathscr{P}$-chain of length $\gamma$} iff 
$X=\cup_{\alpha<\gamma} H_\alpha$ with $H_\alpha\subset H_\beta$ if $\alpha<\beta$ and such
that $H_\alpha$ has property $\mathscr{P}$ for each $\alpha$. 
If $H_\alpha\subsetneq H_\beta$ for each $\alpha<\beta$, we say that the chain is strict.
If the property $\mathscr{P}$ cannot be described in one word, we may use brackets. For instance,
Type I spaces are strict [open-with-Lindel\"of-closure]-chains of length $\omega_1$.
A subspace $E\subset X$ such that $X-E$ satisfies $\mathscr{P}$ is said to be co-$\mathscr{P}$.
We will use the following classical fact.

\begin{lemma}\label{lemma:metricchain}
   Let $Y$ be a metric space. \\
   (a) $Y$ does not contain a strict compact-chain or a strict [co-compact]-chain of length $\ge\omega_1$.\\
   (b) If $Y$ is separable, then $Y$ does contain neither a strict open-chain nor a strict closed-chain of length $\ge\omega_1$.
\end{lemma}
We provide a proof for completeness. Lemma \ref{lemma:metriceq} is used several times implicitely.
\proof[Proof of Lemma \ref{lemma:metricchain}]
We show (b) first.
\\
(b) $Y$ is hereditarily Lindel\"of and hereditarily separable.
   A strict open-chain of length $\omega_1$ is not 
   Lindel\"of and a strict closed-chain of length $\omega_1$ contains a non-separable subspace
   (see the proof of Theorem 3.1 in \cite{Roitman:1984}). \\
(a) If $\langle H_\alpha\,:\,\alpha<\omega_1\rangle$, is a compact-chain of subspaces
    of $Y$, then $H=\cup_{\alpha<\omega_1} H_\alpha$
    is countably compact and hence compact. Indeed, any countable subset of $H$ is contained in some compact $H_\alpha$
    and has thus an accumulation point. Hence $H$ is separable, let $D$ be a countable dense subset.
    Then $D$ is contained in some $H_\alpha$, by closedness
    $H = \wb{D} = H_\alpha$ and the chain cannot be strict.
    If $\langle H_\alpha\,:\,\alpha<\omega_1\rangle$ is a chain of co-compact subspaces of $Y$,
    then $H_\alpha\cap W$ is open in $W = Y-H_0$ which is compact and hence separable. 
    By (b) the open chain $\langle H_\alpha\cap W \,:\,\alpha<\omega_1\rangle$ in $W$ cannot be strict.
    Hence the chain of $H_\alpha$ is not strict either.
\endproof

The following corollary is immediate.
\begin{cor}\label{cor:intervals}
   $\R$ does contain neither a strict connected-chain nor a strict [co-connected]-chain of length $\omega_1$.
\end{cor}

We say that $X$ is a {\em slowpen chain} if is a strict open-chain $X=\cup_{\alpha\in\gamma}U_\alpha$, such that
$\cup_{\beta<\alpha}U_\beta = U_\alpha$ when $\alpha$ is limit,
$U_0=\varnothing$ and $U_{\alpha+1}=U_\alpha\cup E_\alpha$, 
where $E_\alpha$ is Lindel\"of for each $\alpha$.

\begin{lemma}\label{lemma:slowpenchain}
   Let $X=\cup_{\alpha<\gamma}U_\alpha$ be a slowpen chain. 
   If $U_\alpha$ is connected for each $\alpha$, then $\mathsf{L}(X)$ holds.
\end{lemma}
\proof
   If $\gamma<\omega_1$, $X$ is Lindel\"of and there is nothing to prove. Assume that $\gamma\ge\omega_1$.
   Let $f:X\to\R$ be given. Let $U_\alpha,E_\alpha$, $\alpha<\gamma$ be as in the definition of a slowpen chain, hence
   $X = \cup_{\alpha<\gamma}U_\alpha$.
   We prove by induction on $\alpha$ that there is a Lindel\"of subset $L_\alpha\subset X$ such that
   $f(L_\alpha) = f(U_\alpha)$. This yields the theorem when $\alpha = \gamma+1$.\\
   If $\alpha$ is countable, $U_\alpha$ is Lindel\"of hence we may set $L_\alpha =U_\alpha$.
   If $\alpha = \beta + 1$, then the Lindel\"of subset $L_\alpha = L_\beta \cup E_\beta$ satisfies
   $f(U_\alpha) = f(L_\alpha)$.
   If $\text{cof}(\alpha) = \omega$, choose an 
   $\omega$-sequence $\beta_n\nearrow\alpha$ and set $L_\alpha =\cup_{n\in\omega}L_{\beta_n}$.
   Then $L_\alpha$ is Lindel\"of and $f(U_\alpha) = f(L_\alpha)$.
   Assume now that $\text{cof}(\alpha) > \omega$. 
   Since $U_\beta$ is connected for each $\beta<\alpha$, 
   $f(U_\beta)\subset\R$ is connected. Hence, $\langle f(U_\beta)\,:\,\beta<\alpha\rangle$ 
   does not have a cofinal strict subchain by Corollary \ref{cor:intervals}.
   It follows that there is some $\beta<\alpha$ such that $f(U_\alpha) = f(U_\beta) = f(L_\beta)$. 
\endproof

\begin{lemma}\label{lemma:manchain}
  Let $M$ be a manifold. Then $M$ is a 
  slowpen chain $\cup_{\alpha<\gamma}U_\alpha$ of length $\gamma < \mathfrak{c}^+$ with each $U_\alpha$ connected.
  Moreover, each $U_\alpha$ is a submanifold of $M$ and $M$ is metrizable iff $\gamma$ is countable.
\end{lemma}
\proof
  Our argument is very similar to the proof of Theorem 2.9 
  in \cite{Nyikos:1984} which shows that the cardinality of a manifold is the continuum.
  Set $U_0=\varnothing$ and $U_1$ be any open 
  subset homeomorphic to $\R^n$. We proceed by induction and assume that the chain of connected open sets
  $\cup_{\beta<\alpha}U_\beta$ is a slowpen chain.
  If $\alpha$ is limit, we let $U_\alpha = \cup_{\beta<\alpha}U_\beta$. 
  Given $U_\alpha$, if
  $U_\alpha = M$, we are over. If not, by connectedness of $M$ we have that
  $\wb{U_\alpha}\not= U_\alpha$, so 
  we may choose a point $x \in \wb{U_\alpha}- U_\alpha$. 
  We then let $U_{\alpha+1}$ be the union of $U_\alpha$ and an Euclidean connected 
  set $E_\alpha$ containing $x$. 
  Connectedness of $U_\alpha$ for each $\alpha$ is immediate. By construction the chain is slowpen and strict up to $\alpha$.
  Since $U_\alpha$ is open, it is a submanifold of $M$.
  Since $|M| = \mathfrak{c}$, there is some $\gamma < \mathfrak{c}^+$ such that $M = U_\gamma$.
  If $\gamma$ is countable, $M$ is Lindel\"of and hence metrizable. 
  (Recall that (hereditary) Lindel\"ofness and metrizability are equivalent for manifolds,
  see e.g. \cite[Thm 2.1]{GauldBook}.)
  If $\gamma$ is uncountable, it contains the non-Lindel\"of subset
  $\cup_{\alpha<\omega_1}U_\alpha$, which shows that $M$ is not metrizable. 
\endproof

Theorem \ref{thm:manifoldsL} follows immediately by Lemmas \ref{lemma:slowpenchain}--\ref{lemma:manchain}.
Notice however 
that the codomain $\R$ cannot be replaced by $\R^2$, as shown by Example \ref{ex:QRR2} below. 
Also, we cannot replace $\mathsf{L}$ by $\mathsf{BR}$.
\begin{example}\label{ex:prufer}
   There is a separable surface (without boundary) $Q$ such that $\mathsf{L_{cl}}(Q)$ holds while 
   $\mathsf{BR}(Q)$ does not.
\end{example}
\proof[Details]
This example is classical: Take first $P$ to be the Pr\"ufer surface (separable version with boundary) 
$H_0\cup\displaystyle\cup_{a\in\R}\R_a$, where $H_0$ is a copy of $\R\times\R_{>0}$ and each
$\R_a$ is a (distinct) copy of $\R$. See \cite[Ex. 1.25]{GauldBook} for a complete description of the topology. 
This surface can be seen as taking $\R\times\R_{\ge 0}$ and `blowing up' each point 
$\langle a,0\rangle$ in the bottom boundary into the open interval $\R_a$.
This gives again a surface whose boundary components are the $\R_a$, and
$\displaystyle H_0\cup\cup_{a\in A}\R_a$ is open for any $A\subset\R$,
hence a Lindel\"of subset of
$P$ intersects at most countably many $\R_a$.
Any continuous $f:\R\times\R_{\ge 0}\to\R$ yields a continuous 
$f_P:P\to\R$ by letting $f_P$ be constant with value $f(\langle a,0\rangle)$ on $\R_a$ and
equal to $f$ on $H_0$.
Let $\pi:\R\times\R_{\ge 0}\to\R$ be the projection on the first coordinate. 
Then $\pi_P$ shows that $\mathsf{BR}(P)$ does not hold. Indeed, if $Z\subset P$ is Lindel\"of, 
choose $a\in\R$ such that $\R_a\cap Z = \varnothing$ and
set $W = Z\cup \R_a \cup H_0$. Then $W$ is Lindel\"of and $\pi_P(P-Z)\ni a\not\in \pi_P(P-W)$.
We can obtain a boundaryless version $Q$ by 
taking the double of $P$ (two copies with the boundaries identified pointwise),
extending $\pi_P$ the obvious way.
Hence, we have $Q=H_0^0\cup H_0^1\displaystyle\cup_{a\in\R}\R_a$, where $H_0^0$ and $H_0^1$
are two copies of $H_0$.
 \\
To see that $\mathsf{L_{cl}}(Q)$ holds, we show that 
$Q$ is d-connected and apply Lemma \ref{lemma:d-connected}. The proof is easy but a bit 
tedious, so we skip some details.
It might be useful to
recall that a basic neighborhood of a point in $\R_a$ in $P$ is given by 
an interval $(u,v)$ in $\R_a$ union a triangle in $H_0$
situated between two lines of slopes $u,v$ pointing at $\langle a,0\rangle$ and an horizontal line
at height $\epsilon>0$ (see again \cite[Ex. 1.25]{GauldBook}).
In particular, the closure of a vertical segment $\{a\}\times(0,b]\subset H_0$ 
is itself plus the $0$-point $0_a$ of $\R_a$. 
Let $Y = \{x_n\,:\,n\in\omega\}$ be a closed discrete
sequence in $Q$. 
If its intersection with $H_0^i$ is infinite for one $i\in\{0,1\}$, 
then there is a closed discrete subsequence $Y_0$ which is 
monotone in both coordinates. Since $Y_0$ is closed discrete,
joining each consecutive member 
by a straight segment yields a copy of $[0,+\infty)$ inside $H_0^i$ which is closed in $Q$.
(Details are tedious to describe, but the reader can draw some triangles and straight segments
to be convinced.)
If $H_0^i\cap Y$ is finite for $i=0,1$, 
then either $Y$ has infinite intersection with one of the $\R_a$ 
or its intersection with each $\R_a$ is finite for each $a$.
In the former case, there is a subsequence of $Y$ shich is stricly monotone in $\R_a$ 
and we just join the consecutive points of this subsequence.
In the latter case, there is a strictly monotone sequence $a_n\in\R$ such that $Y$ intersects
each $\R_{a_n}$. We may assume that $a_n$ is strictly increasing and that $x_n\in C_{a_n}$.
Choose a smooth
function $c:[a_0,\sup_n a_n)\to \R_{>0}$ that converges to $0$ when $x\to \sup_n a_n$,
with derivative going to $0$ as well. The graph of $c$ in $H_0^0$ is a closed copy of $[0,+\infty)$ in $Q$.
This is obvious if $\sup_n a_n=+\infty$, otherwise notice that any triangle pointing at $\langle a,0\rangle$
(with $a=\sup_n$)
in $H_0^0$ whose height is small enough avoids the graph of $c$.
We may then define the path joining $x_n\in C_{a_n}$ for successive $n$'s as follows:\\
(1) Take a segment joining $x_n$ to $0_{a_n}\in\R_{a_n}$;\\
(2) add the vertical segment $\{a_n\}\times(0,c(a_n))\subset H_0^0$ which joins $0_{a_n}$ to $\langle a_n, c(a_n)\rangle$\\
(3) follow the graph of $c$ until $\langle a_{n+1}, c(a_{n+1})\rangle$;\\
(4) Repeat (2) and (1) backwards with $n+1$ instead of $n$.\\
This gives the required path between each $x_n$.
\endproof

Let us now look further than manifolds and see whether we can say something about 
the interplay between $\mathsf{L}(X)$ and $\mathsf{BR}(X)$ for
more general spaces which are still ``nice'' chains. First, notice that
the following proposition can be proved exactly by the same argument
seen in the proof of Lemma \ref{lemma:slowpenchain}.

\begin{prop}\label{prop:almost_stagnation_BR}
   Let $X=\cup_{\alpha\in\omega_1}H_\alpha$ 
   be a strict Lindel\"of-chain (for instance, a Type I space). \\
   (a) If $H_\alpha$ is connected for uncountably many $\alpha$ then $\mathsf{L}(X)$ holds.
       If $H_\alpha$ is moreover closed for uncountably many $\alpha$, then $\mathsf{L_{cl}}(X)$ holds.\\
   (b) If $X-H_\alpha$ is connected for uncountably many $\alpha$ then $\mathsf{BR}(X)$ holds.
       If $H_\alpha$ is moreover closed for uncountably many $\alpha$, then $\mathsf{BR_{cl}}(X)$ holds.
\end{prop}

In view of our results so far, we may wonder whether connectedness alone is 
enough to ensure that a space is a lazy explorer of, and/or a broken record in $\R$.
This is not the case, and there is a simple example very similar to Example \ref{ex:pasAS}.
Given a cardinal $\kappa$, the $\kappa$-long ray $\LL_\kappa$ is defined as $\kappa\times[0,1)$
with lexicographic order (hence, $\LL_{\ge 0} = \LL_{\omega_1}$). 
Again, we see the ordinal $\kappa$ as a subset of $\LL_\kappa$ by identifying $\alpha$ with $\langle\alpha,0\rangle$.
One sees easily that $\LL_\kappa$ is connected and that for any
$\alpha<\kappa$ the interval $[0,\alpha]\subset\LL_\kappa$ is compact.
Let $\{x_\alpha\,:\,\alpha\in\mathfrak{c}\}$ be a well ordered enumeration of $\R$. We define the following
subspaces of $\LL_{\mathfrak{c}}\times\R$:
\begin{align*}
   H_\mathfrak{c}^+ &= \cup_{\alpha\in\mathfrak{c}}(\alpha,\mathfrak{c})\times\{x_\alpha\},\\
   H_\mathfrak{c}^- &= \cup_{\alpha\in\mathfrak{c}}[0,\alpha]\times\{x_\alpha\}.
\end{align*}
\begin{example}
   \label{ex:connectednoL}
   Let $X$ be the disjoint union of $H_\mathfrak{c}^+$, $H_\mathfrak{c}^-$ and $[0,1]$, where $0$ is identified with 
   $\langle x_0,0\rangle\in H_\kappa^-$ and $1$ with $\langle x_0,1\rangle\in H_\kappa^+$.
   Then $X$ is connected and neither $\mathsf{L}(X)$ nor $\mathsf{BR}(X)$ do hold.
\end{example}
We shall need the following weak version of Fodor's theorem
(see e.g. \cite[Lemma 2.1]{Jech:Handbook} for a proof).
\begin{lemma}
   \label{thm:weakFodor}
   Let $\kappa$ be a cardinal with uncountable cofinality, $S\subset\kappa$ be stationary
   and $f:S\to\kappa$ be regressive, i.e. $f(\alpha)<\alpha$. Then there is a stationary $R\subset S$
   such that $f\upharpoonright R$ is bounded.
\end{lemma}

\begin{proof}[Details of Example \ref{ex:connectednoL}]
  As in Example \ref{ex:pasAS}, the projection on the second coordinate on $H_\mathfrak{c}^+$ and $H_\mathfrak{c}^-$
  (and the identity on $[0,1]$) shows that $\mathsf{L}(X)$ and $\mathsf{BR}(X)$ do not hold.
  We show that $X$ is connected. It is enought to show that both $H_\mathfrak{c}^+$ and $H_\mathfrak{c}^-$ are connected.
  This is clear for $H_\mathfrak{c}^-$ 
  because it is a union of ``rays'' of length $\alpha<\mathfrak{c}$
  all glued at their smallest point to  $\{0\}\times\R$.
  Let us show that $H_\mathfrak{c}^+$ is connected.
  Notice that the projection on the second coordinate of $H_\mathfrak{c}^+$ is onto.
  Since $\mathfrak{c}$ has uncountable cofinality, any countable subset of $\LL_{\mathfrak{c}}$ is bounded, 
  and
  by Lemma \ref{thm:weakFodor} any open set of 
  $\LL_{\mathfrak{c}}$ containing a stationary subset must contain a terminal segment.
  It follows that for any $a\in\R$, any open subset $U$ of 
  $\LL_{\mathfrak{c}}\times\R$ intersecting $\LL_{\mathfrak{c}}\times\{a\}$ in a stationary subset
  must contain $[\beta,\mathfrak{c})\times I$ for some $\beta<\mathfrak{c}$ and open interval $I\ni a$.
  Indeed, $U$ contains $[\beta,\mathfrak{c})\times \{a\}$
  for some $\beta$, if $U$ does not contain 
  $[\beta,\mathfrak{c})\times I$ for some open $I\ni a$,
  there is a sequence $\langle \beta_n,a_n\rangle\in \left([\beta,\mathfrak{c})\times[a-\frac{1}{n},a+\frac{1}{n}]\right)-U$.
  Any cluster point of this sequence is in $\left([\beta,\mathfrak{c})\times \{a\}\right)-U$, a contradiction.
  \\
  Fix non-empty open $U,V$ whose union contains $H_\mathfrak{c}^+$.
  For each $a\in\R$, one of $U,V$ must intersect $\LL_{\mathfrak{c}}\times\{a\}$ 
  in a stationary subset, and hence contain $[\beta(a),\mathfrak{c})\times I(a)$ for some open $I(a)\ni a$
  and $\beta(a)<\mathfrak{c}$.
  Take a countable cover $\{I(a)\,:\,a\in A\}$ of $\R$, and set $\beta=\sup_{a\in A}\beta(a)$.
  Then
  $\beta<\mathfrak{c}$ and $[\beta,\mathfrak{c})\times\{a\}$ is 
  contained either in $U$ or in $V$ for each $a\in\R$. 
  Since $\R$ is connected, there is some $c$ in the intersection
  of the projections of $U\cap\left( [\beta,\mathfrak{c})\times\R \right)$ and 
  $V\cap\left( [\beta,\mathfrak{c})\times\R \right)$. Let $\alpha$ be such that $x_\alpha = c$.
  This shows that $U\cap V$ contains $(\gamma,\mathfrak{c})\times\{c\}$, where $\gamma = \max(\beta,\alpha)$
  and $U\cap V$ is thus non-empty.
\end{proof}

Notice that under {\bf CH}, $X$ is a Type I space.
Finally, along the same lines as Proposition \ref{prop:almost_stagnation_BR}, we have:

\begin{thm}\label{thm:ccQ}
    Let $X$ be a space, $Y$ be a submetrizable space and $\kappa$ be an infinite cardinal.\\
    (a) If $X$ is a strict compact-chain of length $\kappa$, then $\mathsf{L_{cl}}(X,Y)$ holds.\\
    (b) If $X$ is a countably compact 
        Type I space, then both $\mathsf{BR_{cl}}(X,Y)$ and $\mathsf{L_{cl}}(X,Y)$ hold.
\end{thm}
\proof
   We may assume that $Y$ is metrizable by Corollary \ref{obviouscor}. We fix $f:X\to Y$.
   Our old friend Lemma \ref{lemma:metriceq} will again be used implicitely in what follows.
   \\
   (a)
   Let us write $X=\cup_{\alpha<\kappa}H_\alpha$ with each $H_\alpha$ compact.
   There is nothing to prove if $\kappa$ has countable cofinality because $X$ is then Lindel\"of.
   We thus assume that $\text{cf}(\kappa)$ is uncountable.
       Since $f(H_\alpha)$ is compact, Lemma \ref{lemma:metricchain} (a) implies that 
       $f(H_\alpha) = f(X)$ for some $\alpha$.
   \\
   (b) If $X$ is of Type I, then $X=\cup_{\alpha<\omega_1}X_\alpha=\cup_{\alpha<\omega_1}\wb{X_\alpha}$.
       By countable compactness $\wb{X_\alpha}$ is compact, hence
       $\mathsf{L_{cl}}(X,Y)$ follows by (a). 
       We show that $\mathsf{BR}(X,Y)$ holds, which implies $\mathsf{BR_{cl}}(X,Y)$
       by Lemma \ref{lemma:typeIcl}.
       We may assume that $Y=f(X)$, hence $Y$ is compact and thus separable. 
       Since $X_\alpha$ is open,
       $X-X_\alpha$ is countably compact, $f(X-X_\alpha)$ is compact and thus closed. 
       By Lemma \ref{lemma:metricchain} (b) the chain $U_\alpha = Y-f(X-X_\alpha)$
       cannot contain a strict cofinal chain; hence, it must stagnate above some $\alpha$. 
       It follows that $f(X-X_\alpha) = f(X-X_\beta)$ for each $\beta\ge\alpha$.
\endproof

We could have proved that $\mathsf{L_{cl}}(X,Y)$ holds in (b) with the following lemma
(since the image of $X$ in $Y$ is hereditarily separable).

\begin{lemma}\label{lemma:hersep1}
   Let $Y$ be hereditarily separable.
   If $X$ is an $\omega$-bounded space (in particular, a Type I countably compact space), 
   then $\mathsf{L_{cpt}}(X,Y)$ holds and $f(X)$ is compact for any $f:X\to Y$.  
\end{lemma}
\proof
Given $f:X\to Y$, define $E\subset X$ by taking one preimage of each point in a countable dense subset $D$ of $f(X)$. 
By $\omega$-boundedness $\wb{E}$ and thus $f(\wb{E})$ are compact, so in particular $f(\wb{E})$ is closed. 
It follows that $f(\wb{E}) = \wb{f(E)} = \wb{D} = f(X)$, proving that $f(X)$ is compact
and that $\mathsf{L_{cpt}}(X,Y)$ holds.
\endproof

Notice that we cannot weaken the assumption to `$X$ is countably compact', as shown by the next example.

\begin{example}[Folklore]
    \label{ex:ccnotQ}
    There is a countably compact space $X$ such that neither $\mathsf{L_{cl}}(X)$ nor $\mathsf{BR_{cl}}(X)$ hold.
\end{example}
\proof[Details]
    The idea is to obtain a countably compact non-Lindel\"of space with a countable dense subset of isolated points, 
    and to apply Lemma \ref{lemma:exnotcl}.
    Let $\beta\omega$ be the
    \v Cech-Stone compactification of the integers $\omega$. (The integers are given the discrete topology.)
    The closure of any infinite set in $\beta\omega$ has cardinality 
    $2^{\mathfrak{c}}$ (see e.g. \cite[Lemma 0.1]{vanMill:1984}),
    in particular if $p\in\beta\omega - \omega$,
    then $\beta\omega-\{p\}$ is a countably compact non-compact (since non-closed) subspace of $\beta\omega$
    in which the integers are dense.
\endproof

The following question probably has a positive answer, but we did not find any example.
\begin{q}\label{prob:ccnotQ}
    Is there is a countably compact space $X$ such that $\mathsf{L}(X)$ or $\mathsf{BR}(X)$ does not hold~?
\end{q}


\section{A word on products on the codomain}
\label{sec:products}

Let us have a quick look at which properties are preserved under products in the codomain.
The first observation is obvious:
\begin{lemma}
   Let $\mathsf{P}$ be any property in Definition \ref{defgen}, $X$ and $Y_j$ be non empty spaces for $j$ in some index set $J$.
   Then $\mathsf{P}(X,\prod_{j\in J} Y_j) \Rightarrow \mathsf{P}(X, Y_k)$ for each $k\in J$.
\end{lemma}
\proof 
Each $Y_k$ embeds in $\prod_{j\in J} Y_j$. Apply Lemma \ref{obviouslemma}.
\endproof

The next lemma is also almost immediate. 
\begin{lemma}\label{obviouslemma2}
   Let $X$ and $Y_n$, $n\in E$, be spaces.\\
   (a) If $E$ is countable and $\mathsf{EC}(X,Y_n)$ holds for each $n\in E$, then $\mathsf{EC}(X,\prod_{n\in E}Y_n)$ holds.\\
   (b) If $E$ is finite and $\mathsf{EC_{cl}}(X,Y_n)$ holds 
       for each $n\in E$, then $\mathsf{EC_{cl}}(X,\prod_{n\in E}Y_n)$ holds.\\
   (c) If $E$ is finite and $\mathsf{EC_{cpt}}(X,Y_n)$ holds 
       for each $n\in E$, then $\mathsf{EC_{cpt}}(X,\prod_{n\in E}Y_n)$ holds.
\end{lemma}
\proof
   (a) Given $f:X\to \prod_{n\in E}Y_n$, for each $n\in E$ there is a Lindel\"of subspace $Z_n\subset X$
       such that the projection $\pi_n\circ f:X\to Y_n$ is constant outside of $Z_n$. Hence $f$ is constant outside of the
       Lindel\"of set $\cup_{n\in E} Z_n$.
   The proof for (b) and (c) is the same.
\endproof
For $\mathsf{S}$, such a simple argument does not work even for $2$ factors, as it is not clear a priori
whether the retractions corresponding to each projection
can be combined. It is the case for $\omega_1$-trees, as seen in Corollary \ref{lemma:treesproduct} above.
\begin{q}\label{q:Sproduct}
  Are there spaces $X,Y_1,Y_2$ such that
  $\mathsf{S}(X,Y_i)$ holds for $i=1,2$ and $\mathsf{S}(X,Y_1\times Y_2)$ does not~?
\end{q}

On the other hand, $\mathsf{L}$ and $\mathsf{BR}$ behave badly under products, even when the domain space is a manifold.

\begin{example}\label{ex:QRR2}
   Let $M$ be the surface 
   $$ \LL_{\ge_0}\times [-2,1]\, - \, \cup_{\alpha\in\omega_1} \{\alpha\}\times[-1,1].$$
   Then $\mathsf{L_{cl}}(M,\R)$ and $\mathsf{BR_{cl}}(M,\R)$ hold but $\mathsf{L}(M,\R^2)$ and $\mathsf{BR}(M,\R^2)$ do not.
   Moreover, $\mathsf{S}(M,\R)$ does not hold either.
\end{example}

\begin{figure}[h]
  \begin{center}
    \epsfig{ figure = 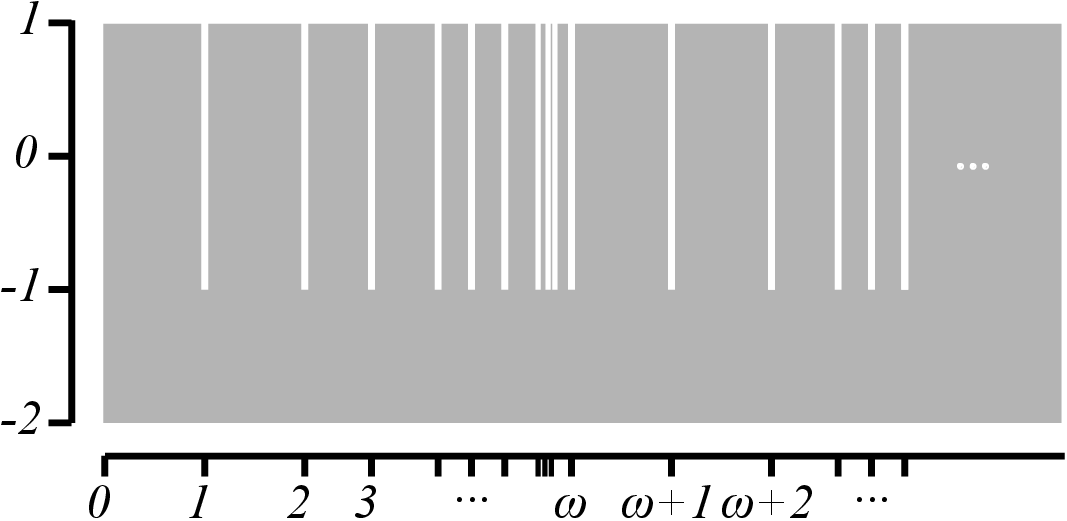, width = .6\textwidth} \\
    \caption{Example \ref{ex:QRR2}.}
    \label{fig:QRR2}
  \end{center}
\end{figure}

(We do not use the short versions $\mathsf{L_{cl}}(M)$, etc, for emphasis on the codomain.)

\proof[Details]
   Since $M$ is a manifold, $\mathsf{L}(M,\R)$ holds by Theorem \ref{thm:manifoldsL} above
   and $\mathsf{L_{cl}}(M,\R)$ follows since $M$ is Type I.
   Define $M_\alpha$ to be the subset of points of $M$ with first coordinate $<\alpha$.
   Since $M-M_\alpha$ is connected for each $\alpha$, 
   $\mathsf{BR}(M,\R)$ (hence, $\mathsf{BR_{cl}}(M,\R)$)
   holds by Proposition \ref{prop:almost_stagnation_BR} (b).
   To show that $\mathsf{L}(M,\R^2)$ and $\mathsf{BR}(M,\R^2)$ do not hold,
   let $\{e_\alpha\,:\,\alpha\in\omega_1\}$ (all distinct) be a subset of the unit circle centered at the origin in $\R^2$.
   Let $I_\alpha$ be the line segment joining the origin and $e_\alpha$.
   Define $f:M\to\R^2$ as follows. If $x\in \LL_{\ge 0}\times[-2,0]\cap M$, $f(x)$ is the origin.
   If $x = \langle u,t\rangle$ with $\alpha<u<\alpha+1$ and $t\in(0,1]$, then $f(x)$ is the point of $I_\alpha$ at distance $t$ from the origin.
   It is easy to check that $f$ is continuous and violates $\mathsf{L}(M,\R^2)$ and $\mathsf{BR}(M,\R^2)$ since 
   the preimage of $I_\alpha$ minus the origin is $(\alpha,\alpha+1)\times(0,1]$. \\
   To finish, observe that $\cup_{\alpha<\omega_1}(\alpha,\alpha+1)\times (0,1)$
   is a discrete collection of open subsets of $M$ (see the first paragraph of Section \ref{sec:SscwH} for a reminder of the definition).
   We show in Section \ref{sec:SscwH} (Lemma \ref{lemma:discretenoS}) that this implies that $\mathsf{S}(M,\R)$
   (and the weaker property $w\mathsf{S}(M,\R)$ defined in Section \ref{sec:lackretr}) does not hold.
\endproof


\section{Type I manifolds: when $\mathsf{EC}\not\Rightarrow\mathsf{S}$.}
\label{sec:lackretr}

As written before, $\mathsf{EC}$ has the looks of a being stronger property than $\mathsf{S}$, 
as being eventually constant is a rather strong assumption for a map. 
But $\mathsf{S}$ asks for retractions, and some spaces do lack them on smaller subspaces.
Actually, even a weaker property (almost stagnation, defined just next) 
might fail while $\mathsf{EC}$ holds.
Let us define two new ways to procrastinate.

\begin{defi}\label{def:wSaS}
   Given spaces $X,Y$, we say that $X$ almost stagnates [resp. weakly stagnates] in $Y$, 
   written $a\mathsf{S}(X,Y)$ [resp. $w\mathsf{S}(X,Y)$], iff
   for each $f:X\to Y$ there is a Lindel\"of $Z\subset X$ and $r:X\to Z$ such that $f\circ r = f$
   [resp. $f\circ r \upharpoonright (X-Z) = f \upharpoonright (X-Z)$].
\end{defi}

Again, $a\mathsf{S}(X)$, $w\mathsf{S}(X)$ are shorthands for $a\mathsf{S}(X,\R)$, $w\mathsf{S}(X,\R)$.
We could have chosen $a\mathsf{S}$ as the ``official'' 
definition for stagnation instead of the one with retractions, but all our counterexamples
in this and the next section are actually counter-examples to $a\mathsf{S}$ and not only to $\mathsf{S}$,
and Suslin trees in Section \ref{sec:Suslin} do satisfy the stronger property.
The next lemma is immediate from the definitions (we again abbreviate $\forall X,Y \,\mathsf{P}(X,Y)$ simply
by $\mathsf{P}$)
\begin{lemma}\label{lemma:arrows2}
   \begin{tikzcd}
  \mathsf{EC} 
    \arrow[r, rightarrow]
   & 
   w\mathsf{S} 
     \arrow[d, rightarrow]
     \arrow[r, leftarrow]
   &
   a\mathsf{S}
     \arrow[r, leftarrow]
   & \mathsf{S}
   \\
   & \mathsf{L}
   \end{tikzcd}
\end{lemma}

It is natural to ask whether some of these arrows reverse. We already know that $\mathsf{S}\not\Rightarrow\mathsf{EC}$;
hence, $w\mathsf{S}\not\Rightarrow\mathsf{EC}$ as well.  
The main result of this section implies that $\mathsf{EC}(M)\not\Rightarrow a\mathsf{S}(M)$
for $M$ a manifold, and hence 
$w\mathsf{S}\not\Rightarrow a\mathsf{S}$ as well.
As a small warm up, 
let us show that some manifolds do lack retractions. This is weaker than what we really need for
the main theorem of this section, but 
the proof contains one idea that will use again, which is contained in the following lemma.
It might be interesting to note that this is the only result from basic
algebraic topology that we need.

\begin{lemma}\label{lemma:trivialhomotopy}
   Let $\mathbb{S}^1 = [0,1]/0\sim 1$ be the circle viewed as the interval with identified endpoints.
   Let $C = [a,b]\times\mathbb{S}^1$ be a cylinder, $\pi:C\to\mathbb{S}^1$ be the projection on the second coordinate,
   and $d_0,d_1:\mathbb{S}^1\to C$ be such that $\pi\circ d_0 = \pi\circ d_1 = id_{\mathbb{S}^1}$.
   Let 
   $I_0=[0,\frac{1}{3}]$, $I_1=[\frac{1}{3},\frac{2}{3}]$ and $I_2=[\frac{2}{3},1]$,
   seen as subsets of $\mathbb{S}^1$.
   Let $r:C\to C$ be given and set $s_i = \pi\circ r \circ d_i:\mathbb{S}^1\to\mathbb{S}^1$, $i=0,1$.
   If $s_0(I_k)\subset I_k$ for $k=0,1,2$, 
   then $s_1$ is not (homotopic to) a constant map.
\end{lemma} 
\begin{proof}
   The assumptions imply that $s_0(I_k) = I_k$. 
   Sneak into the first few classes of a course on the fundamental group to 
   obtain the definition of an homotopy,
   realize that
   $s_0$ is homotopic to the identity of $\mathbb{S}^1$, and conclude that the lemma holds.
\end{proof}

\begin{prop}
   \label{thmnoret0}
   If $X=\cup_{\alpha\in\omega_1}X_\alpha$ is a longpipe such that $\mathsf{EC}(X)$ holds, then
   $X$ does not retract on any $\wb{X_\alpha}$.
\end{prop}
\proof
   In this proof, ``boundary'' means ``topological boundary''.
   First, notice that $\mathsf{EC}(X,\R^2)$ holds by Lemma \ref{obviouslemma2}.
   Suppose that $r:X\to \wb{X_\alpha}$ is a retraction. Since $\wb{X_\alpha}$ embeds in $\R^2$, 
   for some $\gamma$ we have
   $r(X-X_{\gamma})=\{x\}$. 
   Since $r$ is the identity on $\wb{X_\alpha}$, $\gamma>\alpha$.
   The boundary of $\wb{X_{\alpha+1}}$ is homeomorphic to the circle (it might not be true for $\wb{X_{\alpha}}$).
   Fix any successor $0<\beta\le\alpha$.
   The boundary of $\wb{X_{\beta}}$ is also homeomorphic to the circle.
   Take a homeomorphism $\wb{X_{\gamma+1}}\to[0,3]\times\mathbb{S}^1$ that sends $\wb{X_{\beta}}$ and
   $\wb{X_{\gamma+1}}-X_{\gamma}$ respectively 
   to $C_1=[0,1]\times\mathbb{S}^1$ and $C_2=[2,3]\times\mathbb{S}^1$.
   This yields a retraction of the cylinder $[0,3]\times\mathbb{S}^1$ into itself which is the identity on $C_1$ and 
   such that $C_2$ is sent to a point. But this is impossible by Lemma \ref{lemma:trivialhomotopy}.
\endproof

The fact that the longpipes are `cylinders piled up' is important in this proposition.
Each $X_{\alpha+1}$ in a longpipe is homeomorphic to $[0,1)\times\mathbb{S}^1$, $X$ is thus
a surface whose manifold boundary is homeomorphic to $\mathbb{S}^1$.
By sewing a disc in this `hole', we obtain what we call a
{\em sealed longpipe}. Such a sealed longpipe $W$ has each $\wb{W_{\alpha+1}}$ homeomorphic to the closed $2$-disc
and its topological boundary in any $W_\beta$ for higher $\beta$ is homeomorphic to the circle.
(Hence, its topological boundary is equal to its (sub)manifold boundary.)
\begin{lemma}\label{lemma:sealedret}
   Let $X$ be a sealed longpipe, then $X$ retracts on $\wb{X_{\alpha+1}}$ $\forall\alpha\in\omega_1$.
\end{lemma}
In the remaining of this section, we let $B(a)$ denote the closed disk of radius $a$ centered at the origin in $\R^2$. 
   \proof
      There is a homeomorphism $\wb{X_{\alpha+2}}\to B(2)$ 
      that sends $\wb{X_{\alpha+1}}$ to $B(1)$.
      Take a retraction of $B(2)$ on $B(1)$ which sends all the (topological) boundary of $B(2)$ to the origin.
      This yields a retraction $r:\wb{X_{\alpha+2}}\to\wb{X_{\alpha+1}}$ 
      such that the topological boundary of $\wb{X_{\alpha+2}}$ is sent to a point $x\in X_{\alpha+1}$.
      Extend it to the whole $X$ by $r(y)=x$ for all $y\not\in\wb{X_{\alpha+2}}$.
   \endproof

However, the fact that there are retractions `as high as one wants' 
does not ensure that $\mathsf{EC}$ implies $a\mathsf{S}$.
Actually, quite the opposite is true as the main result of this section shows.

\begin{thm}
   \label{thm:ECnotS}
   Let $X$ be either a longpipe or a sealed longpipe. If 
   $\mathsf{EC}(X)$ holds, then $a\mathsf{S}(X)$ does not hold.
\end{thm}

Before proving this theorem, let us show that there are concrete examples of such (sealed) longpipes.
We already encontered one: Theorems \ref{thmDir} and \ref{thm:ECnotS} show 
that the longpipe built with $\clubsuit_C$ in Example \ref{exNy} does not satisfy $\mathsf{S}(Y)$.
But there are examples in {\bf ZFC}.

\begin{example}[{Nyikos, in effect}]
   \label{ex:noret}
   There are sealed longpipes $X=\cup_{\alpha\in\omega_1}X_\alpha$ as in Theorem \ref{thm:ECnotS}
   and Lemma \ref{lemma:sealedret}.
   That is,
   $X$ retracts on each $\wb{X_{\alpha+1}}$, $\mathsf{EC}(X)$ holds
   but $a\mathsf{S}(X)$ does not.
\end{example}
\begin{proof}[Details]
This example is also due to Nyikos and is described in \cite[p. 210]{Nyikos:1992}. 
It is actually a longpipe, but we may sew a disc in $X_1$ to obtain a sealed one.
As in Example \ref{ex:ECnotomega_1cpct} to which we refer for definitions,
we first consider a tangent bundle $T\LL_+$ of $\LL_+$ given by some smoothing.
Write $\pi:T^+\to\LL_+$ for the bundle projection.
As written in Example \ref{ex:ECnotomega_1cpct},
in principle, the details of how the smoothing is built are
important, since $T\LL_+$ and $T^+$ may have quite different topological properties depending 
on the particular construction. This is however irrelevant in our case, as we only use the fact that
maps $T^+\to\R$ are constant on the fibers above some some club $C\subset\LL_+$, which
is true for any smoothing \cite[Corollary 4.15]{Nyikos:1992}. 
Now consider the $\Z$-action on the fibers $\langle x,i\rangle\mapsto 2^i\cdot x$ (see \cite[p. 210]{Nyikos:1992}).
Quotienting by this action, we obtain a longpipe $X$ (each fiber is now a circle), and 
the quotient map $q:T^+\to X$ is actually a covering. There is thus a unique $\wt{\pi}:X\to\LL_+$ such that the left
part of the diagram below commutes.
\begin{center}
  \begin{tikzcd}
  \LL_+
    \arrow[r, "\pi", leftarrow]
   & 
   T^+
     \arrow[d, "q" , rightarrow]
     \arrow[r, "f", rightarrow]
   &
   \R
     \arrow[dl, "\wt{f}" , leftarrow]
   \\
   & X 
     \arrow[ul, "\wt{\pi}", rightarrow]
   \end{tikzcd}
\end{center}
Notice that there is no unbounded map $\LL_+\to X$ (otherwise it could be lifted to $T^+$,
see e.g. \cite[Theorem 5, Chapter 2, Section 4]{Hatcher}).
Now, given $\wt{f}:X\to\R$, $f= \wt{f}\circ q:T^+\to\R$ is constant on the fibers above a club $C\subset\LL_+$
and hence so is $\wt{f}$. It follows that $X$ satisfies $\mathcal{DO}$. By Theorem \ref{thm1EC} (c) 
$\mathsf{EC}(X)$ holds
(recall that longpipes are countably compact). Hence by Theorem \ref{thm:ECnotS}
$a\mathsf{S}(X)$ does not hold.
\end{proof}

The proof of Theorem \ref{thm:ECnotS} is done by exhibiting a function $f:X\to\R$ such that any $r$
with $f\circ r = f$ has ``wrong'' homotopy properties; that is, we may apply 
Lemma \ref{lemma:trivialhomotopy} (as in Proposition \ref{thmnoret0}).
For this purpose, we will define families of maps such that if $f_0,f_1$ are distinct
member of one family, then $f_0\circ r \not= f_1$ for any $r$ with relevant domain and range.
Our first proof used a construction in three stages. At that time we found
convenient to try to be a bit systematic in our treatment, to use (commutative) diagrams
and a general abstract lemma (with an almost trivial proof)
enabling us to pass from a stage to the next.
We then noticed that a much simpler construction was available in just two stages.
We could have given the direct argument but chosed to keep our first approach with the abstract lemma
since we find that it separates the proof into more transparent steps. So, let us state and prove our abstract lemma.

\begin{lemma}\label{lemma:norgeneral}
Let $Z_0,Z_1,X,Y$ be spaces and $f_0,f_1,\wt{f_0},\wt{f_1},\psi_0,\psi_1,\varphi$ be as in the diagram below, with
$\wt{f_0}=f_0\circ\psi_0$, $\wt{f_1}=f_1\circ\psi_1$ and $\psi_1\circ\varphi=id_X$.
\begin{center}
  \begin{tikzcd}
  Z_0  \arrow[r, rightarrow ,"\psi_0"] 
       \arrow[rd , bend right , rightarrow, "\wt{f_0}"'] 
   & X \arrow[r, rightarrow, yshift= 0.4ex ,"\varphi"] 
       \arrow[d, xshift= -0.7ex, "f_0"'] 
       \arrow[d, xshift=0.7ex, "f_1"]
   & Z_1 \arrow[l, rightarrow, yshift= -0.4ex ,"\psi_1"] 
       \arrow[ld , bend left , rightarrow, "\wt{f_1}"]
   \\
   & Y &
   \end{tikzcd}
\end{center}
\noindent
If there is no $r:X\to X$ such that $f_0\circ r=f_1$, then there is no $\wt{r}:Z_1\to Z_0$ such that
$\wt{f_0}\circ\wt{r} = \wt{f_1}$.
\end{lemma}
\proof
  If $\wt{f_0}\circ\wt{r} = \wt{f_1}$ for some $\wt{r}:Z_1\to Z_0$, set 
  $r = \psi_0\circ\wt{r}\circ\varphi$, then $r:X\to X$ and
  $$ f_0\circ r = f_0\circ\psi_0\circ\wt{r}\circ\varphi = 
     \wt{f_0} \circ\wt{r}\circ\varphi = \wt{f_1} \circ \varphi =
     f_1\circ\psi_1\circ\varphi = f_1,$$
  a contradiction.
\endproof

\begin{figure}[h]
  \begin{center}
    \epsfig{ figure = 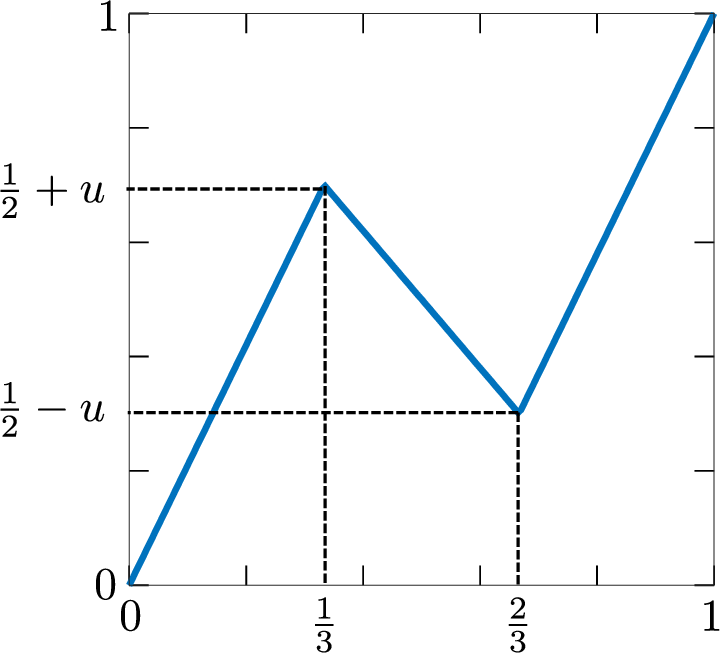, width = .38\textwidth}
    \caption{The family of interval maps $\{f_u\,:\,u\in(0,\frac{1}{2})\}$.}
    \label{fig:f_up_n}
  \end{center}
\end{figure}

Our building block is a simple family of interval maps whose idea was given to us by D. Gauld. 
In the remaining of this section, we denote the closed unit interval $[0,1]$ by $I$.
Let $u\in(0,\frac{1}{2})$. We define $f_u: I \to I $ to be the map 
(depicted in Figure \ref{fig:f_up_n}) that takes values 
$0,\frac{1}{2}+u,\frac{1}{2}-u,1$ at $0,\frac{1}{3},\frac{2}{3},1$,
respectively,
and is linear in between.

\begin{lemma}[D. Gauld]\label{lemma:f_u}
   Let $0<u,v<\frac{1}{2}$ with
   $u\not= v$. Then there is no continuous $r: I \to I $ such that $f_u\circ r = f_v$. 
\end{lemma}
\proof
   Suppose $u>v$. Then $r$ must be increasing up to $1/3$ but $r(1/3)<1/3$, 
   then since $f_v$ must decrease, $r$ is decreasing until $2/3$, 
   then again increasing with $r(1)=1$. But we run into problems when $r(x)$ 
   reaches $1/3$ since then $f_u\circ r$ starts decreasing while $f_v$ is not.
   Suppose now that $u<v$. Then we run into problems as soon as $r(x)$ reaches $1/3$.
\endproof

Let $i: I \to B(1)$, $j:B(1)\to I $ be 
defined as $i(t) = \langle 1-t, 0\rangle$ and $j(x) = 1-|x|$ (where $|x|$ is
the Euclidean norm). Define $\wt{f_u}:B(1)\to I $ as
$\wt{f_u}=f_u\circ j$.

\begin{cor}\label{cor:Bnor}
   $\exists \wt{r}:B(1)\to B(1)$ with $\wt{f_u}\circ\wt{r}=\wt{f_v}$ $\Longleftrightarrow$ $u=v$.
\end{cor}
\proof
   The result follows immediately by the fact that $j\circ i = id_{ I }$,
   the diagram below and Lemmas \ref{lemma:norgeneral} and \ref{lemma:f_u}.
  \begin{center}
  \begin{tikzcd}
  B(1)  \arrow[r, rightarrow ,"j"] 
       \arrow[rd , bend right , rightarrow, "\wt{f_u}"'] 
   & { I } \arrow[r, rightarrow, yshift= 0.4ex ,"i"] 
       \arrow[d, xshift= -0.7ex, "f_u"'] 
       \arrow[d, xshift=  0.7ex, "f_v"]
   & B(1) \arrow[l, rightarrow, yshift= -0.4ex ,"j"] 
       \arrow[ld , bend left , rightarrow, "\wt{f_v}"]
   \\
   & { I } &
   \end{tikzcd}
   \end{center}
\endproof

The family of maps $\wt{f_u}$ enables us to prove the following.

\begin{lemma}\label{lemma:rnothomotopic}
   Let $C=[0,2]\times\mathbb{S}^1$ be the cylinder and $B=B(2)\subset\R^2$.
   Then the following hold.
   \\
   (a) There is $f:C\to I $, constant on $[1,2]\times\mathbb{S}^1$
       such that if $r:C\to C$ satisfies $f\circ r = f$,
       then $r\upharpoonright \{2\}\times\mathbb{S}^1$ is not a constant map.
       \\
   (b) There is $g:B\to I $, constant on $\{x\in B\,:\,|x|\ge 1\}$, such that
       if $s:B\to B$ satisfies $g\circ s = g$,
       then $s$ is not constant on $\partial B = \{x\in B\,:\,|x|=2\}$.
\end{lemma}

\begin{figure}[h]
  \begin{center}
    \epsfig{ figure = 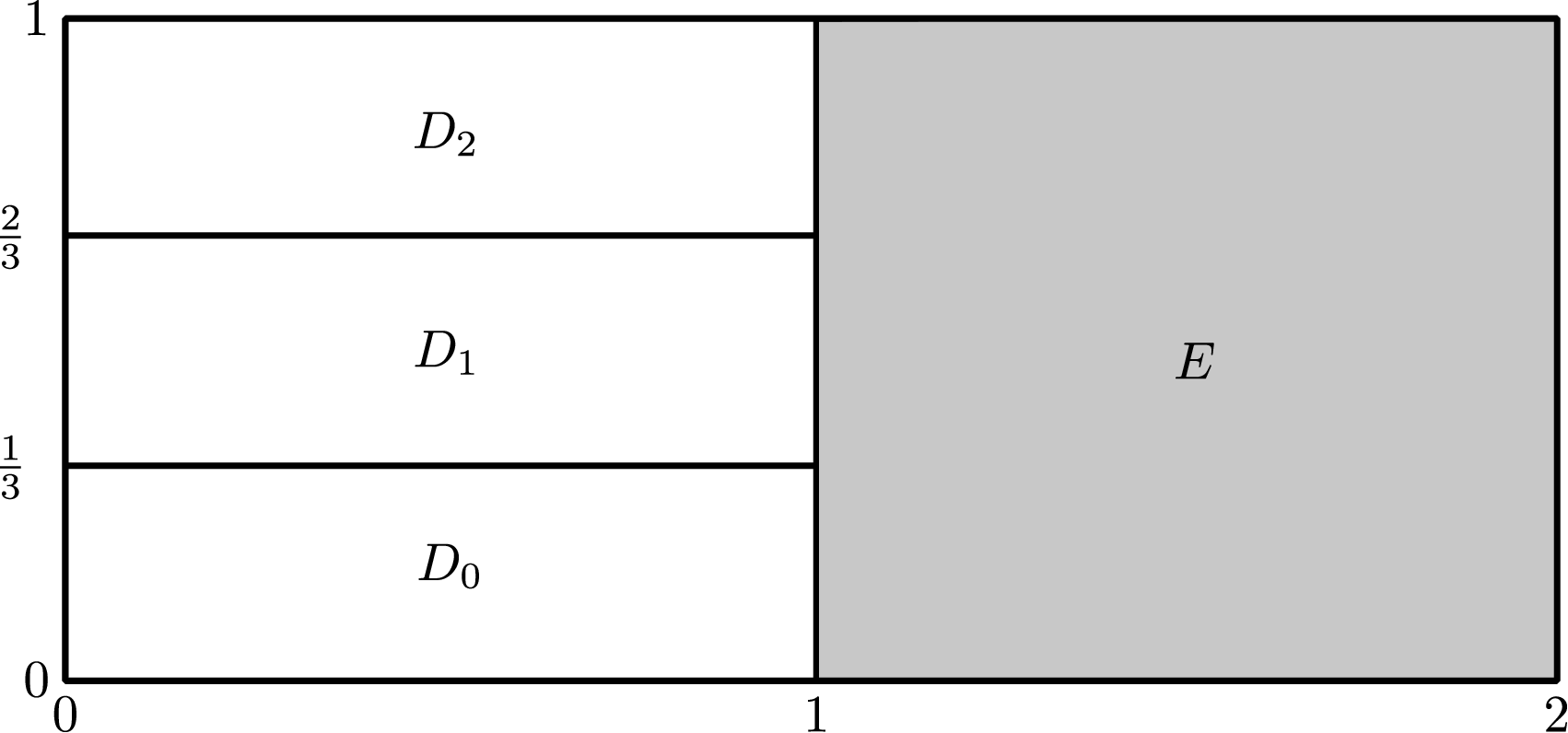, width = .6\textwidth}\quad
    \epsfig{ figure = 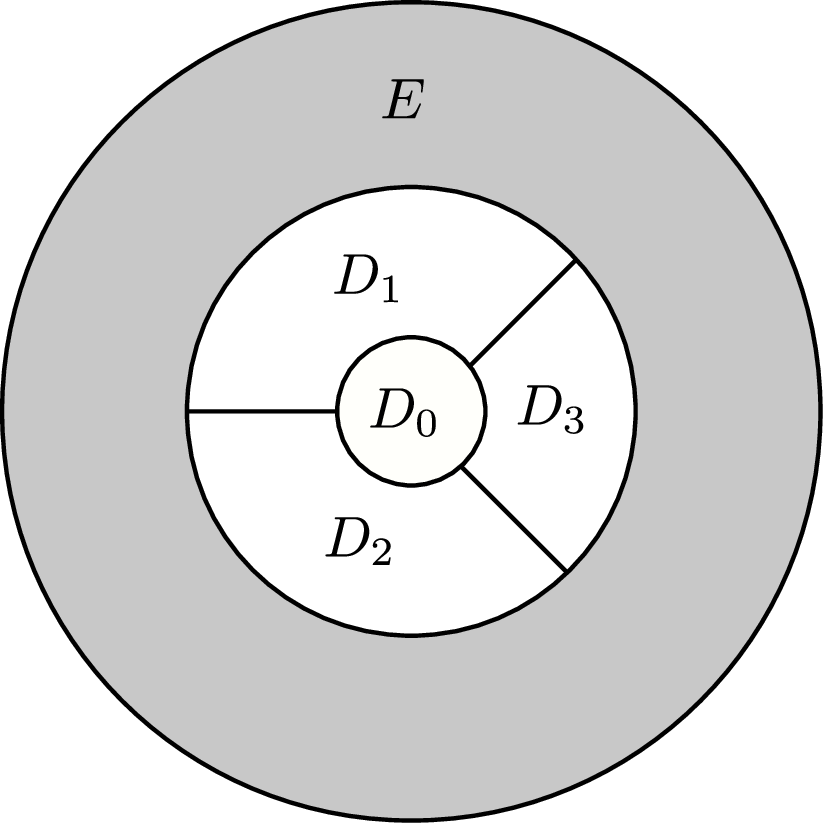, width = .3\textwidth}
    \caption{The subsets for defining the functions $f:C\to I $ and $g:B\to I $.}
    \label{fig:C}
  \end{center}
\end{figure}

\proof
   (a) We see $\mathbb{S}^1$ again as the interval $I$ with endpoints identified.
       Let $E$, $D_0$, $D_1$ and $D_2$ be the closed regions (boundaries are included) of $C$
       depicted in Figure \ref{fig:C} (left).
       Fix any homeomorphisms $\phi_k:D_k\to B(1)$, $k=0,1,2$.
       Fix three distincts points $u_0,u_1,u_2\in(0,\frac{1}{2})$.
       Define $f:C\to I $ to be constant on $0$ on $E$ and equal to $\wt{f}_{u_k}\circ\phi_k$ on $D_k$,
       for $k=0,1,2$.
       By definition, $\wt{f}_{u_k}\circ\phi_k$ takes value $0$ on the topological boundary of $D_k$.
       It follows that $f$ is continuous.
       Let $r:C\to C$ be such that $f\circ r = f$. 
       Since $f^{-1}(\{0\})$ is equal to $E$ union the boundaries of $D_0$, $D_1$ and $D_2$, by connectedness
       $r$ must send the interior of $D_0$ to the interior of $D_0$, $D_1$ or $D_2$.
       Hence, $r\upharpoonright D_0$ has range included in $D_k$ where $k$ is $0$, $1$ or $2$. 
       But we are in the following situation:
   \begin{center}
     \begin{tikzcd}
        D_0  \arrow[r, rightarrow ,"\phi_0"] 
              \arrow[rd , bend right , rightarrow, "\wt{f}_{u_0}\circ\phi_0"'] 
        &  B(1)  \arrow[r, rightarrow, yshift= 0.4ex ,"\phi_k^{-1}"] 
                 \arrow[d, xshift= -0.7ex, "\wt{f}_{u_0}"'] 
                 \arrow[d, xshift=  0.7ex, "\wt{f}_{u_k}"]
        & D_k \arrow[l, rightarrow, yshift= -0.4ex ,"\phi_k"] 
       \arrow[ld , bend left , rightarrow, "\wt{f}_{u_k}\circ\phi_k"]
        \\
        & { I } &
    \end{tikzcd}
   \end{center}
   By Lemma \ref{lemma:norgeneral} and Corollary \ref{cor:Bnor}, we must have $k=0$, and $r$ sends thus 
   $D_0$ into itself. The same argument shows that $r$ sends also 
   $D_1$ and $D_2$ respectively into themselves.
   Let $d_0,d_1:\mathbb{S}^1\to C$ be $d_0(x) = \langle \frac{1}{2},x\rangle$, 
   $d_1(x) = \langle 2,x\rangle$ and 
   $\pi:C\to\mathbb{S}^1$ be the projection on the second coordinate.
   By definition $\pi\circ r\circ d_0$ sends the intervals 
   $[0,\frac{1}{3}],[\frac{1}{3},\frac{2}{3}],[\frac{2}{3},1]\subset\mathbb{S}^1$ respectively into themselves. 
   By Lemma \ref{lemma:trivialhomotopy}, then $\pi\circ r\circ d_1$ cannot be constant, and thus
   $r$ cannot be constant on $\{2\}\times\mathbb{S}^1$.
   \\
   (b) The proof is almost the same, but this time we define $E$, $D_k$ ($k=0,1,2,3$)
       as in the righthandside of Figure \ref{fig:C}.
       Take four distinct points $u_0,u_1,u_2,u_3\in(0,\frac{1}{2})$, 
       fix homeomorphisms $\phi_k:D_k\to B(1)$, and define $g$ exactly as $f$.
       Then, argue as in (a) to show that if $g\circ s = g$, then $s$ 
       sends $D_1\cup D_2\cup D_3$ to itself. By restricting to this subset, we are
       exactly in the same situation than in (a), and we may conclude.
\endproof

We may now prove Theorem \ref{thm:ECnotS}.

\proof[Proof of Theorem \ref{thm:ECnotS}]
   Let $C$ be the cylinder $[0,2]\times\mathbb{S}^1$ and
   $f:C\to\R$ be given by Lemma \ref{lemma:rnothomotopic} (a).
   If $X=\cup_{\alpha<\omega_1}X_\alpha$ is a longpipe, 
   fix a homeomorphism $\Psi:\wb{X_1}\to C$
   and let $\wh{f}$ be defined as $f\circ\Psi$ on $\wb{X_1}$
   and constant on $0$ elsewhere.
   If there is $r:X\to \wb{X_\alpha}$ such that $\wh{f}\circ r = \wh{f}$, then
   since $\wb{X_\alpha}$ embeds in $\R^2$, by $\mathsf{EC}(X,\R^2)$, there is 
   some $\beta$ such that $r$ is constant on $X-X_\beta$.
   We may assume that $\beta>\alpha$.
   Then $r\upharpoonright\wb{X_{\beta+1}}$ is a map sending a cylinder into itself and
   constant on the upper (topological) boundary. But this is impossible
   by Lemma \ref{lemma:rnothomotopic} (a).
   If $X$ is a sealed longpipe, the proof is the same, using $g$ given by Lemma \ref{lemma:rnothomotopic} (b).
\endproof

We end this section with a last example. 
The proofs of the claimed properties are left to the reader, as they are very similar to
what we just did.

\begin{example}\label{ex:wSnotaSnotEC}
   Let $P^0,P^1$ be longpipes such that $\mathsf{EC}(P^k)$ hold for $k=0,1$. Let $P$ be obtained by
   gluing $P^0$ and $P^1$ along their boundary (which is homeomorphic to a circle) in any way.
   Then $w\mathsf{S}(P,\R)$ holds but $a\mathsf{S}(P)$ and $\mathsf{EC}(P)$ do not.
\end{example}


\section{$w\mathsf{S}(M)$, normality, collectionwise Hausdorfness and $\omega_1$-compactness for manifolds}
\label{sec:SscwH}

It happens that, in the world of manifolds, $\mathsf{S}(X)$ and its weakenings $a\mathsf{S}(X)$ 
and $w\mathsf{S}(X)$ (recall Definition \ref{def:wSaS}) 
have an interesting interplay with
normality, collectionwise Hausdorfness and $\omega_1$-compactness. 
Let us first give an example.

\begin{example}(Similar to \cite[Ex. 5.5]{mesziguesnarrow})\\
    Let $X= (\LL_{\ge 0}\times\R) - (\omega_1\times\{0\})$. Then $X$ is a cwH non-normal and non-$\omega_1$-compact surface,
    and $\mathsf{S}(X)$ holds.
\end{example}
Here, of course, $\omega_1$ is seen as a subspace of $\LL_{\ge 0}$.
    All the properties are proved as in \cite[Ex. 5.5]{mesziguesnarrow}
    except $\mathsf{S}(X)$ which is proved exactly as in Example \ref{ex:pasAS}.

\begin{q}\label{prob:normal1}
   If $M$ is a normal manifold such that $\mathsf{S}(M)$, or 
   $a\mathsf{S}(M)$, or $w\mathsf{S}(M)$, hold, is then $M$ $\omega_1$-compact~?
\end{q}

Notice that this question asks for a partial generalization of Theorem \ref{thm:normal} 
(where we have $\mathsf{EC}$ instead of 
$\mathsf{S}$ or its weakenings). The next theorem is a partial answer for the $w\mathsf{S}$-version.

\begin{thm}\label{thmscwh}
   Let $M$ be an $\aleph_1$-scwH manifold. 
   If $w\mathsf{S}(M)$ holds then $M$ is $\omega_1$-compact.
\end{thm}
Recall that a normal $\kappa$-cwH space is $\kappa$-scwH, so Theorem \ref{thmscwh} answers Problem \ref{prob:normal1}
if the following old open problem -- which is probably more fundamental than anything done in the present paper 
and already appeared various times in print -- 
has an affirmative answer:

\begin{prob}\label{prob:normal2}
   Is every normal manifold $\aleph_1$-scwH~?
\end{prob}
As far as we know, there is no consistent counter-example to 
Problem \ref{prob:normal2}, and the answer is affirmative under {\bf V=L}
or in any model obtained after forcing with a Suslin tree, see \cite{Tall:PFAforthemasses}. 

For the proof of Theorem \ref{thmscwh} we will use the fonctions $\wt{f_u}=f_u\circ j$,
where $f_u$ is defined in Figure \ref{fig:f_up_n} and $j$ just before Corollary \ref{cor:Bnor}.
The bulk of the argument for proving Theorem \ref{thmscwh} is done in the next lemma.

\begin{lemma}\label{lemma:discretenoS}
   Let $M$ be a manifold and $D$ be a discrete collection of open sets in $M$. 
   If $D$ is uncountable, then $w\mathsf{S}(M)$ does not hold.
\end{lemma}

\proof 
   Let $B$ be the closed unit ball in $\R^n$, where $n$ is the dimension of $M$.
   We may assume that $D=\{D_\alpha\,:\,\alpha\in\omega_1\}$
   and that each $D_\alpha$ contains some $N_\alpha$ homeomorphic to $B$. 
   Fix homeomorphisms $\phi_\alpha:B\to N_\alpha$.
   Fix a non-continuous $1$-to-$1$ map $\sigma:\omega_1\to(0,\frac{1}{2})$. Define $h:M\to[0,1]$ as follows.
   On the complement of the union of the interiors of the $N_\alpha$, $h$ takes value $0$.
   Then, define 
   $h\upharpoonright N_\alpha $ as $\wt{f}_{\sigma(\alpha)}\circ\phi_\alpha^{-1}$. 
   By construction $h$ is $0$ on the topological boundary of $N_\alpha$; hence, $h$ is continuous
   by discreteness of the $N_\alpha$'s.
   \\
   Suppose that there is some Lindel\"of $Z\subset M$ and $r:M\to Z$ such that 
   $h\circ r\upharpoonright (X-Z) = h\upharpoonright (X-Z)$. Then $Z$ is hereditarily Lindel\"of, hence
   $Z$ intersects at most countably many $N_\alpha$. 
   Fix $\alpha$ such that $Z\cap N_\alpha=\varnothing$.
   Then, by connectedness, the image under
   $r$ of $N_\alpha$ must be contained in some $N_{\beta}\subset Z$, and we have the following diagram.
   \begin{center}
     \begin{tikzcd}
        N_\alpha  \arrow[r, rightarrow ,"\phi_\alpha^{-1}"] 
              \arrow[rd , bend right , rightarrow, "\wt{f}_{\sigma(\alpha)}"'] 
        &  B  \arrow[r, rightarrow, yshift= 0.4ex ,"\phi_\beta"] 
                 \arrow[d, xshift= -0.7ex, "f_{\sigma(\alpha)}"'] 
                 \arrow[d, xshift=  0.7ex, "f_{\sigma(\beta)}"]
        & N_\beta \arrow[l, rightarrow, yshift= -0.4ex ,"\phi_\beta^{-1}"] 
       \arrow[ld , bend left , rightarrow, "\wt{f}_{\sigma(\beta)}"]
        \\
        & { I } &
    \end{tikzcd}
   \end{center}
   This yields a contradiction by Lemmas \ref{lemma:norgeneral} and \ref{lemma:f_u}.
\endproof

\proof[Proof of Theorem \ref{thmscwh}]
   Suppose there is an uncountable closed discrete subset $E$, up to taking a subset, we may assume that $|E|=\aleph_1$.
   Expand $E$ to a discrete collection of closed neighborhoods $\{D_w\,:\,w\in E\}$. 
   Then apply Lemma \ref{lemma:discretenoS}.
\endproof

Note in passing that
since any countable subset
of a manifold is contained in an open set homeomorphic to $\R^n$ (see \cite[Cor. 3.4]{GauldBook}),
an $\omega_1$-compact manifold is cwH.
But it may fail to be scwH (at least consistently),
see Example \ref{ex:lastexample} below, which is moreover non-normal.
This space is eventually constant in $\R$.
Another non-normal manifold $M$ such that $\mathsf{EC}(M)$ holds, but Type I this time,
is Example \ref{ex:ECnotomega_1cpct}.


\section{If `small' means `compact'}\label{sec:cpct}

In this brief section, we look more closely at the properties $\mathsf{P_{cpt}}$.
Let us first gather some trivial facts in a lemma.

\begin{lemma}\label{lemma:cpt-trivialities}
  Let $\mathsf{P}\in\{\mathsf{EC},\mathsf{S},\mathsf{L},\mathsf{BR}\}$.
  Then the following hold. \\
  (a) $ \mathsf{BR_{cpt}} \Longleftarrow \mathsf{EC_{cpt}} \Longrightarrow \mathsf{L_{cpt}} \Longleftarrow \mathsf{S_{cpt}}$.\\
  (b) $\mathsf{P_{cpt}}\Longrightarrow \mathsf{P_{cl}}$.\\
  (c) If $X$ is countably compact and $Y$ any space then
  $\mathsf{P_{cl}}(X,Y)\Longleftrightarrow \mathsf{P_{cpt}}(X,Y)$.\\
  (d) If $\mathsf{L_{cpt}}(X)$ holds then $X$ is pseudocompact.
\end{lemma}

By (d) and Theorem \ref{thm:manifoldsL}, any Type I non pseudocompact manifold 
satisfies $\mathsf{L_{cl}}(X)$ but not $\mathsf{L_{cpt}}(X)$.
Also, 
$\mathsf{S_{cl}}(\LL_+)$ and $\mathsf{BR_{cl}}(\LL_+)$ hold but their cpt-counterparts do not,
for instance.
We note also that the converse of point (d) does not hold.

\begin{example}[{Terasawa, in effect}]\ \\
  There is a pseudocompact Tychonoff $X$ such that $\mathsf{L}(X)$ does not hold.
  \label{ex:Terasawa}
\end{example}
We thank S. Spadaro for showing us this example.
\begin{proof}
  Recall the definition of a $\psi$-space 
  $X=\omega\cup\mathcal{R}$ in Example \ref{ex:notclbis}.
  Notice that the subspace $\mathcal{R}$ is closed discrete and that a subset of $X$ is Lindel\"of
  iff it is countable.
  Terasawa showed in particular in \cite{Terasawa:1980} that there is a $\psi$-space $X$ such that
  $\beta X-X$ is homeomorphic to the Cantor set. 
  Let thus $f:(\beta X-X)\to[0,1]$ have range a Cantor set in $[0,1]$. By Tietze extension theorem, 
  $f$ extends to all of $\beta X$ ($X$ is locally compact hence open in $\beta X$).
  Since $X$ is pseudocompact, $f(X)$ is compact and dense in $f(\beta X)$. It follows that $f(X)=f(\beta X)$.
  Since $f\upharpoonright X$ has uncountable image, no Lindel\"of subspace $W$
  of $X$ satisfies $f(W)=f(X)$.
\end{proof} 
Any $\psi$-space is {\em linearly H-closed}, that is: every chain cover has a dense member
(see \cite{meszigues-od-sel} or \cite{AlasJunqueiraWilson:2019} for a proof), hence $X$ being linearly H-closed 
does not imply $\mathsf{L}(X)$ either. All linearly H-closed spaces are pseudocompact.

\vskip .3cm
In view of points (b)--(d) above, it seems interesting to see whether 
$\mathsf{P_{cl}}(X,Y)\Rightarrow \mathsf{P_{cpt}}(X,Y)$ when $X$ is pseudocompact.
Recall that
a normal pseudocompact space is countably compact 
and a Lindel\"of regular space is normal \cite[Thms 3.10.21 \& 3.8.2]{Engelking}.
A closed subset of a pseudocompact space may fail to be pseudocompact, 
but the following is well known (see e.g. \cite[Exercice 3.10.F(d)]{Engelking}):
\begin{lemma}\label{lemmapseudocompact}
   Let $U$ be an open subset of a Tychonoff pseudocompact space $X$. Then $\wb{U}$ is pseudocompact.
\end{lemma}
This lemma gives almost immediately the following theorems.
We first show that the situation is quite simple for Type I spaces.

\begin{thm}\label{thm:cl-cpt}
   Let $X=\cup_{\alpha\in\omega_1}X_\alpha$ be a regular Type I pseudocompact space and $Y$ be any space. 
   Then $X$ is countably compact, and thus 
   $\mathsf{P}(X,Y)\Longleftrightarrow \mathsf{P_{cpt}}(X,Y)$
   for each $\mathsf{P}\in\{\mathsf{EC},\mathsf{S},\mathsf{L},\mathsf{BR}\}$.
\end{thm}
\proof
   In a regular Type I space $X$, each $\wb{X_\alpha}$ is Lindel\"of regular and hence normal.
   It follows that $X$  is Tychonoff. 
   Hence, each $\wb{X_\alpha}$ is pseudocompact by Lemma \ref{lemmapseudocompact}, thus countably compact,
   thus compact. It follows that
   $X$ is countably compact.
   Conclude with Lemmas \ref{lemma:typeIcl} and \ref{lemma:cpt-trivialities} (c).
\endproof

When $X$ is not of Type I, we can still say something for $\mathsf{EC}$.
\begin{thm}\label{thm:cl-cpt2}
   Let $X$ be a Tychonoff pseudocompact space and $Y$ be any space.
   Then 
   $$\mathsf{EC_{cl}}(X,Y)\Longleftrightarrow \mathsf{EC_{cpt}}(X,Y).$$
\end{thm}
\proof
   The reverse implication is immediate. Suppose that $\mathsf{EC_{cl}}(X,Y)$ holds.
   Given $f:X\to Y$, there is some $y\in Y$ such that $U = f^{-1}(Y-\{y\})$ is contained
   in a closed Lindel\"of subset $Z$ of $X$. Then $Z$ is (regular Lindel\"of hence) normal, and $\wb{U}\subset Z$ as well.
   By Lemma \ref{lemmapseudocompact}, $\wb{U}$ is countably compact and hence compact.
   By definition, $f$ is constant outside of $\wb{U}$. This shows that $\mathsf{EC_{cpt}}(X,Y)$ holds.
\endproof

We note that the implication $\mathsf{EC}(X)\Rightarrow \mathsf{EC_{cpt}}(X)$ for pseudocompact $X$
does not hold, as shown by Example \ref{ex:notclbis}.
If $\mathfrak{b}=\omega_1$, there is even a manifold counter-example.
(This example was also alluded to in Section \ref{sec:interplays}.)
Recall that $\mathfrak{b}$ is the smallest cardinality of an $<^*$-unbounded family of functions $\omega\to\omega$, where
$f <^* g$ iff there is some $n\in\omega$ such that $f(m)<g(m)$ when $m\ge n$.
We can assume that such an unbounded family is well ordered by $<^*$ (see e.g. \cite[Theorem 3.3]{vanDouwen:1984}).
Recall in passing that $\omega_1\le\mathfrak{p}\le\mathfrak{b}\le 2^\omega$ and that each inequality may be strict
depending on the model of set theory.

\begin{example}[{Nyikos}]
   \label{ex:lastexample}
   An $\omega_1$-compact surface $S$ which is not countably compact,
   such that $\mathsf{EC}(X)$
   holds but $\mathsf{EC_{cpt}}(X)$ and $\mathsf{EC_{cl}}(X)$ do not.
   If $\mathfrak{b}=\omega_1$, then $S$ can be made pseudocompact and not scwH.
\end{example}  
\proof[Idea of the construction] We only give a sketch, as
the general construction is detailed in \cite[Example 1.29]{GauldBook}, and appeared for the first
time in \cite[Ex. 6.3]{Nyikos:1990}.
The idea is quite similar to (a version) of Example \ref{ex:notcl}.
Start with
an $<^*$-well ordered and $<^*$ unbounded
family of functions $f_\alpha:\omega\to\omega$, $\alpha\in\mathfrak{b}$. 
We might assume that each $f_\alpha$ is strictly increasing and $f_\alpha(0)=0$.
Now consider (the graphs of) the strictly increasing maps $[0,1)\to[0,1)$ (with supremum $1$) given by 
first embedding $\omega\times\omega$ in $[0,1)^2$, sending $\langle n,m\rangle$ to 
$\langle 1-\frac{1}{n+1},1-\frac{1}{m+1}\rangle$
and interpolating linearly in between to obtain maps $\wh{f}_\alpha:[0,1)\to[0,1)$. 
The surface can then be seen as the unit square $[0,1]^2$ with $\langle 1,1\rangle$ removed, to which is 
attached a copy of $\LL_{\ge 0}$ in such a way that $\lim_{x\to 1}\wh{f}_\alpha(x) = \alpha\in\LL_{\ge 0}$
when $\alpha<\omega_1$. 
(The actual construction by Nyikos is actually slightly different,
but only on a superficial level.)
$\mathsf{EC}(S)$ holds because any real valued map on $\LL_{\ge 0}$ is 
eventually constant and the remainder of the space is Lindel\"of.
The construction is made such that the subset $[0,1)\times\{1\}$ is closed (hence $S$ is not countably compact),
and $\{1\}\times[0,1)$ is ``attached'' at the start of the copy of $\LL_{\ge 0}$.
Since any uncountable subset has a cluster point either in $[0,1]^2-\{\langle 1,1\rangle\}$
or in $\LL_{\ge 0}$, $S$ is $\omega_1$-compact.\\
The real valued map consisting of the projection on the second factor on $[0,1]^2-\{\langle 1,1 \rangle\}$
and constant on $1$ on $\LL_{\ge 0}$ contradicts $\mathsf{EC_{cl}}(X)$. 
\\
If $\mathfrak{b} >\omega_1$, then 
$\lim_{x\to 1} f_\alpha(x)$ does not exist when $\alpha\ge\omega_1$,
but if on the contrary $\mathfrak{b} =\omega_1$,
Nyikos showed that the resulting surface is pseudocompact.  
Indeed, in this case, any countable sequence 
in $[0,1]\times[0,1)$ has a cluster point.
Moreover, the closed discrete subspace $\{ \langle 1-1/n,1\rangle\,:\,n\in\omega\}$ 
cannot be expanded to a discrete open collection, because taking one point in each open set intersected with $[0,1]\times[0,1)$ yields a cluster point.
Hence, $S$ is not scwH.\\
Nyikos original construction does not have a manifold boundary but ours does. To get rid of it, we may 
take two copies and glue their boundaries pointwise.
\endproof

Notice that $\mathsf{L}(S), \mathsf{BR}(S)$ and $w\mathsf{S}(S)$ hold since $\mathsf{EC}(S)$ does. 
We do not know whether $S$ is d-connected, and if 
any of $\mathsf{L_{cl}}(S)$,
$a\mathsf{S}(S)$, $\mathsf{S}(S)$ or $\mathsf{BR_{cl}}(S)$ do hold.

\section{Summary and tables}\label{sec:tables}

In this section, we summarize some of our results in a concise (and somewhat incomplete) form.
First, and just for the pleasure of drawing an overly complicated diagram, the following
implications hold by Lemmas \ref{lemma:arrows}, \ref{lemma:arrows2} and \ref{lemma:cpt-trivialities}
(whose proofs are all immediate). Recall that by $\mathsf{P}_1\Rightarrow\mathsf{P}_2$
we mean ``for all Hausdorff spaces $X,Y$, $\mathsf{P}_1(X,Y)\Rightarrow\mathsf{P}_2(X,Y)$''.

\begin{center}
  \begin{tikzcd}
  \mathsf{S_{cpt}} 
    \arrow[r, rightarrow] 
    \arrow[dd, rightarrow]
  &
  a\mathsf{S_{cpt}} 
    \arrow[r, rightarrow] 
    \arrow[dd, rightarrow]
  &
  w\mathsf{S_{cpt}} 
    \arrow[r, leftarrow] 
    \arrow[d, rightarrow]
    \arrow[dd, bend left = 60, rightarrow]
  &
  \mathsf{EC_{cpt}} 
    \arrow[r, rightarrow] 
    \arrow[dd, rightarrow]
  &
  \mathsf{BR_{cpt}} 
    \arrow[dd,rightarrow]
  \\
  & & 
  \mathsf{L_{cpt}}  
    \arrow[dd, bend right = 60, rightarrow]
  \\
  \mathsf{S_{cl}} 
    \arrow[very thick, color=red,r, rightarrow] 
    \arrow[dd, leftrightarrow]
  &
  a\mathsf{S_{cl}} 
    \arrow[r, rightarrow] 
    \arrow[very thick, color=red, dd, rightarrow]
  &
  w\mathsf{S_{cl}} 
    \arrow[r, leftarrow] 
    \arrow[d, rightarrow]
    \arrow[dd, bend left = 60, rightarrow]
  &
  \mathsf{EC_{cl}} 
    \arrow[r, rightarrow] 
    \arrow[dd, rightarrow]
  &
  \mathsf{BR_{cl}} 
    \arrow[dd,rightarrow]
  \\
  & & 
  \mathsf{L_{cl}}  
    \arrow[dd, bend right = 60, rightarrow]
  \\
  \mathsf{S} 
    \arrow[very thick, color=red,r, rightarrow] 
  &
  a\mathsf{S} 
    \arrow[r, rightarrow] 
  &
  w\mathsf{S} 
    \arrow[r, leftarrow] 
    \arrow[d, rightarrow]
  &
  \mathsf{EC} 
    \arrow[r, rightarrow] 
  &
  \mathsf{BR} 
  \\
  & & 
  \mathsf{L}  
  \\
   \end{tikzcd}
   \end{center}

If the domain space is Type I, the vertical arrows between $\mathsf{P}$ and $\mathsf{P_{cl}}$ reverse
(Lemma \ref{lemma:typeIcl}) and if the domain space is Type I and countably compact, the vertical arrows 
between $\mathsf{P_{cl}}$ and $\mathsf{P_{cpt}}$ reverse
(Theorem \ref{thm:cl-cpt}).
The thick red arrows are the only ones for which we do not know a counterexample to their converse,
as shown by the tables below.
These red arrows yield the following questions
that an attentive reader has probably already formulated in their head, and
maybe even solved, which is not our case.

\begin{q}
  \label{prob:SwS}
  Are there spaces $X,Y$ such that $a\mathsf{S}(X,Y)$ holds but $\mathsf{S}(X,Y)$ does not~? 
  Is there an example with $Y=\R$ and/or $X$ a manifold~?
\end{q}

\begin{q}
  Are there spaces $X,Y$ such that $a\mathsf{S}(X,Y)$ 
  holds but $a\mathsf{S_{cl}}(X,Y)$ does not~? 
  Is there an example with $Y=\R$ and/or $X$ a manifold~?
\end{q}

The tables below compile most of our (counter-)examples (and some positive results).
We believe that our choice of notation is self explanatory.
Table \ref{table:1} summarizes quickly some of the results of section \ref{sec:isocompact}.
In the same vein, the contents of section \ref{sec:Suslin} and some of section \ref{sec:BRL} are summarized 
in the Tables \ref{table:2}--\ref{table:3}.
Finally, Tables \ref{table:4}--\ref{table:last} scan through the whole paper for
(counter)-examples to $\mathsf{P}(X)$. We denote the discrete space of cardinality $\aleph_1$
by $D_{\aleph_1}$. The question mark ``?'' means that we do not know whether the 
property in question holds (it might however not be difficult to find out).
The fact that $a\mathsf{S}(X)$ does not hold for Examples \ref{ex:notclbis}--\ref{ex:notcl}
(Table \ref{table:4})
follows easily with the same functions we used to show that $\mathsf{L_{cl}}(X)$ does not hold. 
The fact that $w\mathsf{S}(X)$ does not hold for Example \ref{ex:prufer} follows from Lemma \ref{lemma:discretenoS}
(Table \ref{table:4}).
In Table \ref{table:last}, Lemma \ref{lemma:cpt-trivialities} (d) is responsible for the cpt-properties not holding, 
since the spaces in question are not pseudocompact.

\vskip .3cm
\begin{table}[H]
\begin{center}
  \begin{tabular}{|c|c|c|c|c|c|c|}
     \hline
    \multirow{3}{*}{Reference} & 
    \multirow{3}{*}{Axiom} & \multicolumn{5}{|c|}{Property of $Y$} \\
                         \cline{3-7}
                         & & \multirow{2}{*}{isocompact}
                           & \multirow{2}{*}{$G_\delta$ points}
                           & \multirow{2}{*}{countably tight} & loc. compact & \multirow{2}{*}{$\mathsf{EC}(\omega_1,Y)$} \\
                         & & & & & or Type I & \\
    \hline
    Ex. \ref{ex:quotientomega_1} &  & \checkmark & $\times$ & \checkmark  & \checkmark & $\times$\\
    Ex. \ref{exNy} & $\clubsuit_C$ & $\times$ & \checkmark  & \checkmark  & \checkmark & \checkmark \\
    Thm \ref{thm:iso_PFA} & {\bf PFA} & implied &  & assumed & assumed & assumed\\
    \hline
  \end{tabular}
  \
  \vskip .3cm
  \caption{Section \ref{sec:isocompact} -- about $\mathsf{EC}(\omega_1,Y)$} 
  \label{table:1}
\end{center}
\end{table}
\vskip .3cm

\vskip .3cm
\begin{table}[H]
\begin{center}
  \begin{tabular}{|c|c|c|c|c|}
    \hline 
    Reference & Property of $Y$ & Property of $T$ & $\mathsf{S}(T,Y)$ & $\mathsf{BR}(T,Y)$ \\
    \hline
    Thm \ref{thm:Steprans} \& Prop. \ref{prop:stepransBR} & submetrizable & 
                           $\omega_1$-cpct & \checkmark & \checkmark \\
    Prop. \ref{prop:stepransBR} & uncountable & & assumed & implied \\
    Lemma \ref{ex:BRnotS} & $\Q$ & c-special & $\times$ & \checkmark \\
    Ex. \ref{ex:suslinline} & her. separable & Suslin & $\times$ & $\times$  \\
    \hline
  \end{tabular}
  \
  \vskip .3cm
  \caption{Section \ref{sec:Suslin} -- $\mathsf{S}(T,Y)$ and $\mathsf{BR}(T,Y)$ when $T$ is an $\omega_1$-tree} 
  \label{table:2}
\end{center}
\end{table}
\vskip .3cm

\begin{table}[H]
\begin{center}
  \begin{tabular}{|c|c|c|c|}
    \hline 
    Reference & $X$ is Type I  & $\mathsf{L_{cl}}(X)$ & $\mathsf{BR_{cl}}(X)$ \\
    \hline
    Thm \ref{thm:ccQ} (b) & assumed & implied & implied \\
    Ex. \ref{ex:ccnotQ} & $\times$ & $\times$ & $\times$ \\
    \hline
  \end{tabular}
  \
  \vskip .3cm
  \caption{Section \ref{sec:BRL} -- $\mathsf{BR_{cl}}(X)$ and $\mathsf{L_{cl}}(X)$ when $X$ is countably compact}
  \label{table:3}
\end{center}
\end{table}
\vskip .3cm

\vskip .3cm
\begin{table}[H]
\begin{center}
  \begin{tabular}{|c|c|c|c|c|c|c|c|c|c|c|}
    \hline
    \multirow{2}{*}{$X$} & 
                                 \multicolumn{9}{|c|}{Property $\mathsf{P}$} & $X$ is a  \\ 
                                 \cline{2-10}
                                 & $\mathsf{L}$ & $\mathsf{L_{cl}}$ & $\mathsf{BR}$ & $\mathsf{BR_{cl}}$ 
                                 & $\mathsf{EC}$  & $\mathsf{EC_{cl}}$
                                 & $w\mathsf{S}$ & $a\mathsf{S}$ & $\mathsf{S}$ & manifold \\
    \hline
    $D_{\aleph_1}$ &  $\times$ & $\times$ & $\times$ & $\times$ & $\times$ & $\times$ & $\times$ & $\times$ & $\times$ 
                        & $\times$\\
    $\LL_+$  
                 & \checkmark & \checkmark  & \checkmark & \checkmark & \checkmark & \checkmark & \checkmark
                 & \checkmark & \checkmark & \checkmark \\ 
    $\LL$ &  \checkmark & \checkmark  & \checkmark & \checkmark & $\times$ & $\times$ & \checkmark & \checkmark & \checkmark
          & \checkmark \\ 
    Ex. \ref{ex:notclbis} 
                  &  \checkmark & $\times$ & \checkmark & $\times$ & \checkmark & $\times$               
                  & \checkmark & $\times$ & $\times$ & $\times$
                  \\
    Ex. \ref{ex:notcl} 
                  & \checkmark & $\times$ & \checkmark & $\times$ & \checkmark & $\times$               
                  & \checkmark & $\times$ & $\times$ & $\times$
                  \\
    Ex. \ref{ex:pasAS} (a) 
                  & \checkmark & \checkmark & $\times$ & $\times$ & $\times$ & $\times$
                  & \checkmark & \checkmark & \checkmark & $\times$ \\
    Ex. \ref{ex:pasAS} (b) 
                  & $\times$ & $\times$ & \checkmark & \checkmark & $\times$ & $\times$
                  & $\times$ & $\times$ & $\times$ & $\times$ \\
    Ex. \ref{ex:ECnotomega_1cpct} 
                  & \checkmark & \checkmark & \checkmark & \checkmark & \checkmark & \checkmark 
                  & \checkmark & \checkmark & \checkmark & \checkmark \\                  
    Ex. \ref{ex:prufer} 
                  & \checkmark & \checkmark & $\times$ & $\times$ & $\times$ & $\times$ & $\times$
                  & $\times$ & $\times$ & \checkmark \\
    Ex. \ref{ex:connectednoL} 
                  & $\times$ & $\times$ & $\times$ & $\times$ & $\times$ & $\times$ & $\times$
                  & $\times$ & $\times$ & $\times$ \\              
    Ex. \ref{ex:QRR2} &  \checkmark & \checkmark & \checkmark & \checkmark 
                  & $\times$ & $\times$ & $\times$ & $\times$ & $\times$ & \checkmark \\
    Ex. \ref{ex:noret} 
                  & \checkmark & \checkmark & \checkmark & \checkmark & \checkmark & \checkmark
                           & \checkmark & $\times$ & $\times $ & \checkmark\\
    Ex. \ref{ex:wSnotaSnotEC} 
                  & \checkmark & \checkmark & \checkmark & \checkmark
                          & $\times$ & $\times$ & \checkmark & $\times$ & $\times$ & \checkmark \\
    Ex. \ref{ex:lastexample} 
                  & \checkmark & ? & \checkmark & ? & \checkmark & $\times$               
                  & \checkmark & ? & ? & \checkmark
                  \\    
    \hline
  \end{tabular}
  \
  \vskip .3cm
  \caption{$\mathsf{P}(X)$ and $\mathsf{P_{cl}}(X)$}
  \label{table:4}
\end{center}
\end{table}

\vskip .3cm
\begin{table}[H]
\begin{center}
  \begin{tabular}{|c|c|c|c|c|c|c|c|c|c|c|c|c|}
    \hline
    \multirow{2}{*}{$X$} & 
                                 \multicolumn{12}{|c|}{Property $\mathsf{P}$} \\ 
                                 \cline{2-13}
                                 &  $\mathsf{L_{cl}}$ & $\mathsf{L_{cpt}}$ & $\mathsf{BR_{cl}}$ & $\mathsf{BR_{cpt}}$ 
                                 & $\mathsf{EC_{cl}}$  & $\mathsf{EC_{cpt}}$
                                 & $w\mathsf{S_{cl}}$ & $w\mathsf{S_{cpt}}$ 
                                 & $a\mathsf{S_{cl}}$ & $a\mathsf{S_{cpt}}$
                                 & $\mathsf{S_{cl}}$ & $\mathsf{S_{cpt}}$ \\
    \hline
    $\LL_{\ge 0}$  
            & \checkmark & \checkmark & \checkmark & \checkmark & \checkmark & \checkmark 
            & \checkmark & \checkmark & \checkmark & \checkmark & \checkmark & \checkmark \\ 
    $\LL_+$  
            & \checkmark & $\times$  & \checkmark & $\times$ & \checkmark & $\times$ 
            & \checkmark & $\times$  & \checkmark & $\times$ & \checkmark & $\times$ \\ 
    $\LL_+\times\R$ 
            & \checkmark & $\times$  & \checkmark & $\times$ & $\times$ & $\times$ 
            & \checkmark & $\times$ & \checkmark & $\times$ & \checkmark & $\times$               
                  \\
    Ex. \ref{ex:noret} 
            & \checkmark & \checkmark  & \checkmark & \checkmark & \checkmark & \checkmark 
            & \checkmark & \checkmark   & $\times$ & $\times$ & $\times$ & $\times$ 
            \\ 
    Ex. \ref{ex:noret}$-\{pt\}$               
            & \checkmark & $\times$ & \checkmark & $\times$ & \checkmark & $\times$               
            & \checkmark & $\times$ & $\times$ & $\times$ & $\times$ & $\times$ 
            \\
    Ex. \ref{ex:wSnotaSnotEC}       
            & \checkmark & \checkmark & \checkmark & \checkmark & $\times$ & $\times$               
            & \checkmark & $\times$ & $\times$ & $\times$ & $\times$ & $\times$ 
            \\
    Ex. \ref{ex:wSnotaSnotEC} $-\{pt\}$       
            & \checkmark & $\times$ & \checkmark & $\times$ & $\times$ & $\times$               
            & \checkmark & $\times$ & $\times$ & $\times$ & $\times$ & $\times$ 
            \\
    \hline
  \end{tabular}
  \
  \vskip .3cm
  \caption{$\mathsf{P_{cl}}(X)$ and $\mathsf{P_{cpt}}(X)$ for manifolds (with boundary for $\LL_{\ge 0}$)}
  \label{table:last}
\end{center}
\end{table}

\section{Questions}
\label{sec:questions}

For the convenience of the reader, we re-state here the questions that appeared in the text, and add a few new ones.

\begin{q}
    Is there a manifold $M$ such that $\mathsf{S}(M)$ holds
    but not $\mathsf{BR}(M)$~? 
\end{q}

\begin{q}
   Is there an example in {\bf ZFC} of a countably compact space $X$ such that
   $\mathsf{EC}(X)$ holds but $\mathsf{EC_{cl}}(X)$ does not~?
\end{q}

\begin{q}
    Is there is a countably compact space $X$ such that $\mathsf{L}(X)$ or $\mathsf{BR}(X)$ does not hold~?
\end{q}

\begin{q}
  Are there spaces $X,Y_1,Y_2$ such that
  $\mathsf{S}(X,Y_i)$ holds for $i=1,2$ and $\mathsf{S}(X,Y_1\times Y_2)$ does not~?
\end{q}

\begin{q}
   If $M$ is a normal manifold such that $\mathsf{S}(M)$, or 
   $a\mathsf{S}(M)$, or $w\mathsf{S}(M)$, hold, is then $M$ $\omega_1$-compact~?
\end{q}

\begin{q}
  Are there spaces $X,Y$ such that $a\mathsf{S}(X,Y)$ holds but $\mathsf{S}(X,Y)$ does not~? 
  Is there an example with $Y=\R$ and/or $X$ a manifold~?
\end{q}

\begin{q}
  Are there spaces $X,Y$ such that $a\mathsf{S}(X,Y)$ 
  holds but $a\mathsf{S_{cl}}(X,Y)$ does not~? 
  Is there an example with $Y=\R$ and/or $X$ a manifold~?
\end{q}

We now gather some questions about relationships between our properties and other 
recent generalizations
of Lindel\"ofness and/or compactness. Some of them might have easy answers.

\begin{q}
  If $X$ is linearly Lindel\"of, does any (all~?) properties $\mathsf{P}(X),\mathsf{P_{cl}}(X)$ hold~?
\end{q}

Recall that a space is {\em cellular Lindel\"of} [resp. {\em strongly cellular Lindel\"of}\;] (resp. 
{\em cellular compact})
iff given any collection of disjoint open subsets, there is a Linde\"of [resp. closed Lindel\"of] (resp. compact)
subset that intersects each member of the collection.
Cellular Lindel\"of spaces were introduced by S. Spadaro and A. Bella in \cite{BellaSpadaro:2018}
and cellular compact ones by Tkachuk and Wilson in 
\cite{TkachukWilson:2019}.
Of course, cellular compact $\Rightarrow$ strongly cellular Lindel\"of $\Rightarrow$ cellular Lindel\"of.
Notice that $\mathsf{EC_{cpt}}(\omega_1)$ holds but $\omega_1$ is not cellular Lindel\"of.
Separable spaces, in particular $\psi$-spaces, are cellular Lindel\"of, hence Example \ref{ex:Terasawa}
shows that $X$ cellular compact does not imply $\mathsf{L}(X)$.
Cellular compact spaces are linearly H-closed
\cite[Proposition 3.6]{TkachukWilson:2019}, but 
$\psi$-spaces are not cellular compact.

\begin{q}
  If $X$ is cellular compact, does $\mathsf{L}(X)$ hold~? Does $\mathsf{L_{cpt}}(X)$ hold~?
\end{q}

\begin{q}
  If $X$ is strongly cellular Lindel\"of, does $\mathsf{L}(X)$ hold~? Does $\mathsf{L_{cl}}(X)$ hold~?
\end{q}

Finally, a very vague and general question about generalizations of 
our properties. 
Let $\mathscr{Q}$ be a property similar to Lindel\"ofness, as those described just above.
Define $\mathsf{P}_\mathscr{Q}(X,Y)$ for $\mathsf{P}\in\{\mathsf{EC},\mathsf{L},\mathsf{S},\mathsf{BR}\}$
as in Definition \ref{defgen} but asking $Z$ to have property $\mathscr{Q}$.

\begin{q}
   Are there instances of topological properties $\mathscr{Q}$
   such that $\mathsf{P}_\mathscr{Q}(X,Y)$ are interesting 
   and significantly different from the ones we investigated in this note~?
\end{q}

\bibliographystyle{plain}
\bibliography{../biblio}

\begin{thebibliography}{10}

\bibitem{AlasJunqueiraWilson:2019}
O.T. Alas, L.R. Junqueira, and R.G. Wilson.
\newblock On linearly {H}-closed spaces.
\newblock {\em Topology Appl.}, 258:161--171, 2019.

\bibitem{AlexandroffUrysohn}
P.~Alexandroff and P.~Urysohn.
\newblock {\em {M\'emoire} sur les espaces topologiques compacts}.
\newblock Number~14 in Proceedings of the section of mathematical sciences.
  Koninklijke Nederlandse Akademie van Wetenschappen te Amsterdam, 1929.

\bibitem{Bacon:1970}
P.~Bacon.
\newblock The compactness of countably compact spaces.
\newblock {\em Pac. J. Math.}, 32(3):587--592, 1970.

\bibitem{meszigues-od-sel}
M.~Baillif.
\newblock Notes on linearly {H}-closed spaces and od-selection principles.
\newblock {\em Top. Proc.}, 54:109--124, 2019.

\bibitem{mesziguesnarrow}
M.~Baillif.
\newblock Relative (functionally) {Type I} spaces and narrow subspaces.
\newblock {\em Top. Proc.}, 62:217--258, 2023.

\bibitem{mesziguesContractibility}
M.~Baillif.
\newblock Non-metrizable manifolds and contractibility.
\newblock {\em Appl. Gen. Topol.}, 26(1):303--339, 2025.
\newblock doi: 10.4995/agt.2025.21796.

\bibitem{MesziguesNyikosExample}
M.~Baillif.
\newblock Filling the gaps in an unpublished example of {N}yikos: a countably
  compact non-compact manifold under {$\clubsuit_C$}, 2026.
\newblock ArXiv preprint {\em https://arxiv.org/abs/2607.05535}.

\bibitem{Balogh:1989}
Z.~Balogh.
\newblock On compact {H}ausdorff spaces of countable tightness.
\newblock {\em Proc. Amer. Math. Soc.}, 105:755--764, 1989.

\bibitem{BellaSpadaro:2018}
A.~Bella and S.~Spadaro.
\newblock On the cardinality of almost discretely {Lindel\"of} spaces.
\newblock {\em Monatsh. Math.}, 186:345--353, 2018.

\bibitem{BennettLutzer:1997}
H.~Bennett and D.~Lutzer.
\newblock Diagonal conditions in ordered spaces.
\newblock {\em Fund. Math.}, 153(2):99--123, 1997.

\bibitem{ChoPark}
M.H. Cho and W.W. Park.
\newblock Isocompactness and related topics of weak covering property.
\newblock {\em Bull. Korean Math. Soc.}, 39(2):347--357, 2002.

\bibitem{vanDouwen:1984}
E.~K.~Van Douwen.
\newblock The integers and topology.
\newblock In {\em Handbook of {S}et-{T}heoretic {T}opology (Kenneth Kunen and
  Jerry~E. Vaughan, eds.)}, pages 111--167. North-Holland, Amsterdam, 1984.

\bibitem{Eisworth:2002}
T.~Eisworth.
\newblock On perfect pre-images of $\omega_1$.
\newblock {\em Topology Appl.}, 125(2):263--278, 2002.

\bibitem{Engelking}
R.~Engelking.
\newblock {\em General topology}.
\newblock Heldermann, Berlin, 1989.
\newblock Revised and completed edition.

\bibitem{GauldBook}
D.~Gauld.
\newblock {\em Non-metrisable manifolds}.
\newblock Springer-Verlag, New York-Berlin, 2014.

\bibitem{Hart-Souslin}
K.P. Hart.
\newblock More remarks on {S}ouslin properties and tree topologies.
\newblock {\em Topology Appl.}, 15:151--158, 1983.

\bibitem{Hatcher}
A.~Hatcher.
\newblock {\em Algebraic topology}.
\newblock Cambridge University Press, Cambridge, 2002.
\newblock (A free electronic version is available).

\bibitem{IsmailNyikos}
M.~Ismail and P.~Nyikos.
\newblock On spaces in which countably compact sets are closed, and hereditary
  properties.
\newblock {\em Topology Appl.}, 11(3):281--292, 1980.

\bibitem{Jech:Handbook}
T.~Jech.
\newblock Stationary sets.
\newblock In {\em Handbook of {S}et {T}heory (Foreman, M. and Kanamori, A.
  (eds))}, pages 93--128. Springer, Dordrecht, 2010.

\bibitem{JustWeeseII}
W.~Just and M.~Weese.
\newblock {\em Discovering Modern Set Theory. {II}: {S}et-theoretic, tools for
  every mathematician}, volume~18 of {\em Graduate Studies in Mathematics}.
\newblock American Mathematical Society, Providence, 1997.

\bibitem{Kunen}
K.~Kunen.
\newblock {\em Set theory, an introduction to independance proofs (reprint)}.
\newblock Number 102 in Studies in Logic and the Foundations of Mathematics.
  North-Holland, Amsterdam, 1983.

\bibitem{KunenLarsonSteprans:2012}
K.~Kunen, J.~Larson, and J.~{Stepr\=ans}.
\newblock Continuous maps on {Aronszajn} trees.
\newblock {\em Order}, 29(2):311--316, 2012.

\bibitem{MardaniThesis}
A.~Mardani.
\newblock Topics in the general topology of non-metric manifolds.
\newblock The University of Auckland. Ph.D Thesis, 2014.

\bibitem{Mrowka:1977}
S.~Mr\'owka.
\newblock Some set-theoretic constructions in topology.
\newblock {\em Fund. Math.}, 94:83--92, 1977.

\bibitem{Nyikos:Antidiamond}
P.~Nyikos.
\newblock Antidiamond and anti-{PFA} axioms and topological applications.
\newblock Preliminary Draft on his webpage.

\bibitem{Nyikos:1984}
P.~Nyikos.
\newblock The theory of nonmetrizable manifolds.
\newblock In {\em Handbook of {S}et-{T}heoretic {T}opology (Kenneth Kunen and
  Jerry~E. Vaughan, eds.)}, pages 633--684. North-Holland, Amsterdam, 1984.

\bibitem{Nyikos:countablycompact1986}
P.~Nyikos.
\newblock Progress on countably compact spaces.
\newblock In Z.~Frolik, editor, {\em General {T}opology and its {R}elations to
  {M}odern {A}nalysis and {A}lgebra {VI}}, pages 379--406, Berlin, 1988.
  Heldermann Verlag.

\bibitem{Nyikos:1990}
P.~Nyikos.
\newblock On first countable, countably compact spaces {III}: The problem of
  obtaining separable non-compact examples.
\newblock In {\em Open Problems in Topology}, pages 127--161. North-Holland,
  Amsterdam, 1990.

\bibitem{Nyikos:HereditaryNormality}
P.~Nyikos.
\newblock Hereditary normality versus countable tightness in countably compact
  spaces.
\newblock {\em Topology Appl.}, 44:271--292, 1992.

\bibitem{Nyikos:1992}
P.~Nyikos.
\newblock Various smoothings of the long line and their tangent bundles.
\newblock {\em Adv. Math.}, 93(2):129--213, 1992.

\bibitem{Nyikos:trees}
P.~Nyikos.
\newblock Various topologies on trees.
\newblock In P.R. Misra and M.~Rajagopalan, editors, {\em Proceedings of the
  Tennessee Topology Conference}, pages 167--198. World Scientific, 1997.

\bibitem{Roitman:1984}
J.~Roitman.
\newblock Basic {S} and {L}.
\newblock In {\em Handbook of {S}et-{T}heoretic {T}opology (Kenneth Kunen and
  Jerry~E. Vaughan, eds.)}, pages 295--326. North-Holland, Amsterdam, 1984.

\bibitem{CEIT}
L.~A. Steen and J.~A.~Seebach Jr.
\newblock {\em Counterexamples in topology}.
\newblock Springer Verlag, New York, 1978.

\bibitem{StepransTrees}
J.~{Stepr\=ans}.
\newblock Trees and continuous mapping into the real line.
\newblock {\em Topology Appl.}, 12:181--185, 1981.

\bibitem{Tall:PFAforthemasses}
F.~D. {Tall}.
\newblock {PFA(S)[S] for the masses}.
\newblock {\em Topology Appl.}, 232:13--21, 2017.

\bibitem{Terasawa:1980}
J.~Terasawa.
\newblock Spaces {$N\cup\mathscr{R}$} and their dimensions.
\newblock {\em Topology Appl.}, 11:93--102, 1980.

\bibitem{TkachukWilson:2019}
V.V. Tkachuk and R.G. Wilson.
\newblock Cellular-compact spaces and their applications.
\newblock {\em Acta Math. Hungar.}, 159(2):674--688, 2019.

\bibitem{vanMill:1984}
J.~van Mill.
\newblock An introduction to $\beta\omega$.
\newblock In {\em Handbook of {S}et-{T}heoretic {T}opology (Kenneth Kunen and
  Jerry~E. Vaughan, eds.)}, pages 503--568. North-Holland, Amsterdam, 1984.

\end{thebibliography}

\end{document}